\documentclass[a4paper,12pt]{article}
\usepackage{amsmath,amssymb}
\usepackage[mathscr]{eucal}
\def\qr#1#2{\left(\dfrac{#1}{#2}\right)}
\def\meqv#1#2{\equiv{#1}\ ({\rm mod}\ {#2})}
\newcommand{\Z}{{\bf Z}} \newcommand{\Q}{{\bf Q}} 
\newcommand{\Or}{{\mathcal O}} \newcommand{\G}{{\mathscr{G}}}
\newcommand{\ve}{\varepsilon} \newcommand{\Ze}{{\mathsf{Z}}}
\newcommand{\bpi}{\boldsymbol{\pi}}

\newcommand{\prodp}{\sideset{}{'_{}}\prod}

\renewcommand{\theequation}{\Roman{section}.\arabic{equation}}

\makeatletter
\renewcommand{\section}{\@startsection
{section}{2}{0mm}{0mm}{3mm}{\noindent\bfseries\large}}
\makeatother
\makeatletter
\renewcommand{\subsection}{\@startsection
{subsection}{2}{0mm}{10mm}{3mm}{\noindent\bfseries\normalsize}}
\makeatother

\newtheorem{Theorem}{Theorem}

\newtheorem{Lemma}{Lemma}

\newtheorem{Example}{Example}

\newtheorem{Corollary}{Corollary}

\begin{document}

\begin{center}
{\large\bf Elliptic Gauss Sums and Hecke \boldmath{$L$}-values at 
\boldmath{$s=1$}}

\medskip

{\itshape Dedicated to Professor Tomio Kubota}

\medskip

Tetsuya Asai\footnote{\texttt{apeiros@thn.ne.jp}
\ / Department of Mathematics, Shizuoka University (\textit{retired})}

\end{center}

\bigskip

\noindent
{\bf Introduction}

\medskip

It seems that some retro-fashioned but still fascinating formulas 
lead us to consider the \textit{elliptic Gauss sum}. 
The classical formulas concerned are the following : 
\[ \dfrac{1}{2}\sum_{k=1}^{p-1}\qr{k}{p}\cot\dfrac{\pi k}{p}=h(-p)\sqrt{p}
\quad\ \ \mbox{and} \quad\ \ 
\dfrac{1}{2}\sum_{k=1}^{p-1}\qr{k}{p}\sec\dfrac{2\pi k}{p}=h(-p)\sqrt{p}, \]
where $p\ (>3)$ is a rational prime such that\ \ 
$p\meqv{3}{4}$\ \ and\ \ $p\meqv{1}{4}$, respectively ; 
$(\tfrac{k}{p})$ is the Legendre symbol 
and $h(-p)$ is the class number of the quadratic field 
$\Q(\sqrt{-p})$. 
The formulas are apparently related to the Dirichlet $L$-values at $s=1$. 

To get a typical elliptic Gauss sum, 
we have only to replace the Legendre symbol 
by the cubic or the quartic residue character, 
and the trigonometric function by a suitable elliptic function. 
The notion of elliptic Gauss sum was first introduced by G.\,Eisenstein 
for a concern of higher reciprocity laws, but since then it has been regarded 
seemingly as a minor object of study. (cf. [L], p.311)

We shall here try to reconsider it. 
Especially, we treat the problem of \textit{rationality of the 
coefficient}, so we call, of the elliptic Gauss sum, 
which is an analogy of the coefficient $h(-p)$ in the above classical case. 
A typical example is as follows.

Let ${\rm sl}(u)$ be the lemniscatic sine of Gauss so that 
${\rm sl}((1-i)\varpi\,u)$ is an elliptic function with the period lattice 
$\Z[i]$, where\ 
$\displaystyle{\varpi=2\,\int_0^1\dfrac{dx}{\sqrt{1-x^4}}}$.\ \ 
Let $\boldsymbol{\pi}$ be a primary prime in $\Z[i]$; 
$\boldsymbol{\pi}\meqv{1}{(1+i)^3}$ and assume\ \ 
$p=\boldsymbol{\pi}\,\overline{\boldsymbol{\pi}}\meqv{13}{16}$. 
We consider the sum 
\[ \G_{\boldsymbol{\pi}}
=\dfrac{1}{4}\sum_{\nu\in (\Or/(\boldsymbol{\pi}))^{\times}}
\qr{\nu}{\boldsymbol{\pi}}_{\!\!4}
{\rm sl}((1-i)\varpi\,\nu/\boldsymbol{\pi}). \]

It is not difficult to see that it is expressible as 
$\G_{\pi}=\alpha_{\boldsymbol{\pi}}(\sqrt[4]{-\boldsymbol{\pi}})^3$ by an 
integer $\alpha_{\boldsymbol{\pi}}$ in $\Z[i]$. 
We can, furthermore, find remarkable facts by some experimental observation of 
the exact value\ $\alpha_{\boldsymbol{\pi}}$\ after having chosen 
the \textit{canonical quartic root}\ \ $\widetilde{\boldsymbol{\pi}}
=\sqrt[4]{-\boldsymbol{\pi}}$.\ \ Namely, 

\smallskip

{}\hfill 
$\G_{\pi}=\alpha_{\boldsymbol{\pi}}\,\widetilde{\boldsymbol{\pi}}^3$ 
\textit{with a rational integer} $\alpha_{\boldsymbol{\pi}}$ ; 
\textit{the magnitude of which is rather small.} \hfill {}

\smallskip

\noindent
This is not so trivial, but we can now prove the rationality. 
It proceeds as follows. 

Let $\widetilde{\chi}_{\boldsymbol{\pi}}$ be the Hecke character 
of weight one induced by the quartic residue character to the modulus 
$\boldsymbol{\pi}$. 
As is well known, the associated Hecke $L$-series 
$L(s,\,\widetilde{\chi}_{\boldsymbol{\pi}})$ has the functional equation. 
We have in particular the \textit{central value equation}\ \ 
$L(1,\,\widetilde{\chi}_{\boldsymbol{\pi}})
=C(\widetilde{\chi}_{\boldsymbol{\pi}})\,
\overline{L(1,\,\widetilde{\chi}_{\boldsymbol{\pi}})}$\ \ at\ $s=1$ ; 
the constant $C(\widetilde{\chi}_{\boldsymbol{\pi}})$ is the so-called 
\textit{root number}. Then it will be first shown that the value of 
$L(1,\,\widetilde{\chi}_{\boldsymbol{\pi}})$ is expressed by 
the elliptic Gauss sum $\G_{\boldsymbol{\pi}}$. 
Secondly, it will be seen the root number 
$C(\widetilde{\chi}_{\boldsymbol{\pi}})$ coincides with the 
classical quartic Gauss sum $G_4(\boldsymbol{\pi})$ in this case, 
and fortunately the explicit formula of the value is known 
owing to Cassels-Matthews. 
Finally, the accordant expression of $\widetilde{\boldsymbol{\pi}}$ 
and $G_4(\boldsymbol{\pi})$ combined with the central value equation 
proves immediatly the fact 
$\alpha_{\boldsymbol{\pi}}=\overline{\alpha}_{\boldsymbol{\pi}}$, 
that is, the coefficient $\alpha_{\boldsymbol{\pi}}$ of 
$\G_{\boldsymbol{\pi}}$ is a rational integer. 
As a corollary we shall obtain a new formula on the value 
$L(1,\,\widetilde{\chi}_{\boldsymbol{\pi}})$. It also should be remarked that 
we shall bring out a better understanding of Matthews' formula 
by considering together with the elliptic Gauss sum.

On the other hand, 
Mr.\,Naruo Kanou has observed these coefficients $\alpha_{\boldsymbol{\pi}}$ 
for many primes by computer. According to his result, it holds\ \ 
$-49\le \alpha_{\boldsymbol{\pi}}\le 49$\ \ 
for $35,432$ primes in the interval $13\le p\le 3999949$. 
At present, however, the reason is not completely clarified, 
so we shall not touch on this topic of the small magnitude, 
while some persons say that 
$\alpha_{\boldsymbol{\pi}}^2$ might closely relate to the order of a certain 
Tate-Shafarevich group.

\medskip

In this paper, we describe detail of the proof, 
mainly of rationality of the coefficient. There are two cases : 
the cubic character case and the quartic one. 
Although the idea is common and most of discussion goes in parallel, 
we would like to treat the two cases separately 
in the parts I and II to avoid possible confusion and complication. 
Since the same notations appear with each different meaning, 
we hope the reader would read carefully and also 
he would tolerate some redundant and overlapped description. 

\bigskip

\noindent
{\small 
{\bf Acknowledgements}\ \ 

\smallskip

I am deeply grateful to Prof.\,Yoshikazu Baba, 
who called my attention to a theorem of heptagon, 
which has made a beginning of the study. 
Special thanks are due to my younger friends for their constant support, and 
the task has never been completed without the collaboration with them. 
Especially, Naruo Kanou has first made and provided a big table of 
the coefficients, 
which still appeals and waits for better elucidation; 
Yoshihiro $\widehat{\rm O}$nishi himself has discovered the congruence 
between the coefficients of elliptic Gauss sums and 
the generalized Hurwitz numbers, 
which has brought me a great inspiration and encouragement; 
Yoshinori Mizuno has suggested accurately to use the doubly periodic 
but non-analytic function in computation of Hecke $L$-values, 
and I followed him fully in this point. 
While it was not till a week after the meeting at Kyoto 
that I completed the proof by using the functional equation, 
Mr.\,Seidai Yasuda became independently aware of the same idea just 
after my meeting talk. 
He tenderly informed me of this after my realization. 
I thank him for this and permitting my own publication.}

\bigskip

\normalsize
\setlength{\baselineskip}{14pt}

\noindent
{\bf References}
\begin{list}{}{}
\item[{[IR]}] K.\,Ireland and M.\,Rosen, \textit{A Classical Introduction to 
Modern Number Theory}, Springer, 1982 \vspace{-3mm} 
\item[{[L]}] F.\,Lemmermeyer, \textit{Reciprocity Laws : from Euler to Eisenstein}, Springer, 2000 \vspace{-3mm} 
\item[{[M1]}] C.\,R.\,Matthews, Gauss Sums and Elliptic Functions : 
I. Kummer Sum,\ \ \textit{Inventiones math}. {\bf 52}, 163-185 (1979) 
\vspace{-3mm} 
\item[{[M2]}] C.\,R.\,Matthews, Gauss Sums and Elliptic Functions : 
II. The Quartic Sum,\ \ \textit{Inventiones math}. {\bf 54}, 23-52 (1979) \vspace{-3mm} 
\item[{[W]}] A.\,Weil, 
\textit{Elliptic Functions According to Eisenstein and 
Kronecker}, Springer, 1976
\end{list}

\newpage

\section{The Cubic Character Case}

\bigskip

Let\ $\rho$\ be the cubic root of unity\ $e^{2\pi i/3}$. 
Throughout the part {\bf I}, the field $\Q(\rho)$ and the ring $\Z[\rho]$ 
are abbreviated to $F$ and $\Or$, respectively. 
The set $\Or$ or its constant multiple appears also as a period 
lattice for elliptic functions. 
We use the notations as well : 
$W=\{\pm 1,\,\pm \rho,\,\pm \overline{\rho}\},\ \ 
W'=\{1,\,\rho,\,\overline{\rho}\}$.

\subsection{Special elliptic functions with complex multiplication}

\noindent
{\bf 1.1.}\quad We shall define some special functions which play 
the leading role in our argument. 
We denote by\ $\varpi_1$\ the real period given by 
\[ \varpi_1 \doteq \int_0^1\dfrac{d\,x}{\sqrt[3]{(1-x^3)^2}}
= \sqrt[3]{2}\int_0^1\dfrac{d\,x}{\sqrt{1-x^3}}
=1.76663875 \cdots, \]
and let\ $\wp(u)$\ denote Weierstrass' $\wp$ with the period lattice\  
$\varpi_1\Or$,\ \ so that\ \ $\wp'^2=4\wp^3-27$.\ \ 
Then it is obvious\ $\wp(\varpi_1u)$\ and\ $\wp'(\varpi_1u)$\ 
are elliptic functions of the period\ $\Or$. 
Further, we can get another doubly periodic function 
by a slight modification of Weierstrass' $\zeta$ : 

\medskip

\noindent
{\bf Definition.}\quad The non-analytic but doubly periodic function\ 
$\Ze(u)$\ is defined by
\begin{equation}
 \Ze(u) \doteq 
\zeta(\varpi_1u)-\dfrac{2\pi}{\sqrt{3}\,\varpi_1}\,\overline{u}, \label{cZ}
\end{equation}
where\ $\zeta(u)$ is Weierstrass' $\zeta$ to the period 
$\varpi_1\Or$. Double periodicity of $\Ze(u)$ relative to $\Or$ 
is easily verified by a usual formula of $\zeta$.

\medskip

The addition formula of\ $\Ze$,\ that follows immediately from one 
of $\zeta$,\ is useful :
\begin{equation}
\Ze(u+v)=\Ze(u)+\Ze(v)+\dfrac{1}{2}\,\dfrac{\wp'(\varpi_1u)-
\wp'(\varpi_1v)}{\wp(\varpi_1u)-\wp(\varpi_1v)}. \label{cZadd}
\end{equation}

\medskip

The significance of the following two functions will be clear later 
when we see that they are closely related to some Hecke's $L$-values at $s=1$. 
In fact, they are the corresponding functions with the lemniscatic sine and 
cosine of the quartic case. Anyway we shall find 
that these functions are automatically 
introduced by the associated Hecke $L$-series.

\medskip

\noindent
{\bf Definition.}\quad The elliptic functions\ $\varphi(u)$\ and\ 
$\psi(u)$\ to the period $\Or$ are defined by
\begin{equation} 
\varphi(u) \doteq \dfrac{1}{3}\left\{
\mathsf{Z}\left(u-\dfrac{1}{3}\right)+
\overline{\rho}\,\mathsf{Z}\left(u-\dfrac{\rho}{3}\right)+
\rho\,\mathsf{Z}\left(u-\dfrac{\overline{\rho}}{3}\right) \right\}, \label{cphi}
\end{equation}
\begin{equation}
\psi(u) \doteq -\dfrac{1}{3}\left\{
\mathsf{Z}\left(u-\dfrac{1}{3}\right)+
\rho\,\mathsf{Z}\left(u-\dfrac{\rho}{3}\right)+
\overline{\rho}\,\mathsf{Z}\left(u-\dfrac{\overline{\rho}}{3}\right) \right\}. 
\label{cpsi}
\end{equation}

\medskip

From the addition formula (\ref{cZadd}) 
we can easily derive the following expressions. 
\begin{equation}
\varphi(u)=\dfrac{6\,\wp(\varpi_1u)}{9+\wp'(\varpi_1u)},\quad 
\psi(u)=\dfrac{9-\wp'(\varpi_1u)}{9+\wp'(\varpi_1u)}. \label{ppphipsi}
\end{equation}

We need also the following formula.
\begin{equation}
\varphi(u)^{-1}+\varphi(-u)^{-1}=\dfrac{3}{\wp(\varpi_1u)}=
\dfrac{1}{\sqrt{-3}}
\left\{\mathsf{Z}\left(u-\dfrac{1}{\sqrt{-3}}\right) 
- \mathsf{Z}\left(u+\dfrac{1}{\sqrt{-3}}\right) \right\}. \label{cpeinv}
\end{equation}

\bigskip

It seems that these functions 
are highly basic in the theory of elliptic functions 
relative to the lattice $\Or$, 
especially in the theory of complex multiplication. 
We here list some fundamental properties of these functions, 
which are easily derived from the definitions and 
by usual theory of elliptic functions. 
Sometimes we need further properties, 
including the addition and multiplication formulas, 
which, however, we shall collect in the end of the part (Appendix I.2) 
for descriptive simplicity. 

\[ {\rm div}\,(\varphi)=(0)+((1-\rho)/3)+((\rho-1)/3)-(1/3)
-(\rho/3)-(\overline{\rho}/3), \]
\[ {\rm div}\,(\psi)=(-1/3)+(-\rho/3)+(-\overline{\rho}/3)
-(1/3)-(\rho/3)-(\overline{\rho}/3), \]
\[ \Ze(\rho u)=\overline{\rho}\,\Ze(u),\quad 
\varphi(\rho\,u)=\rho\,\varphi(u),\quad \psi(\rho\,u)=\psi(u), \]
\[ \Ze(-u)=-\Ze(u),\quad \varphi(-u)=-\varphi(u)\,\psi(u)^{-1},
\quad \psi(-u)=\psi(u)^{-1}, \]
\[ \psi(u)=\varphi(-u-1/3),\quad \varphi(u)^{-1}=\varphi(-u+1/3), \]
\[ \varphi'(u)=-3\varpi_1\,\psi(u)^2,\quad 
\psi'(u)=3\varpi_1\,\varphi(u)^2,\quad \varphi(u)^3+\psi(u)^3=1. \]

\bigskip

\bigskip

\noindent
{\bf 1.2.}\quad 
Let $\bpi$ be a complex prime in $\Or=\Z[\rho]$ so that 
$p=\bpi\,\overline{\bpi}\meqv{1}{3}$, 
and we assume also $\bpi\meqv{1}{3}$. 
Then we have\ \ $(\Or/(\bpi))^{\times} 
\cong (\Z/p\Z)^{\times}$. 
We often abbreviate as\ \ $\nu\,({\rm mod}\,\bpi)$\ in such a case 
when $\nu$ runs over $(\Or/(\bpi))^{\times}$.

\medskip

Class field or complex multiplication theory tells us 
that such a division value\ 
$\varphi(1/\bpi)$\ or\ $\psi(1/\bpi)$\ generates 
an abelian extension field of $F=\Q(\rho)$. Namely, we have 

\begin{Lemma} 
Let\ \ $L=F(\varphi(1/\bpi))$\ \ and\ \ 
$L_1=F(\psi(1/\bpi))$. One has
\begin{enumerate}
\item[{\rm (i)}] $L/F$ is a cyclic extension of degree $p-1$, and 
$L_1$ is a subfield ; $[L : L_1]=3$.
\item[{\rm (ii)}] ${\rm Gal}(L/F) \cong 
(\Or/(\bpi))^{\times}$\ 
by corresponding\ $\sigma_{\mu}$ to $\mu$, and it holds\ 
$\varphi(\nu/\bpi)^{\sigma_{\mu}}
=\varphi(\mu\nu/\bpi)$,\ \ 
$\psi(\nu/\bpi)^{\sigma_{\mu}}=\psi(\mu\nu/\bpi)$\ \ 
for arbitrary\ \ $\mu, \nu\in (\Or/(\bpi))^{\times}$.
\item[{\rm (iii)}] $\varphi(1/\bpi),\ \psi(1/\bpi)$ are algebraic integers, 
and particularly\ \ $\psi(1/\bpi)$ is a unit.
\item[{\rm (iv)}] The prime ideal $(\bpi)$ splits completely in $L$ :\ 
$(\bpi)=\mathfrak{P}^{p-1}$ where\ \ $\mathfrak{P}=(\varphi(1/\bpi))$.
\end{enumerate}
\end{Lemma}

In fact, $L$ and $L_1$ are called the ray class fields of the conductors\ 
$(3 \bpi)$ and $(\sqrt{-3}\bpi)$, respectively. 
We omit a proof, but in the appendix we give a brief comment on the 
complex multiplication formula which is crucial for the theory of division 
values. For a general consultation and a background we can refer to [L].
We here only give some numerical examples.

\medskip

\begin{Example}{\normalfont \rmfamily 
\ In each case tabulated below, we have the 
$\bpi$-multiplication formula :
\[ \varphi(\bpi\,u)
=\varphi(u)\cdot \dfrac{U(\varphi(u))}{R(\varphi(u))},\quad  
R(x)=x^{p-1}\,U(x^{-1}).\]
We can also verify 
\[ U(x)=\prodp_{\nu\,({\rm mod}\,\bpi)} (x-\varphi(\nu/\bpi)),\quad 
V(x)^3=\prodp_{\nu\,({\rm mod}\,\bpi)} (x-\psi(\nu/\bpi)), \]
\[ x\,U(x)-R(x)=(x-1)\,V(x)^3. \]
Furthermore it is easy to check 
\[ U(x),\ V(x) \in \Or[x],\quad U(x)\meqv{x^{p-1}}{(\bpi)},\quad 
U(0)=\bpi,\ \ V(0)=1. \]
In fact,\ \ $U(x)$ and $V(x)$ are the minimal polynomials of 
$\varphi(1/\bpi)$ and $\psi(1/\bpi)$, respectively.

\[ \begin{array}{cc|l}
p & \bpi & \quad U(x)\\ \hline
7 & 1+3 \rho & x^6 - \bpi\,x^3 + \bpi \\
13 & 4+3 \rho & x^{12} +(1+3\rho)\,\bpi\,x^9 -3\rho\,\bpi\,x^6-2\bpi\,x^3+\bpi \\
19 & -2+3\rho & 
x^{18}-3\rho\,\bpi\,x^{15}-(3-9\rho)\,\bpi\,x^{12}+(5-12\rho)\,\bpi\,x^9
+6\rho\,\bpi\,x^6-3\bpi\,x^3+\bpi \\[2mm]
p & \bpi & \quad V(x) \\ \hline
7 & 1+3\rho & x^2+\rho\,x+1 \\
13 & 4+3 \rho & x^4-(1+\rho)\,x^3-(1+2\rho)\,x^2-(1+\rho)\,x+1 \\
19 & -2+3 \rho & x^6+(1-\rho)\,x^5+(1+2\rho)\,x^4-x^3+(1+2\rho)\,x^2
+(1-\rho)\,x+1 \\
\end{array} \]
}\end{Example}

\medskip

\noindent
{\bf Remark.}\quad As for the division value $\Ze(1/\bpi)$, we can prove 
the following. (cf. II.1)
\[ L=F(\varphi(1/\bpi))=F(\Ze(1/\bpi)),\quad 
\Ze(\nu/\bpi)^{\sigma_{\mu}}=\Ze(\mu\nu/\bpi)\quad 
(\,\mu, \nu \in (\Or/(\bpi))^{\times}\,). \]

\vspace{5mm}

\subsection{\boldmath{$L$}-series for Hecke characters of weight one}

\noindent
{\bf 2.1.}\quad 
In this section we recall some fundamental facts about the topic, 
which will give a basis and a framework of our whole discussion.

Let $\widetilde{\chi}$ denote a Hecke character of weight one relative to 
the modulus $(\beta)\subset \Or$, so that it is a multiplicative function 
on the ideal group of $\Or$ 
of the following form :
\[ \widetilde{\chi}((\nu))=\chi_1(\nu)\overline{\nu},\ \ 
\chi_1 : (\Or/(\beta))^{\times} \rightarrow {\bf C}^{\times},\ \ 
\chi_1(\ve)=\ve\ \ (\,\ve\in W\,), \]
where $\chi_1$ is an ordinary residue class character to the modulus 
$(\beta)$, and $(\beta)$ is called the coductor 
of\ \ $\widetilde{\chi}$\ \ if\ \ $\chi_1$\ \ is a primitive character 
to the modulus $(\beta)$. 

It is well known the associated $L$-series has the analytic continuation and 
satisfies the functional equation. 

We here follow Weil's argument and his notation.
\begin{eqnarray*}
L(s,\ \widetilde{\chi}) &=& \sum_{\mathfrak{a}} 
\widetilde{\chi}(\mathfrak{a})\,N\mathfrak{a}^{-s}
= \dfrac{1}{6} \sum_{\nu\in \Or} \chi_1(\nu)\,\overline{\nu}\,|\nu|^{-2s} 
\notag \\
&=& \dfrac{1}{6} \sum_{\lambda\,({\rm mod}\,\beta)} \chi_1(\lambda)\,
\sum_{\mu\in\Or} (\overline{\lambda}+\overline{\mu}\overline{\beta})\,
|\lambda + \mu\beta|^{-2s}. 
\end{eqnarray*}
Therefore we have
\begin{equation}
L(s,\ \widetilde{\chi}) = \beta^{-1}\,
|\beta|^{2-2s}\cdot \dfrac{1}{6}
\sum_{\lambda\,({\rm mod}\,\beta)} \chi_1(\lambda)\,
K_1(\lambda/\beta,\ 0,\ s). \label{cL1g}
\end{equation}
Here the function $K_1$ is defined for ${\rm Re}\,s>3/2$ as follows, and 
it is analytically continued to the whole $s$-plane and satisfies 
the own functional equation :  (cf. [W], VIII)
\begin{equation*}
K_1(u, u_0, s) = \sum_{\mu \in \Or}
e^{\tfrac{2\pi}{\sqrt{3}}\,(\overline{u}_0\mu-u_0\overline{\mu})}\,
(\overline{u}+\overline{\mu})\,|u+\mu|^{-2s}, 
\end{equation*}
\begin{equation*}
\left(\dfrac{2\pi}{\sqrt{3}}\right)^{-s}
\varGamma(s)\, K_1(u, u_0, s)
= e^{\tfrac{2\pi}{\sqrt{3}}\,(\overline{u}_0 u-u_0\overline{u})}\, 
\left(\dfrac{2\pi}{\sqrt{3}}\right)^{s-2}
\varGamma(2-s)\, K_1(u_0, u, 2-s).
\end{equation*}

\medskip

If $(\beta)$ is the conductor of $\widetilde{\chi}$,\ \ 
we can apply a usual computation of Gauss sum, and thus we obtain 
the functional equation of Hecke $L$-series in this case : 
\[ \varLambda(s,\ \widetilde{\chi})
=C(\widetilde{\chi})\,\varLambda(2-s,\ \overline{\widetilde{\chi}}), \]
\[ \textit{where}\quad \varLambda(s,\ \widetilde{\chi})
=\left(\dfrac{2\pi}{\sqrt{3\cdot N(\beta)}}\right)^{-s}\,
\varGamma(s)\,L(s,\ \widetilde{\chi}), \]
\[ \textit{and}\quad 
C(\widetilde{\chi})=-\rho\,\beta^{-1}\,\sum_{\lambda\,({\rm mod}\,\beta)}
 \chi_1(\lambda)\,e^{2\pi i\,S(\lambda /\beta)}. \]
In the above we use the following abbreviation :
\[ S(\lambda)=a\ \ \mbox{for}\ \ 
\lambda=a + b\,\overline{\rho}\ \ (a,\ b\in \Q),
\ \ \ \textrm{i.e.}\ \ \ 2\pi i\,S(\lambda) = 
\tfrac{2}{\sqrt{3}}\,\pi\,
(\lambda\,\rho-\overline{\lambda}\,\overline{\rho}). \]

\bigskip

In particular we have a simple equality of Hecke $L$-values at $s=1$, 
which we call \textit{the central value equation}, 
and the constant $C(\widetilde{\chi})$ is called \textit{the root number} :
\begin{Lemma} Let $\widetilde{\chi}$ be a Hecke character 
of weight $1$ with the conductor $(\beta)$. Then
\begin{equation}
L(1,\ \widetilde{\chi})=C(\widetilde{\chi})\,
\overline{L(1,\ \widetilde{\chi})}, \label{ccenteq} 
\end{equation}
\begin{equation}
\mbox{where}\quad 
C(\widetilde{\chi})=-\rho\,\beta^{-1}\,
\sum_{\lambda\,({\rm mod}\,\beta)}
 \chi_1(\lambda)\,e^{2\pi i\,S(\lambda /\beta)}. \label{cthe_root}
\end{equation}
\end{Lemma}
\noindent
{\bf Remark.}\quad $L(1,\ \overline{\widetilde{\chi}})=
\overline{L(1,\ \widetilde{\chi})}$.

\bigskip

\bigskip

\noindent
{\bf 2.2.}\quad On the other hand, we can notice that 
the value $L(1,\ \widetilde{\chi})$ relates to some elliptic functions. 
As is remarked in [W] (VIII,\,\S 14) or in others, the following is valid. 
\[ E_1^*(u) \doteq K_1(u,\ 0,\ 1)=\varpi_1\,\zeta(\varpi_1u)
-\dfrac{2\pi}{\sqrt{3}}\,\overline{u}. \]

By the definition (\ref{cZ}) the right-hand side is nothing but our function\ 
$\varpi_1\,\Ze(u)$. Combining this and the equation (\ref{cL1g}) at\ $s=1$,\ 
we obtain the following formula. 
\begin{Lemma} Under the same condition of the preceding Lemma, it holds 
\begin{equation}
L(1,\ \widetilde{\chi})=\dfrac{\varpi_1}{6\,\beta} 
\sum_{\lambda\,({\rm mod}\,\beta)} \chi_1(\lambda)\,\Ze(\lambda/\beta).
\label{cEGSg}
\end{equation}
\end{Lemma}

\medskip

As we shall later discuss, the sum appeared in the right-hand of (\ref{cEGSg}) 
is a prototype of \textit{elliptic Gauss sum}. 
When a Hecke character $\widetilde{\chi}$ is given in a suitably explicit 
form, we may evaluate both the elliptic Gauss sum and the root number more 
explicitly, and the central values equation will give some relation between 
the two. In particular, from the value of the elliptic Gauss sum, if 
non-vanishing, we can know the value of the root number. 
This is the case of our cubic characters, that is the point of this report.

\bigskip

\begin{Example}{\normalfont \rmfamily 
\ \ The following is probably the simplest case and the derived formula\ \ 
$L(1,\,\widetilde{\chi}_0)=\tfrac{\varpi_1}{3}$\ \ may be compared with 
the classical formula : 
$1-\tfrac{1}{3}+\tfrac{1}{5}-\tfrac{1}{7}+\cdots=\tfrac{\pi}{4}$.

\medskip

The conductor is the ideal $(3)$. The Hecke character\ $\widetilde{\chi}_0$\ 
mod $(3)$\ is given as follows :
\[ \widetilde{\chi}_0((\nu)) \doteq \chi_0(\nu)\,\overline{\nu},\ \ 
\textit{where}\ \ \chi_0\ :\ (\Or/(3))^{\times} \cong W\ \ 
\textit{is the natural isomorphism.} \]
Then we can evaluate 
the $L$-value at $s=1$\ directly by (\ref{cEGSg}) :
\[ L(1,\ \widetilde{\chi}_0)=\dfrac{\varpi_1}{18}
\sum_{\varepsilon\in W} \varepsilon\,\Ze(\varepsilon/3)
=\dfrac{\varpi_1}{3}\,\Ze(1/3)=\dfrac{\varpi_1}{3}, \]
because\ $\Ze(1/3)=1$\ by usual theory of elliptic functions. 
Also we can easily check\ \ $\displaystyle{
C(\widetilde{\chi}_0)=-\dfrac{\rho}{3} \sum_{\varepsilon\in W}\,
\varepsilon\,e^{2\pi i\,S(\varepsilon/3)}=1}$\ \ as is expected.
}\end{Example}

\vspace{5mm}

\subsection{Elliptic Gauss sums for cubic characters}

\noindent
{\bf 3.1.}\quad 
Let $\bpi$ be a primary prime in $\Or$ ; namely\ $\bpi\meqv{1}{3}$. 
Let $\chi_{\bpi}$ be the cubic residue character to the modulus $\bpi$ ; 
the notation will be fixed throughout the part {\bf I} : 
\[ \chi_{\bpi}(\nu)=\qr{\nu}{\bpi}_{\!\!3}\ :\ 
\chi_{\bpi}(\nu)^3=1\ \ \mbox{and}\ \ 
\chi_{\bpi}(\nu)\meqv{\nu^{(p-1)/3}}{\bpi}\ \ (\nu\in (\Or/(\bpi))^{\times}). \]

Let $f(u)$ be a certain elliptic function with the periods $\Or$, which 
we specify below.

\medskip

\noindent
{\bf Definition.}\quad The following is called an \textit{elliptic Gauss sum}.
\begin{equation}
\G_{\bpi}(\chi_{\bpi},\,f) \doteq 
\dfrac{1}{3}\,\sum_{\nu\,({\rm mod}\,\bpi)}\,\chi_{\bpi}(\nu)\,f(\nu/\bpi).
\end{equation}

\medskip

In the part {\bf I}, we deal with the elliptic Gauss sums\ 
$\G_{\bpi}(\chi_{\bpi},\,\varphi),\ \G_{\bpi}(\chi_{\bpi},\,\varphi^{-1})$\ 
and\ $\G_{\bpi}(\chi_{\bpi},\ \psi)$\ only, where\ $\varphi(u)$\ and\ $\psi(u)$ are the special elliptic functions defined by (\ref{cphi}), (\ref{cpsi}). 
So we understand that 
$f(u)$ denotes an arbitrary one of these functions $\varphi(u), 
\varphi(u)^{-1}$ and $\psi(u)$ in the subsequence. 
In these cases, unless \textit{``the parity condition"}, so we call,\ 
$\chi_{\bpi}(\rho\,\nu)\,f(\rho\,\nu/\bpi)=\chi_{\bpi}(\nu)\,f(\nu/\bpi)$\ 
is satisfied, $\G_{\bpi}(\chi_{\bpi},\,f)$ vanishes trivially. 
Since $\varphi(\rho\,u)=\rho\,\varphi(u),\ \psi(\rho\,u)=\psi(u)$\ and\ 
$\chi_{\bpi}(\rho)=\rho^{(p-1)/3}$, we can easily check that 
$\G_{\bpi}(\chi_{\bpi},\,f)$ is not trivial in the following only
three cases. 
The parity condition, however, is not sufficient for non-vanishing of 
the elliptic Gauss sum as we shall see later.

\medskip

\textit{The elliptic Gauss sums which we shall consider are the followings} : 

\medskip

(a)\ \ \ $\G_{\bpi}(\chi_{\bpi},\,\varphi)$\quad \quad \ \textit{in the case}\ \ 
$p=\bpi\,\overline{\bpi}\meqv{7}{18}$,

\medskip

(b)\ \ \ $\G_{\bpi}(\chi_{\bpi},\,\varphi^{-1})$\quad  \ \textit{in the case}\ \ 
$p=\bpi\,\overline{\bpi}\meqv{13}{18}$,

\medskip

(c)\ \ \ $\G_{\bpi}(\chi_{\bpi},\,\psi)$\quad \quad \ \textit{in the case}\ \ 
$p=\bpi\,\overline{\bpi}\meqv{1}{18}$.

\medskip

As noted in Lemma I.1, the division values $\varphi(\nu/\bpi), \psi(\nu/\bpi)$ 
are algebraic integers in $L=F(\varphi(1/\bpi))$, and 
${\rm Gal}(L/F)\cong (\Or/(\bpi))^{\times}$, 
hence we can immediately see the following. 
\begin{equation}
\G_{\bpi}(\chi_{\bpi},\,f)^{\sigma_{\mu}}=\overline{\chi}_{\bpi}(\mu)\,
\G_{\bpi}(\chi_{\bpi},\,f)\quad (\,\mu\in (\Or/(\bpi))^{\times}\,) \label{clagr}
\end{equation}

In particular, $\G_{\bpi}(\chi_{\bpi},\,f)^3$ is an element of $F$, and 
furthermore we have
\begin{Lemma}\ \ 
$\G_{\bpi}(\chi_{\bpi},\,f)^3$\ is an algebraic integer in $\Or$.
\end{Lemma}

\noindent
\textit{Proof.}\ \ We show only the integrity of 
$\G_{\bpi}(\chi_{\bpi},\,\varphi^{-1})^3$. Since $\varphi(1/\bpi)\,
\varphi(\nu/\bpi)^{-1}$ is a unit,\ 
$\mathfrak{P}^3\,(\G_{\bpi}(\chi_{\bpi},\,\varphi^{-1})^3)$ is an integral 
ideal, where $\mathfrak{P}=(\varphi(1/\bpi))$. On the other hand, 
$(\bpi)=\mathfrak{P}^{p-1}$ and $p-1\ge 12$, hence 
$(\G_{\bpi}(\chi_{\bpi},\,\varphi^{-1})^3)$ must be integral by itself.

\begin{Lemma}\ \ $\G_{\bpi}(\chi_{\bpi},\,\varphi)^3\meqv{1}{\sqrt{-3}},\ \ 
\G_{\bpi}(\chi_{\bpi},\,\varphi^{-1})^3\meqv{-1}{\sqrt{-3}}$\ \ and\ \ 
$\G_{\bpi}(\chi_{\bpi},\,\psi)^3\meqv{0}{\sqrt{-3}}$,\ \ if\ \ 
$p=\bpi\,\overline{\bpi}\meqv{7}{18}$,\ \ 
$p=\bpi\,\overline{\bpi}\meqv{13}{18}$\ \ and\ \ 
$p=\bpi\,\overline{\bpi}\meqv{1}{18}$,\ \ respectively. 
\end{Lemma}

\noindent
\textit{Proof.}\quad First, we quote the $\sqrt{-3}$ multiplication formula of 
$\varphi(u)$ (cf. Appendix I.2) :
\[ \varphi(\sqrt{-3}\,u)=\dfrac{\sqrt{-3}\,\varphi(u)\,\psi(u)}
{1+\overline{\rho}\,\varphi(u)^3}. \]
Put $u=\nu/\bpi$, then we know $\mathfrak{P}=(\varphi(\nu/\bpi))
=(\varphi(\sqrt{-3}\,\nu/\bpi))$ and $\psi(\nu/\bpi)$ is a unit. 
Hence an ideal equality\ $(\varphi(\nu/\bpi)^3+\rho)=(\sqrt{-3})$\ holds. 
Therefore we have
\begin{equation}
\varphi(\nu/\bpi)^3\meqv{-1}{\sqrt{-3}}. \label{ccong-1}
\end{equation}

Next, let $S$ be an arbitrary third set 
of\ $(\Or/(\bpi))^{\times}$, 
namely,\ \ $(\Or/(\bpi))^{\times} = S \cup \rho\,S 
\cup \overline{\rho}\,S$. 
By virtue of the parity condition, we have\ \ 
$\displaystyle{\G_{\bpi}(\chi_{\bpi},\,f)
=\sum_{\nu\in S} 
\chi_{\bpi}(\nu)\,f(\nu/\bpi)}$.\ \ So we obtain
\[ \G_{\bpi}(\chi_{\bpi},\,\varphi)^3\equiv 
\sum_{\nu\in S} \varphi(\nu/\bpi)^3
\equiv \dfrac{p-1}{3}\cdot (-1)\equiv 1\!\!\pmod{\sqrt{-3}}. \]

For other cases of\ \ $\G_{\bpi}(\chi_{\bpi},\,\varphi^{-1})$\ \ and\ \ 
$\G_{\bpi}(\chi_{\bpi},\,\psi)$,\ \ 
the same argument holds 
by using 
\[ \bpi\cdot\varphi(\nu/\bpi)^{-3}\meqv{-1}{\sqrt{-3}}
\quad \mbox{and} \quad 
\psi(\nu/\bpi)^3\meqv{-1}{\sqrt{-3}}, \] 
instead of (\ref{ccong-1}), respectively. 
Thus we complete the proof of Lemma I.5.

\bigskip

\bigskip

\noindent
{\bf 3.2.}\quad Obviously from (\ref{clagr}), the value 
$\G_{\bpi}(\chi_{\bpi},\,f)$ belongs to 
the cubic extension over $F$. 
So it is necessary to give an suitable cubic root of $\bpi$ for 
the precise investigation of the value of $\G_{\bpi}(\chi_{\bpi},\,f)$. 

\medskip

Let $S$ be an arbitrary third set of $(\Or/(\bpi))^{\times}$. 
Let $\gamma(S)$  be the cubic root of unity such that 
$\displaystyle{\gamma(S)\equiv -\prod_{\nu\in S} \nu \mod (\bpi)}$. 
 (\,cf. [M1], (1.6)\,)

\medskip

\noindent
{\bf Definition.}\quad The following is called the \textit{canonical cubic 
root} of\ $\bpi$.
\begin{equation}
\widetilde{\bpi} \doteq \gamma(S)^{-1}\,\prod_{\nu\in S}
\varphi(\nu/\bpi).
\end{equation}

\medskip

Because of the property $\varphi(\rho\,u)=\rho\,\varphi(u)$,\ 
$\widetilde{\bpi}$ is independent of the choice of $S$, and also we can 
easily show the following by the theory of complex multiplication. 
(cf. Appendix)
\[ \widetilde{\bpi}^3 = \prodp_{\nu\,({\rm mod}\,\bpi)} \varphi(\nu/\bpi)=\bpi. \]

The following is a fundamental property of the cubic residue symbol :
\begin{equation}
\widetilde{\bpi}^{\sigma_{\mu}} = \chi_{\bpi}(\mu)\,\widetilde{\bpi}, \quad 
(\,\mu\in (\Or/(\bpi))^{\times}\,) \label{cuberes}
\end{equation}
which is also easily verified in view of\ \ 
$\gamma(\mu\,S)=\chi_{\bpi}(\mu)\,\gamma(S)$. 

\bigskip

\noindent
{\bf Definition.}\quad 
The following is called the \textit{coefficient} of the elliptic Gauss sum\ 
$\G_{\bpi}(\chi_{\bpi},\,f)$.
\begin{equation}
\alpha_{\bpi} \doteq \widetilde{\bpi}^{-2}\,
\G_{\bpi}(\chi_{\bpi},\,f).
\end{equation}

\medskip

\begin{Theorem}\ \ 
Let $\widetilde{\bpi}$ be 
the canonical cubic root of $\bpi$. Then the elliptic 
Gauss sum is expressible as follows:
\[ \G_{\bpi}(\chi_{\bpi},\,f)=\alpha_{\bpi}\,\widetilde{\bpi}^2, \]
where the coefficient $\alpha_{\bpi}$ is an algebraic integer in $\Or$. 
Futher, it holds
\begin{equation}
\alpha_{\bpi} \equiv 
\begin{cases} 
1 \!\!\pmod{\sqrt{-3}}\ \ \ 
& \mbox{if}\ \ p=\bpi\,\overline{\bpi}\meqv{7}{18}, \\
-1 \!\!\pmod{\sqrt{-3}}\ \ \ 
& \mbox{if}\ \ p=\bpi\,\overline{\bpi}\meqv{13}{18}, \\
0 \!\!\pmod{\sqrt{-3}}\ \ \ 
& \mbox{if}\ \ p=\bpi\,\overline{\bpi}\meqv{1}{18}.
\end{cases} \end{equation}
\end{Theorem}

\noindent
{\textit{Proof.}}\quad By the definition of the coefficient\ $\alpha_{\bpi}$\ 
and by virtue of the properties (\ref{clagr}) and (\ref{cuberes}),\ \ 
$\alpha_{\bpi}^{\sigma_{\mu}}=\alpha_{\bpi}\ \ 
(\,\mu\in (\Or/(\bpi))^{\times}\,)$, and hence $\alpha_{\bpi}\in F$. 
For the integrity, we can check it similarly to the proof of Lemma I.4 : 
$\alpha_{\bpi}\,\widetilde{\bpi}^2$ is 
an algebraic integer, while $\bpi=\widetilde{\bpi}^3$ is a prime in $\Or$ ; 
it means $\alpha_{\bpi}$ itself is already an integer. 
The latter part of Theorem I.1 is immediately observed by Lemma I.5.

\bigskip

\noindent
{\bf Remark.}\quad When $p\meqv{7 \mbox{ or } 13}{18}$, we can take 
$S=\ker \chi_{\bpi}$ as a typical third set of $(\Or/(\bpi))^{\times}$ ; 
namely, $S$ is the subgroup consisting of all cubic residues mod $\bpi$. 
This choice has some advantages. Particularly, it is valid 
\[ \G_{\bpi}(\chi_{\bpi},\,f)=\sum_{\nu\in S} f(\nu/\bpi),\quad 
\gamma(S)=1, \quad \widetilde{\bpi}=\prod_{\nu\in S} \varphi(\nu/\bpi). \]

\medskip

\begin{Example}{\normalfont \rmfamily 
Consider the case of $\bpi=4+3\rho\ (p=13)$, 
and we shall show $\alpha_{\bpi}=-\overline{\rho}$. 

Take $S=\ker \chi_{\bpi}$. 
Then $S=\{ \pm 1,\,\pm 5\}
=\{1,\,-1,\,1-\overline{\rho},\,-2\overline{\rho}\}$,\ \ and so we have 
\[ \widetilde{\bpi}=\prod_{\nu\in S} \varphi(\nu/\bpi)
=\varphi(1/\bpi)\,\varphi(-1/\bpi)\,\varphi((1-\overline{\rho})/\bpi)\,
\varphi(-2\overline{\rho}/\bpi). \]
By using suitable multiplication formulas (cf. Appendix) we can compute the 
right-hand to the following form : 
\[ \widetilde{\bpi}=\dfrac{(\rho-\overline{\rho})\,\varphi(1/\bpi)^4\,
(\varphi(1/\bpi)^3-1)}{(1+\rho\,\varphi(1/\bpi)^3)(1-2\,\varphi(1/\bpi)^3)}. \]
Therefore $\varphi(1/\bpi)$ is a solution of the following equaition.
\[ (1+2\rho)\,x^7+2\rho\,\widetilde{\bpi}\,x^6-2(1+2\rho)\,x^4+
(2-\rho)\,\widetilde{\bpi}\,x^3-\widetilde{\bpi}=0. \]
The equation is decomposed as follows.
\[ ((1+2\rho)\,x^3+2\,\widetilde{\bpi}\,x^2+\overline{\rho}\,
\widetilde{\bpi}^2\,x-1)
(x^4-2\overline{\rho}\,\widetilde{\bpi}\,x^3
-(1-\rho)\,\widetilde{\bpi}^2\,x^2
+\overline{\rho}\,\widetilde{\bpi}^3\,x
+\widetilde{\bpi})=0. \]
The second factor must be the minimal polynomial of 
$\varphi(1/\bpi)$ over 
$F(\widetilde{\bpi})=\Q(\rho,\,\sqrt[3]{\bpi})$. 
So $\varphi(\nu/\bpi)^{-1}\ \ (\nu\in S)$\ \ are the four 
roots of the reciprocal equation.
\[ x^4+\overline{\rho}\,\widetilde{\bpi}^2\,x^3-(1-\rho)\,\widetilde{\bpi}\,x^2
-2\overline{\rho}\,x+\widetilde{\bpi}^{-1}=0. \]
Comparing the second coefficient, we have
\[ \G_{\bpi}(\chi_{\bpi},\,\varphi^{-1})=\sum_{\nu\in S}\varphi(\nu/\bpi)^{-1}
=-\overline{\rho}\,\widetilde{\bpi}^2
\quad \therefore\ \ \alpha_{\bpi}=-\overline{\rho}, \]
which satisfies certainly $\alpha_{\bpi}\meqv{-1}{\sqrt{-3}}$.

\medskip

In general, it seems pretty hard to compute the value of the coefficient 
$\alpha_{\bpi}$ by hand. 
Some examples by computor are given in the table Appendix I.1.
}\end{Example}

\vspace{5mm}

\subsection{The cubic Hecke characters and 
{\boldmath$L$}-values at {\boldmath$s=1$}}

\noindent
{\bf 4.1.}\quad 
We introduce a Hecke character $\widetilde{\chi}_{\bpi}$ induced by 
the cubic residue character $\chi_{\bpi}$.\ \ 
As mentioned before, it is of the form\ \ $\widetilde{\chi}_{\bpi}((\nu))
=\chi_1(\nu)\,\overline{\nu}$\ \ with a residue class character $\chi_1$. 
For the purpose we first modify the character $\chi_{\bpi}$ 
into $\chi_1$ satisfying $\chi_1(-\rho)=-\rho$.\ \ 
After the preparation of supplementary characters\ 
$\chi_0$ and $\chi_0'$,\ we shall treat the three cases separately 
in view of\ $\chi_{\bpi}(-\rho)=\rho^{(p-1)/3}$.

\medskip

Let $\chi_0$ be the character of $(\Or/(3))^{\times}$ defined by
\[ \chi_0(\nu) \doteq \varepsilon\ \ \mbox{for}\ \ 
\nu\meqv{\varepsilon}{3},\ \ \varepsilon\in 
W=\{\pm 1,\,\pm \rho,\,\pm \overline{\rho}\}. \]
We here should notice that $\chi_0$ gives the natural isomorphism\ 
$(\Or/(3))^{\times} \cong W$.

Let $\chi_0'$ be the character of $(\Or/(\sqrt{-3}))^{\times}$ defined by
\[ \chi_0'(\nu) \doteq \delta\ \ \mbox{for}\ \ 
\nu\meqv{\delta}{\sqrt{-3}},\ \ \delta\in \{\pm 1\}. \]
We also should notice that $\chi_0'$ gives the natural isomorphism\ 
$(\Or/(\sqrt{-3}))^{\times} \cong \{\pm 1\}$.

\medskip

\noindent
{\bf Definition.}\quad 
For each primary prime $\bpi$ in $\Or$, the Hecke character 
$\widetilde{\chi}_{\bpi}$ is defined and fixed throughout the part {\bf I} 
as follows : 
\begin{equation}
\widetilde{\chi}_{\bpi}((\nu)) \doteq \chi_1(\nu)\,\overline{\nu},\quad 
\chi_1 \doteq 
\begin{cases}
\chi_{\bpi}\cdot \overline{\chi}_0 \quad 
& \mbox{for}\quad p=\bpi\,\overline{\bpi}\meqv{7}{18}, \\
\chi_{\bpi}\cdot \chi_0' \quad 
& \mbox{for}\quad p=\bpi\,\overline{\bpi}\meqv{13}{18}, \\
\chi_{\bpi}\cdot \chi_0 \quad 
& \mbox{for}\quad p=\bpi\,\overline{\bpi}\meqv{1}{18}. 
\end{cases}
\end{equation}

\medskip

For later use, we present a list of the circumstance of each case.

\medskip

(a)\ \ \ The case\ \ $p=\bpi\,\overline{\bpi}\meqv{7}{18}$.\ \ 
\begin{equation}
(\Or/(\beta))^{\times} \cong (\Or/(\bpi))^{\times}\,\times\,W\ \ 
\mbox{by}\ \ \lambda \ \ \mbox{to}\ \ (\kappa,\ \varepsilon) :\ 
\lambda\equiv 3\,\kappa + \bpi\,\varepsilon
\!\!\pmod{\beta}. \label{c7}
\end{equation}

\quad The conductor of $\widetilde{\chi}_{\bpi}$ is $(\beta)$\ 
where\ $\beta=3\,\bpi$,\ \ and we have 
$\chi_1(\lambda)=\chi_{\bpi}(3)\,\chi_{\bpi}(\kappa)\,\overline{\varepsilon}$.

\medskip

(b)\ \ \ The case\ \ $p=\bpi\,\overline{\bpi}\meqv{13}{18}$.\ \ 
\begin{equation}
(\Or/(\beta))^{\times} \cong (\Or/(\bpi))^{\times}\,\times\,\{\pm 1\}\ \ 
\mbox{by}\ \ \lambda \ \ \mbox{to}\ \ (\kappa,\ \delta) :\ 
\lambda\equiv \sqrt{-3}\,\kappa + \bpi\,\delta
\!\!\pmod{\beta}. \label{c13}
\end{equation}

\quad The conductor of $\widetilde{\chi}_{\bpi}$ is $(\beta)$\ 
where\ $\beta=\sqrt{-3}\,\bpi$,\ \ and we have 
$\chi_1(\lambda)=\overline{\chi}_{\bpi}(3)\,\chi_{\bpi}(\kappa)\,\delta$.

\medskip

(c)\ \ \ The case\ \ $p=\bpi\,\overline{\bpi}\meqv{1}{18}$.\ \ 
\begin{equation}
(\Or/(\beta))^{\times} \cong (\Or/(\bpi))^{\times}\,\times\,W\ \ 
\mbox{by}\ \ \lambda \ \ \mbox{to}\ \ (\kappa,\ \varepsilon) :\ 
\lambda\equiv 3\,\kappa + \bpi\,\varepsilon
\!\!\pmod{\beta}. \label{c1}
\end{equation}

\quad The conductor of $\widetilde{\chi}_{\bpi}$ is $(\beta)$\ 
where\ $\beta=3\,\bpi$,\ \ and we have 
$\chi_1(\lambda)=\chi_{\bpi}(3)\,\chi_{\bpi}(\kappa)\,\varepsilon$.

\bigskip

\bigskip

\noindent
{\bf 4.2.}\quad 
We are now ready to 
evaluate the value of the associated\ $L$-series, especialy at\ $s=1$,\ 
and we show that $L(1,\,\widetilde{\chi}_{\bpi})$ is expressed 
by the corresponding elliptic Gauss sum. 

\begin{Theorem} Let $\widetilde{\chi}_{\bpi}$ be the Hecke character 
for $\bpi$. Then 
\begin{equation}
\varpi_1^{-1}\,L(1,\,\widetilde{\chi}_{\bpi}) = 
\begin{cases}
-\chi_{\bpi}(3)\,\bpi^{-1}\,
\G_{\bpi}(\chi_{\bpi},\,\varphi)\ \ 
& \mbox{if}\quad p=\bpi\,
\overline{\bpi}\meqv{7}{18}, \\
-\overline{\chi}_{\bpi}(3)\,\bpi^{-1}\,
\G_{\bpi}(\chi_{\bpi},\,\varphi^{-1})\ \ 
& \mbox{if}\quad p=\bpi\,
\overline{\bpi}\meqv{13}{18}, \\
\quad \chi_{\bpi}(3)\,\bpi^{-1}\,
\G_{\bpi}(\chi_{\bpi},\,\psi)\ \ 
& \mbox{if}\quad p=\bpi\,\overline{\bpi}\meqv{1}{18}. 
\end{cases} \end{equation}
\end{Theorem}

\noindent
{\textit{Proof.}}\ \ We follow the formula (\ref{cEGSg}) of Lemma I.3, 
and refer to (\ref{c7}), (\ref{c13}) and (\ref{c1}).

\medskip

(a)\ \ \ The case $p=\bpi\,\overline{\bpi}\meqv{7}{18}$.\ \ 
In view of (\ref{c7}), we have
\begin{eqnarray*}
L(1,\ \widetilde{\chi}_{\bpi}) & = & 
\dfrac{\varpi_1}{6\,\beta}\cdot \sum_{\lambda\,({\rm mod}\,\beta)}\,
\chi_1(\lambda)\,\Ze(\lambda/\beta) \\
&=& \dfrac{\varpi_1}{18\,\bpi}\,\chi_{\bpi}(3)\cdot
\sum_{\kappa\,({\rm mod}\,\bpi)} \chi_{\bpi}(\kappa)
\sum_{\varepsilon\in W} \overline{\ve}\,\Ze(\kappa/\bpi+\varepsilon/3) \\
{} &=& -\dfrac{\chi_{\bpi}(3)\,\varpi_1}{6\cdot \bpi} \cdot \sum_
{\kappa\,({\rm mod}\,\bpi)} \chi_{\bpi}(\kappa)
\left\{\varphi(\kappa/\bpi)+\varphi(-\kappa/\bpi)\right\} \\
&=& -\dfrac{\chi_{\bpi}(3)\,\varpi_1}{\bpi}\cdot \dfrac{1}{3}\,\sum_
{\kappa\,({\rm mod}\,\bpi)} \chi_{\bpi}(\kappa)\,
\varphi(\kappa/\bpi),
\end{eqnarray*}
by the definition (\ref{cphi}).

\medskip

(b)\ \ \ The case $p=\bpi\,\overline{\bpi}\meqv{13}{18}$.\ \ 
In view of (\ref{c13}), we have
\begin{eqnarray*}
L(1,\ \widetilde{\chi}_{\bpi})
&=& \dfrac{\varpi_1}{6\sqrt{-3}\,\bpi}\,\overline{\chi}_{\bpi}(3)\cdot
\sum_{\kappa\,({\rm mod}\,\bpi)} \chi_{\bpi}(\kappa)
\sum_{\delta=\pm 1} \delta\,\Ze(\kappa/\bpi+\delta/\sqrt{-3}) \\
{} &=& -\dfrac{\overline{\chi}_{\bpi}(3)\,\varpi_1}{6\cdot \bpi} \cdot \sum_
{\kappa\,({\rm mod}\,\bpi)} \chi_{\bpi}(\kappa)
\left\{\varphi(\kappa/\bpi)^{-1}+\varphi(-\kappa/\bpi)^{-1}\right\} \\
&=& -\dfrac{\overline{\chi}_{\bpi}(3)\,\varpi_1}{\bpi}\cdot \dfrac{1}{3}\,\sum_
{\kappa\,({\rm mod}\,\bpi)} \chi_{\bpi}(\kappa)\,\varphi(\kappa/\bpi)^{-1},
\end{eqnarray*}
by the formula (\ref{cpeinv}).

\medskip

(c)\ \ \ The case $p=\bpi\,\overline{\bpi}\meqv{1}{18}$.\ \ 
In view of (\ref{c1}), we have
\begin{eqnarray*}
L(1,\ \widetilde{\chi}_{\bpi}) & = & 
\dfrac{\varpi_1}{18\,\bpi}\,\chi_{\bpi}(3)\cdot
\sum_{\kappa\,({\rm mod}\,\bpi)} \chi_{\bpi}(\kappa)
\sum_{\varepsilon\in W} \varepsilon\,\Ze(\kappa/\bpi+\varepsilon/3) \\
{} &=& \dfrac{\chi_{\bpi}(3)\,\varpi_1}{6\cdot \bpi} \cdot \sum_
{\kappa\,({\rm mod}\,\bpi)} \chi_{\bpi}(\kappa)
\left\{\psi(\kappa/\bpi)+\psi(-\kappa/\bpi)\right\} \\
&=& \dfrac{\chi_{\bpi}(3)\,\varpi_1}{\bpi}\cdot \dfrac{1}{3}\,\sum_
{\kappa\,({\rm mod}\,\bpi)} \chi_{\bpi}(\kappa)\,
\psi(\kappa/\bpi),
\end{eqnarray*}
by the definition (\ref{cpsi}). 
Thus the proof of Theorem I.2 is finished. 

\bigskip

It may be noteworthy that the special elliptic functions 
$\varphi(u), \varphi(u)^{-1}$ and $\psi(u)$ appear naturally and 
automatically in these $L$-series ; consequently the associated $L$-series 
would introduce those elliptic functions the division values of 
which generate some abelian extensions of the field $F$.

\vspace{5mm}

\subsection{The explicit formula of the root number\ \ 
{\boldmath$C(\widetilde{\chi}_{\bpi})$}}

\bigskip

\noindent
{\bf 5.1.}\quad 
We require an important formula about 
the classical cubic Gauss sum. 
Let $\bpi$ be a primary prime in $\Or$ ; $\bpi\meqv{1}{3}$, 
and set $p=\bpi\,\overline{\bpi}$ as before. 
The cubic Gauss sum, often called the Kummer sum, 
is defined and is denoted by 
\begin{equation}
G_3(\bpi) 
\doteq \sum_{r\,({\rm mod}\,p)} \chi_{\bpi}(r)\,e^{2\pi ir/p}
= \sum_{r=1}^{p-1} \qr{r}{\bpi}_{\!\!3}\,e^{2\pi ir/p}.
\end{equation}
Also we here should recall our definition of 
the canonical cubic root $\widetilde{\bpi}$ of $\bpi$ : 
\[ \widetilde{\bpi}
=\gamma(S)^{-1}\prod_{\nu\in S}\varphi(\nu/\bpi)\quad 
\mbox{where}\quad \gamma(S)^3=1\ \ \mbox{and}\ \ 
\gamma(S)\meqv{-\prod_{\nu\in S}\nu}{\bpi}, \]
where $S$ is an arbitrary third set of modulus $(\bpi)$ ; 
$(\Or/(\bpi))^{\times} = S \cup \rho\,S \cup \overline{\rho}\,S$. 

\begin{Lemma}[$G_3$-formula]\ \ 
\begin{equation}
G_3(\bpi) = - \chi_{\bpi}(3)\,
\widetilde{\bpi}^2\,\overline{\widetilde{\bpi}}. 
\label{G3}
\end{equation}
\end{Lemma}

\noindent
{\textit{Proof.}}\ \ This is only a slight modification of the celebrated 
Cassels-Matthews formula. 
They use the lattice $\theta\,\Or$ instead of our $\varpi_1\,\Or$, 
where $\theta=\sqrt{3}\,\varpi_1=3.05990807\cdots$.\ \ Let\ \ $\wp_1(u)$ 
denote Weierstrass' $\wp$ with the period lattice\ $\theta\,\Or$. 
Hence the relation $\wp(\varpi_1u)=3\,\wp_1(\theta\,u)$ holds, 
so that $\wp_1'(u)^2=4\,\wp_1(u)^3-1$. Their formula states 

\medskip

\noindent
{\bf Formula} (Cassel-Matthews,\ [M1], Th.1)\ \ 
\begin{equation}
G_3(\bpi)=-\gamma(S)^{-1}\,\bpi\,p^{1/3}\,
\prod_{\nu\in S} \wp_1(\theta\,\nu/\bpi). \label{cMatthews}
\end{equation}

We shall show the formula (\ref{G3}) from (\ref{cMatthews}). 
Indeed it will be seen they are equivalent. 

\noindent
First, by using the following two identities : 
the latter being the $\sqrt{-3}$ multiplication,
\[ \varphi(u)=\dfrac{6\,\wp_1(\theta\,u)}
{\sqrt{3}(\wp_1'(\theta\,u)+\sqrt{3})}\quad \mbox{\textit{and}} \quad 
\wp_1(\sqrt{-3}\,u)=-\dfrac{\wp_1'(u)^2-3}{12\,\wp_1(u)^2}, \]
we have 
\begin{equation*}
\varphi(u)^{-1}\varphi(-u)^{-1} = \wp_1(\sqrt{-3}\,\theta\,u).
\end{equation*}
Next, substitute $u=\nu/\bpi$ 
and make the product over $\nu\in S$, then we have 
\[ \prod_{\nu\in S}\varphi(\nu/\bpi)^{-1}\cdot 
\prod_{\nu\in -S}\varphi(\nu/\bpi)^{-1}
=\prod_{\nu\in \sqrt{-3}\,S}\wp_1(\theta\,\nu/\bpi). \]
Finally, multyply the factor $\gamma(S)^{-1}$ to the both sides and notice 
such properties as below, then we can see that 
Cassels-Matthews' formula easily turns to our formula (\ref{G3}).

In fact, on the one hand\ \ 
$\gamma(S)^{-1}=\gamma(S)^2=\gamma(S)\cdot \gamma(-S)$,\ \ 
and on the other hand\ \ $\gamma(S)^{-1}
=\chi_{\bpi}(\sqrt{-3})\,\gamma(\sqrt{-3}\,S)^{-1}
=\overline{\chi}_{\bpi}(3)\,\gamma(\sqrt{-3}\,S)^{-1}$. 

\bigskip

\noindent
{\bf Remark.}\ \ As is immediately observed by\ $G_3$-formula,\ \ 
$G_3(\bpi)^3=-\bpi^2\,\overline{\bpi}$\ \ holds. 

\bigskip

\bigskip

\noindent
{\bf 5.2.}\quad We are now ready to give the value of the root number 
$C(\widetilde{\chi}_{\bpi})$ explicitly. 

\begin{Theorem}\ \ Let $\widetilde{\chi}_{\bpi}$ be the Hecke character 
for $\bpi$. Then 
\begin{equation}
C(\widetilde{\chi}_{\bpi})=
\begin{cases}
\chi_{\bpi}(3)\,\widetilde{\bpi}^{-1}
\overline{\widetilde{\bpi}}\quad\ \ \ & \mbox{if}\quad 
p=\bpi\,\overline{\bpi}\meqv{7}{18}, \\
\overline{\chi}_{\bpi}(3)\,\widetilde{\bpi}^{-1}
\overline{\widetilde{\bpi}}\quad\ \ \ & \mbox{if}\quad 
p=\bpi\,\overline{\bpi}\meqv{13}{18}, \\
-\chi_{\bpi}(3)\,\widetilde{\bpi}^{-1}
\overline{\widetilde{\bpi}}\quad\ & \mbox{if}\quad 
p=\bpi\,\overline{\bpi}\meqv{1}{18}.
\end{cases} \end{equation}
\end{Theorem}

\noindent
\textit{Proof.}\ \ As a preparation we shall evaluate some simple Gauss sums. 
The first three are easily verified by direct calculation :
\begin{equation} \begin{split}
g(\chi_0)\doteq &
\sum_{\varepsilon\in W} \varepsilon\,e^{2\pi iS(\varepsilon/3)}
=-3\overline{\rho},\quad 
g(\overline{\chi}_0)\doteq 
\sum_{\ve\in W} \overline{\ve}\,e^{2\pi iS(\ve/3)}=3\,\rho \\
 & \mbox{\textit{and}}\quad 
g(\chi_0')\doteq 
\sum_{\delta=\pm 1} \delta\,e^{2\pi iS(\delta/\sqrt{-3})}
=\sqrt{-3}.
\end{split} \end{equation}
The next sum is nothing but the cubic Gauss sum :
\begin{equation}
g(\chi_{\bpi}) \doteq \sum_{\kappa\,({\rm mod}\,\bpi)} 
\chi_{\bpi}(\kappa)\,e^{2\pi iS(\kappa/\bpi)}
=-\overline{\chi}_{\bpi}(\rho)\,\widetilde{\bpi}^2\,
\overline{\widetilde{\bpi}}.  \label{cgchi}
\end{equation}
In fact, we first replace the sum over\ $\kappa\,({\rm mod}\,\bpi)$\ 
by one over\ $r\,({\rm mod}\,p)$,\ and then, by using\ 
$S(r\,\overline{\bpi}/p)=ar/p$\ where $\bpi=a+b\rho\ (a,\,b\in \Z)$, 
we can calculate as follows :
\[ g(\chi_{\bpi})=\sum_{r\,({\rm mod}\,p)} 
\chi_{\bpi}(r)\,e^{2\pi iS(r\,\overline{\bpi}/p)}
=\sum_{r\,({\rm mod}\,p)} \chi_{\bpi}(r)\,e^{2\pi i a r/p} 
= \overline{\chi}_{\bpi}(a)\,G_3(\bpi). \]
Further, since we know $\overline{\chi}_{\bpi}(a)
=\chi_{\bpi}(1-\rho)$\ \ (cf.\ \ [IR] Chap.\,9, Exerc.\,24,\,26.), 
it follows\ $\overline{\chi}_{\bpi}(a)=\overline{\chi}_{\bpi}(\rho)\,\overline{\chi}_{\bpi}(3)$. 
Finally, by applying $G_3$-formula, we can obtain (\ref{cgchi}).

\medskip

We now follow the formula (\ref{cthe_root}) of the root number in Lemma I.2, 
and we treat the three cases separately as in 4.1, especially in view of 
(\ref{c7}), (\ref{c13}) and (\ref{c1}).

\medskip

(a)\ \ The case\ \ $p=\bpi\,\overline{\bpi}\meqv{7}{18}$.\ \ 

\medskip

\noindent
Since\ $\beta=3\,\bpi$,\ 
$\lambda\meqv{3\,\kappa+\bpi\,\varepsilon}{\beta}$\ and\ 
$\chi_1(\lambda)=\chi_{\bpi}(3)\,\chi_{\bpi}(\kappa)\,\overline{\ve}$,\ 
we have 
\begin{eqnarray*}
C(\widetilde{\chi}_{\bpi}) &=& 
-\rho\,\beta^{-1}\,
\sum_{\lambda\,({\rm mod}\,\beta)}
 \chi_1(\lambda)\,e^{2\pi i\,S(\lambda /\beta)} \\
{} &=& -\rho\,(\bpi\cdot 3)^{-1} \chi_{\bpi}(3)\,g(\overline{\chi}_0)\,
g(\chi_{\bpi}) = \chi_{\bpi}(3)\,\widetilde{\bpi}^{-1}
\overline{\widetilde{\bpi}}.
\end{eqnarray*}

\medskip

(b)\ \ The case\ \ $p=\bpi\,\overline{\bpi}\meqv{13}{18}$.\ \ 

\medskip

\noindent
Since\ $\beta=\sqrt{-3}\,\bpi$,\ 
$\lambda\meqv{\sqrt{-3}\,\kappa+\bpi\,\delta}{\beta}$\ and\ 
$\chi_1(\lambda)=\overline{\chi}_{\bpi}(3)\,\chi_{\bpi}(\kappa)\,\delta$,\ \ 
we have 
\[ C(\widetilde{\chi}_{\bpi}) = -\rho\,(\sqrt{-3}\,\bpi)^{-1}\cdot 
\chi_{\bpi}(\sqrt{-3})\cdot g(\chi_0')\,g(\chi_{\bpi})
= \overline{\chi}_{\bpi}(3)\,\widetilde{\bpi}^{-1}\overline{\widetilde{\bpi}}. \]

(c)\ \ The case\ \ $p=\bpi\,\overline{\bpi}\meqv{1}{18}$.\ \ 

\medskip

\noindent
Since\ $\beta=3\,\bpi$,\ 
$\lambda\meqv{3\,\kappa+\bpi\,\varepsilon}{\beta}$\ and\ 
$\chi_1(\lambda)=\chi_{\bpi}(3)\,\chi_{\bpi}(\kappa)\,\varepsilon$,\ 
we have 
\[ C(\widetilde{\chi}_{\bpi}) = 
-\rho\,(\bpi\cdot 3)^{-1}\cdot \chi_{\bpi}(3)\cdot g(\chi_0)\,g(\chi_{\bpi})
= -\chi_{\bpi}(3)\,\widetilde{\bpi}^{-1}\,\overline{\widetilde{\bpi}}. \]
These complete the proof of Theorem I.3.

\vspace{5mm}

\subsection{Rationality of the elliptic Gauss sum coefficient}

\noindent
{\bf 6.1.}\quad 
In Theorem I.1 we have seen that each coefficient $\alpha_{\bpi}$ 
is an algebraic integer in $\Or$. Now we can mention about their 
$\Q$-rationality. More precisely, the coefficient itself is not always 
rational, but it will be seen that 
the essential factor of this is certainly a rational integer. 
The next is our main theorem of the part {\bf I}.
\begin{Theorem}\ 
For a primary prime $\bpi$ in $\Or$ 
there exists a rational integer\ \ $a_{\bpi}$,\ 
and the coefficient\ \ $\alpha_{\bpi}$\ of the 
elliptic Gauss sum is expressed by\ $a_{\bpi}$\ as follows. 
\begin{equation}
\alpha_{\bpi}= \begin{cases}
\chi_{\bpi}(3)\,a_{\bpi}\quad 
\mbox{and}\quad a_{\bpi}\meqv{1}{3} & 
\mbox{if}\quad p=\bpi\,\overline{\bpi}\meqv{7}{18}, \\
\overline{\chi}_{\bpi}(3)\,a_{\bpi}\quad 
\mbox{and}\quad a_{\bpi}\meqv{-1}{3} & 
\mbox{if}\quad p=\bpi\,\overline{\bpi}\meqv{13}{18}, \\
\chi_{\bpi}(3)\,a_{\bpi}\,\sqrt{-3}\quad  
&\mbox{if}\quad p=\bpi\,\overline{\bpi}\meqv{1}{18}. 
\end{cases} \end{equation}
\end{Theorem}

\noindent
\textit{Proof.}\ \ By the theorems I.1,\,I.2 and I.3 
we know already both the explicit values 
of $L(1,\,\widetilde{\chi}_{\bpi})$ and $C(\widetilde{\chi}_{\bpi})$. 
To prove Theorem I.4, we have only to substitute them for the both sides of 
the central value equation (\ref{ccenteq}) of Lemma I.2. 
There are three cases : 

\medskip

(a)\ \ The case\ \ $p=\bpi\,\overline{\bpi}\meqv{7}{18}$.\ \ 
In this case we have
\[ \varpi_1^{-1}\,L(1,\,\widetilde{\chi}_{\bpi})
=-\chi_{\bpi}(3)\,\widetilde{\bpi}^{-1}\,\alpha_{\bpi}\quad \mbox{and}\quad 
C(\widetilde{\chi}_{\bpi})=\chi_{\bpi}(3)\,\widetilde{\bpi}^{-1}\,
\overline{\widetilde{\bpi}}. \]
Hence from the central value equation\ \ 
$L(1,\,\widetilde{\chi}_{\bpi})=C(\widetilde{\chi}_{\bpi})\,
\overline{L(1,\,\widetilde{\chi}_{\bpi})}$,\ \ we can deduce 
\[ -\chi_{\bpi}(3)\,\widetilde{\bpi}^{-1}\,\alpha_{\bpi}
= \chi_{\bpi}(3)\,\widetilde{\bpi}^{-1}\,\overline{\widetilde{\bpi}}\cdot 
(-1)\,\overline{\chi}_{\bpi}(3)\,\overline{\widetilde{\bpi}}^{-1}\,
\overline{\alpha_{\bpi}}\quad \therefore\ \ 
\alpha_{\bpi}=\overline{\chi}_{\bpi}(3)\,\overline{\alpha_{\bpi}}. \]
This means\ \ $\overline{\chi}_{\bpi}(3)\,\alpha_{\bpi}
=\chi_{\bpi}(3)\,\overline{\alpha_{\bpi}} \in \Or \cap {\bf R}$,\ \ 
which we may denote by $a_{\bpi}$,\ \ so that
\[ \alpha_{\bpi}=\chi_{\bpi}(3)\,a_{\bpi}\quad 
\mbox{where}\ \ a_{\bpi}\in \Z\ \ \mbox{and}\ \ a_{\bpi}\meqv{1}{3}. \]
The last congruence follows from Theorem I.1.

\medskip

(b)\ \ The case\ \ $p=\bpi\,\overline{\bpi}\meqv{13}{18}$.\ \ 
Since we have
\[ \varpi_1^{-1}\,L(1,\,\widetilde{\chi}_{\bpi})
=-\overline{\chi}_{\bpi}(3)\,\widetilde{\bpi}^{-1}\,\alpha_{\bpi}\quad \mbox{and}\quad 
C(\widetilde{\chi}_{\bpi})=\overline{\chi}_{\bpi}(3)\,\widetilde{\bpi}^{-1}\,
\overline{\widetilde{\bpi}}, \]
we can deduce quite similarly to the above
\[ \alpha_{\bpi}=\overline{\chi}_{\bpi}(3)\,a_{\bpi}\quad 
\mbox{where}\ \ a_{\bpi}\in \Z\ \ \mbox{and}\ \ a_{\bpi}\meqv{-1}{3}. \]

(c)\ \ The case\ \ $p=\bpi\,\overline{\bpi}\meqv{1}{18}$.\ \ 
We know in this case
\[ \varpi_1^{-1}\,L(1,\,\widetilde{\chi}_{\bpi})
=\chi_{\bpi}(3)\,\widetilde{\bpi}^{-1}\,\alpha_{\bpi}\quad \mbox{and}\quad 
C(\widetilde{\chi}_{\bpi})=-\chi_{\bpi}(3)\,\widetilde{\bpi}^{-1}\,
\overline{\widetilde{\bpi}}, \]
so that we have\ \ 
$\overline{\chi}_{\bpi}(3)\,\alpha_{\bpi}=-
\chi_{\bpi}(3)\,\overline{\alpha_{\bpi}}$,\ \ 
which means 
\[ \alpha_{\bpi}=\chi_{\bpi}(3)\,
a_{\bpi}\,\sqrt{-3}\ \ 
\mbox{with some}\ \ a_{\bpi}\in \Z. \]
Thus the proof is completed.

\medskip

\begin{Example}{\normalfont \rmfamily 
We follow Example 3, where we evaluated the coefficient 
of the elliptic Gauss sum\ $\G_{\pi}(\chi_{\pi},\,\varphi^{-1})$ : 
$\alpha_{\bpi}=-\overline{\rho}$\ \ in the case $\bpi=4+3\rho,\ p=13$. 
Since we find\ \ $\chi_{\bpi}(3)=\rho$\ \ in this case, 
we can represent this as\ \ 
$\alpha_{\bpi}= -\overline{\rho}=\overline{\chi}_{\bpi}(3)\cdot (-1)$,\ \ 
thus we get $a_{\bpi}=-1$,\ \ which satisfies obviously 
the expected congruence\ \ $a_{\bpi}\meqv{-1}{3}$. 
Other examples by computor are given in the table in Appendix I.1.
}\end{Example} 

\noindent
{\bf Remark.}\ \ By tracing the process of the proof 
we can observe a remarkable fact. 
Under the theorems I.1, I.2 and I.3, 
the assertions of Theorem I.4 and Lemma I.6 ($G_3$-formula) 
are equivalent to each other. 
Therefore if the rationality of the elliptic Gauss sum coefficient 
could be independently proved beforehand, 
we can get Cassels-Matthews' formula as a corollary. 
It might be a natural proof of $G_3$-formula. 

\bigskip

\bigskip

\noindent
{\bf 6.2.}\ \ The substance of Theorem I.4 can be stated by the language of 
Hecke $L$-values. The following may be simply regarded 
as a precise form of Damerell's general result in a very special case. 
At the same time, however, it shows that there is a direct relation 
between the values $L(1,\,\widetilde{\chi}_{\bpi})$ and $G_3(\bpi)$, 
especially between their arguments.

\begin{Theorem} 
Let $a_{\bpi}$ be a rational integer as given in Theorem I.4. 
\begin{equation}
\varpi_1^{-1}\,L(1,\,\widetilde{\chi}_{\bpi})=
\begin{cases}
p^{1/3}\,G_3(\bpi)^{-1}\,a_{\bpi}  &
\textit{if}\quad p=\bpi\,\overline{\bpi}\meqv{7\ {\rm or}\ 13}{18}, \\
-p^{1/3}\,\sqrt{-3}\,G_3(\bpi)^{-1}\,a_{\bpi}  &
\textit{if}\quad p=\bpi\,\overline{\bpi}\meqv{1}{18}. 
\end{cases} \end{equation}
\end{Theorem}

\noindent
{\textit{Proof.}}\ \ Combining Theorems I.1--I.4 and $G_3$-formula, it can be easily verified.

\bigskip

\noindent
{\bf Corollary I.1.}\quad $L(1,\,\widetilde{\chi}_{\bpi})\ne 0$\quad 
\textit{if}\quad $p=\bpi\,
\overline{\bpi}\meqv{7\ {\rm or}\ 13}{18}$.

\medskip

\noindent
{\textit{Proof.}}\ \ Because $a_{\pi}\meqv{\pm 1}{3}$ in these cases.

\medskip

\noindent
{\bf Remark.}\ \ On the other hand, we can observe that 
$L(1,\,\widetilde{\chi}_{\bpi})$ happens often to vanish 
in the case $p=\bpi\,\overline{\bpi}\meqv{1}{18}$. 
For examples it is the case for each prime as follows : 

$p=73,\,271,\,307,\,523,\,577,\,919,\,1531,\,1549,\,1783,\,
2179,\,2287,\,2971,\,3079,\,3529,\,\ldots,$

\noindent while any reason or any rule is not known yet.

\bigskip

\bigskip

For convenience' and interest's sake, we append a table 
and a brief list of formulas. 

\medskip

\noindent
{\bf Appendix I.1.}\ \ 

\medskip

A small table of the coefficients of elliptic Gauss sums is given. 
In the table, the coefficient $\alpha_{\bpi}$ is 
expressed as\ \ $\alpha_{\bpi}=a_{\bpi}\cdot 
\chi_{\bpi}(3),\ 
a_{\bpi}\cdot \overline{\chi}_{\bpi}(3),\ 
a_{\bpi}\cdot \chi_{\bpi}(3)\cdot \sqrt{-3}$,\ \ 
for the case\ \ 
$p\meqv{7}{18},\ p\meqv{13}{18},\ p\meqv{1}{18}$,\ \ respectively. 
We can observe that the size of $a_{\bpi}$ is remarkably small. 
The computation was made by {\sc ubasic}.

\bigskip

\noindent
{\bf Appendix I.2.}\ \ 

\medskip

Addition and multiplication formulas of the functions 
$\varphi(u)$ and $\psi(u)$ are selected. Proofs are omitted, while 
a brief comment on the general complex multiplication formula is given. 
Those formulas listed without proof might be less familiar 
in comparison with the lemniscatic function case. 
It, however, is not difficult to obtain them. For example, 
we first deduce the expressions of $\wp(\varpi_1\,u)$ and $\wp'(\varpi_1\,u)$ 
in $\varphi(u), \psi(u)$ from (\ref{ppphipsi}), and substitute them 
into an ordinary addition formula of $\wp,\ \wp'$, 
e.g. the determinant formula, to derive the addition formula 
of $\varphi,\ \psi$, and so forth. 
For a general survey, one may refer to the book [L]. 
While the lemniscatic case is mainly treated there, 
the cubic case proceeds quite analogously.

\newpage

{} \hfill 
{\bf Appendix I.1.\ \ Table of the coefficients of elliptic Gauss sums}
\hfill {}

\vfill

\scriptsize
\[ \G_{\bpi}(\chi_{\bpi},\,\varphi)=\alpha_{\bpi}\,\widetilde{\bpi}^2
\quad \quad \quad 
\G_{\bpi}(\chi_{\bpi},\,\varphi^{-1})=\alpha_{\bpi}\,\widetilde{\bpi}^2
\quad \quad \quad 
\G_{\bpi}(\chi_{\bpi},\,\psi)=\alpha_{\bpi}\,\widetilde{\bpi}^2 \]
\[ \begin{array}{rrr|rrr|rrr}
p  & \boldsymbol{\pi}\quad & 
 \alpha_{\boldsymbol{\pi}} & 
p & \boldsymbol{\pi}\quad & 
 \alpha_{\boldsymbol{\pi}} & 
p & \boldsymbol{\pi}\quad & \alpha_{\boldsymbol{\pi}}\qquad \\
\hline
 7 & 1 + 3 \rho & 1 \cdot\rho & 13 & 4 + 3 \rho &-1 \cdot \overline{\rho} & 19  & -2 + 3 \rho & -1  \cdot\rho \cdot \sqrt{-3} \\
 43 & 7 + 6 \rho & 1 \cdot \overline{\rho} & 31 & 1 + 6 \rho &-1 \cdot\rho & 37  &  7 + 3 \rho & -1  \cdot\rho \cdot \sqrt{-3} \\
 61 & 4 + 9 \rho & 1 \cdot 1 & 67 & 7 + 9 \rho & 2 \cdot 1 & 73  &  1 + 9 \rho &  0  \cdot\rho \cdot \sqrt{-3} \\
 79 & 10 + 3 \rho & 1 \cdot\rho & 103 &-2 + 9 \rho & 2 \cdot 1 & 109  &  7 + 12 \rho & -1  \cdot\rho \cdot \sqrt{-3} \\
 97 &-8 + 3 \rho &-2 \cdot\rho & 139 & 13 + 3 \rho & 2 \cdot \overline{\rho} & 127  &  13 + 6 \rho &  1  \cdot \overline{\rho} \cdot \sqrt{-3} \\
 151 &-5 + 9 \rho &-2 \cdot 1 & 157 & 13 + 12 \rho &-1 \cdot \overline{\rho} & 163  & -11 + 3 \rho & -1  \cdot\rho \cdot \sqrt{-3} \\
 223 &-11 + 6 \rho & 1 \cdot \overline{\rho} & 193 & 16 + 9 \rho & 2 \cdot 1 & 181  &  4 + 15 \rho &  1  \cdot \overline{\rho} \cdot \sqrt{-3} \\
 241 & 16 + 15 \rho & 1 \cdot \overline{\rho} & 211 & 1 + 15 \rho &-1 \cdot\rho & 199  &  13 + 15 \rho &  1  \cdot \overline{\rho} \cdot \sqrt{-3} \\
 277 & 19 + 12 \rho & 1 \cdot\rho & 229 &-5 + 12 \rho &-1 \cdot \overline{\rho} & 271  &  19 + 9 \rho &  0  \cdot\rho \cdot \sqrt{-3} \\
 313 & 19 + 3 \rho & 4 \cdot\rho & 283 & 19 + 6 \rho & 2 \cdot\rho & 307  &  1 + 18 \rho &  0 \cdot\rho \cdot \sqrt{-3} \\
 331 & 10 + 21 \rho &-5 \cdot\rho & 337 & 13 + 21 \rho &-1 \cdot \overline{\rho} & 379  &  22 + 15 \rho & -2  \cdot \overline{\rho} \cdot \sqrt{-3} \\
 349 &-17 + 3 \rho & 4 \cdot\rho & 373 & 4 + 21 \rho &-4 \cdot \overline{\rho} & 397  & -11 + 12 \rho & -1  \cdot\rho \cdot \sqrt{-3} \\
 367 & 22 + 9 \rho & 1 \cdot 1 & 409 &-8 + 15 \rho &-4 \cdot\rho & 433  &  13 + 24 \rho &  1  \cdot \overline{\rho} \cdot \sqrt{-3} \\
 421 & 1 + 21 \rho & 4 \cdot\rho & 463 & 22 + 21 \rho & 2 \cdot \overline{\rho} & 487  & -2 + 21 \rho & -1  \cdot\rho \cdot \sqrt{-3} \\
 439 &-5 + 18 \rho & 4 \cdot 1 & 499 & 25 + 18 \rho & 2 \cdot 1 & 523  & -17 + 9 \rho &  0  \cdot\rho \cdot \sqrt{-3} \\
 457 & 7 + 24 \rho & 1 \cdot \overline{\rho} & 571 &-5 + 21 \rho &-1 \cdot \overline{\rho} & 541  &  25 + 21 \rho &  2  \cdot\rho \cdot \sqrt{-3} \\
 547 & 13 + 27 \rho &-2 \cdot 1 & 607 &-23 + 3 \rho & 5 \cdot \overline{\rho} & 577  &  19 + 27 \rho &  0  \cdot\rho \cdot \sqrt{-3} \\
 601 & 25 + 24 \rho &-2 \cdot \overline{\rho} & 643 &-11 + 18 \rho & 2 \cdot 1 & 613  &  28 + 9 \rho & -3  \cdot 1 \cdot \sqrt{-3} \\
 619 & 22 + 27 \rho & 4 \cdot 1 & 661 &-20 + 9 \rho & 2 \cdot 1 & 631  & -14 + 15 \rho &  1  \cdot \overline{\rho} \cdot \sqrt{-3} \\
 673 &-8 + 21 \rho & 4 \cdot\rho & 733 & 31 + 12 \rho &-1 \cdot \overline{\rho} & 739  &  7 + 30 \rho & -1  \cdot\rho \cdot \sqrt{-3} \\
 691 & 19 + 30 \rho & 1 \cdot\rho & 751 & 31 + 21 \rho & 5 \cdot \overline{\rho} & 757  &  28 + 27 \rho &  3  \cdot 1 \cdot \sqrt{-3} \\
 709 & 28 + 3 \rho &-2 \cdot\rho & 769 &-17 + 15 \rho &-1 \cdot\rho & 811  &  31 + 6 \rho &  1  \cdot \overline{\rho} \cdot \sqrt{-3} \\
 727 & 31 + 18 \rho & 1 \cdot 1 & 787 &-2 + 27 \rho &-4 \cdot 1 & 829  &  13 + 33 \rho &  1  \cdot \overline{\rho} \cdot \sqrt{-3} \\
 853 & 31 + 27 \rho & 4 \cdot 1 & 823 & 19 + 33 \rho &-4 \cdot\rho & 883  &  34 + 21 \rho & -1  \cdot\rho \cdot \sqrt{-3} \\
 907 & 7 + 33 \rho &-5 \cdot \overline{\rho} & 859 & 10 + 33 \rho & 2 \cdot\rho & 919  & -17 + 18 \rho &  0  \cdot\rho \cdot \sqrt{-3} \\
 997 & 13 + 36 \rho & 1 \cdot 1 & 877 & 31 + 3 \rho &-1 \cdot \overline{\rho} & 937  & -29 + 3 \rho &  2  \cdot\rho \cdot \sqrt{-3} \\
 1033 & 37 + 21 \rho &-2 \cdot\rho & 967 & 34 + 27 \rho & 2 \cdot 1 & 991  & -26 + 9 \rho &  3  \cdot 1 \cdot \sqrt{-3} \\
 1051 &-29 + 6 \rho &-2 \cdot \overline{\rho} & 1021 & 25 + 36 \rho & 2 \cdot 1 & 1009  & -8 + 27 \rho & -3  \cdot 1 \cdot \sqrt{-3} \\
 1069 & 37 + 12 \rho & 1 \cdot\rho & 1039 & 37 + 15 \rho &-1 \cdot\rho & 1063  &  34 + 3 \rho & -1  \cdot\rho \cdot \sqrt{-3} \\
 1087 &-17 + 21 \rho & 7 \cdot\rho & 1093 & 7 + 36 \rho & 2 \cdot 1 & 1117  &  37 + 9 \rho &  3  \cdot 1 \cdot \sqrt{-3} \\
 1123 & 34 + 33 \rho & 1 \cdot \overline{\rho} & 1129 &-32 + 3 \rho &-1 \cdot \overline{\rho} & 1153  &  16 + 39 \rho & -4  \cdot\rho \cdot \sqrt{-3} \\
 1213 & 28 + 39 \rho & 1 \cdot\rho & 1201 & 40 + 21 \rho &-4 \cdot \overline{\rho} & 1171  &  25 + 39 \rho & -1  \cdot\rho \cdot \sqrt{-3} \\
 1231 & 10 + 39 \rho & 4 \cdot\rho & 1237 & 37 + 33 \rho &-1 \cdot\rho & 1279  & -5 + 33 \rho & -2  \cdot \overline{\rho} \cdot \sqrt{-3} \\
 1249 & 40 + 27 \rho &-5 \cdot 1 & 1291 &-26 + 15 \rho & 5 \cdot\rho & 1297  &  7 + 39 \rho & -1  \cdot\rho \cdot \sqrt{-3} \\
 1303 &-14 + 27 \rho &-2 \cdot 1 & 1327 & 19 + 42 \rho &-1 \cdot\rho & 1423  &  31 + 42 \rho &  1  \cdot \overline{\rho} \cdot \sqrt{-3} \\
 1321 & 40 + 9 \rho &-2 \cdot 1 & 1381 & 4 + 39 \rho & 2 \cdot \overline{\rho} & 1459  &  43 + 30 \rho & -1  \cdot\rho \cdot \sqrt{-3} \\
 1429 & 43 + 15 \rho &-2 \cdot \overline{\rho} & 1399 & 43 + 18 \rho & 2 \cdot 1 & 1531  &  19 + 45 \rho &  0  \cdot\rho \cdot \sqrt{-3} \\
 1447 & 37 + 39 \rho & 7 \cdot\rho & 1453 &-23 + 21 \rho &-1 \cdot \overline{\rho} & 1549  &  28 + 45 \rho &  0  \cdot\rho \cdot \sqrt{-3} \\
 1483 & 1 + 39 \rho &-2 \cdot\rho & 1471 &-35 + 6 \rho &-1 \cdot\rho & 1567  & -38 + 3 \rho &  2  \cdot\rho \cdot \sqrt{-3} \\
 1609 & 13 + 45 \rho & 7 \cdot 1 & 1489 & 40 + 3 \rho &-7 \cdot \overline{\rho} & 1621  & -35 + 9 \rho & -3  \cdot 1 \cdot \sqrt{-3} \\
 1627 & 43 + 6 \rho & 1 \cdot \overline{\rho} & 1543 & 43 + 9 \rho & 2 \cdot 1 & 1657  & -23 + 24 \rho &  1  \cdot \overline{\rho} \cdot \sqrt{-3} \\
 1663 &-26 + 21 \rho &-2 \cdot\rho & 1579 & 37 + 42 \rho &-1 \cdot\rho & 1693  &  43 + 39 \rho &  2  \cdot\rho \cdot \sqrt{-3} \\
 1699 &-17 + 30 \rho & 1 \cdot\rho & 1597 & 43 + 36 \rho & 2 \cdot 1 & 1747  & -14 + 33 \rho &  1  \cdot \overline{\rho} \cdot \sqrt{-3} \\
 1753 & 19 + 48 \rho &-2 \cdot\rho & 1669 &-20 + 27 \rho & 2 \cdot 1 & 1783  &  46 + 9 \rho &  0  \cdot\rho \cdot \sqrt{-3} \\
 1789 &-35 + 12 \rho & 4 \cdot\rho & 1723 & 1 + 42 \rho &-1 \cdot\rho & 1801  &  49 + 24 \rho &  1  \cdot \overline{\rho} \cdot \sqrt{-3} \\
 1861 & 4 + 45 \rho & 1 \cdot 1 & 1741 &-5 + 39 \rho & 8 \cdot \overline{\rho} & 1873  &  49 + 33 \rho &  1  \cdot \overline{\rho} \cdot \sqrt{-3} \\
 1879 &-23 + 27 \rho & 1 \cdot 1 & 1759 & 7 + 45 \rho & 2 \cdot 1 & 1999  & -5 + 42 \rho &  1  \cdot \overline{\rho} \cdot \sqrt{-3} \\
 1933 & 49 + 36 \rho & 1 \cdot 1 & 1777 & 31 + 48 \rho &-1 \cdot \overline{\rho} & 2017  &  7 + 48 \rho &  2  \cdot\rho \cdot \sqrt{-3} \\
\end{array} \]

\vfill

\vfill

\newpage

\setlength{\textheight}{638pt}
\renewcommand{\theequation}{\arabic{equation}}
\setcounter{equation}{0}

\normalsize
\setlength{\baselineskip}{14pt}

{} \hfill 
{\bf Appendix I.2.\ \ Formulas of special elliptic functions}\ \ 
(cf. {\bf \S I.1})
\hfill {}

\bigskip

\noindent
{\bf Addition Formula}
\begin{eqnarray} 
&&\varphi(u+v)=
\dfrac{\varphi(u)^2\,\psi(v)-\varphi(v)^2\,\psi(u)}
{\varphi(u)\,\psi(v)^2-\varphi(v)\,\psi(u)^2}
=\dfrac{\varphi(v)+\varphi(u)\,\psi(u)\,\psi(v)^2}
{\psi(u)+\varphi(u)^2\,\varphi(v)\,\psi(v)} \label{addfr} \\
&&\psi(u+v)=
\dfrac{\varphi(u)\,\psi(u)-\varphi(v)\,\psi(v)}
{\varphi(u)\,\psi(v)^2-\varphi(v)\,\psi(u)^2}
=\dfrac{\psi(u)^2\,\psi(v)-\varphi(u)\,\varphi(v)^2}
{\psi(u)+\varphi(u)^2\,\varphi(v)\,\psi(v)} \\
&&\varphi(u-v)=
\dfrac{\varphi(u)^2\,\psi(v)-\varphi(v)^2\,\psi(u)}
{\varphi(u)+\varphi(v)\,\psi(v)\,\psi(u)^2}=
\dfrac{\varphi(u)\,\psi(u)-\varphi(v)\,\psi(v)}
{\psi(u)\,\psi(v)^2-\varphi(u)^2\,\varphi(v)} \\
&&\psi(u-v)=
\dfrac{\psi(u)^2\,\psi(v)-\varphi(u)\,\varphi(v)^2}
{\psi(u)\,\psi(v)^2-\varphi(u)^2\,\varphi(v)}
=\dfrac{\varphi(v)+\varphi(u)\,\psi(u)\,\psi(v)^2}
{\varphi(u)+\varphi(v)\,\psi(v)\,\psi(u)^2} 
\end{eqnarray}

\noindent
{\bf Multiplication Formula}
\begin{eqnarray} 
&&\varphi(\rho\,u)=\rho\,\varphi(u),\quad \psi(\rho\,u)=\psi(u),\quad 
\Ze(\rho\,u)=\overline{\rho}\,\Ze(u) \\
&&\varphi(-u)=-\varphi(u)\,\psi(u)^{-1},\quad 
\psi(-u)=\psi(u)^{-1},\quad \Ze(-u)=-\Ze(u) \\
&&\varphi(-2u)=\varphi(u)\cdot \dfrac{\varphi(u)^3-2}{1-2\,\varphi(u)^3},\quad
\psi(-2u)=\psi(u)\cdot \dfrac{\psi(u)^3-2}{1-2\,\psi(u)^3} \\
&&\varphi(\sqrt{-3}\,u) = 
\dfrac{\sqrt{-3}\,\varphi(u)\psi(u)}{1+\overline{\rho}\,\varphi(u)^3},\quad 
\psi(\sqrt{-3}\,u) = \dfrac{\rho+\psi(u)^3}{1+\rho\,\psi(u)^3} \\
&& \Ze((1-\rho)\,u)=(1-\overline{\rho})\,\Ze(u)
+(1-\overline{\rho})\bigl\{\varphi(u)^{-1}-\varphi(-u)^{-1} \bigr\}
\end{eqnarray}

\medskip

\noindent
{\bf Primary Prime Multiplication :}
\quad $p=\boldsymbol{\pi}\,\overline{\boldsymbol{\pi}},\ \ 
\boldsymbol{\pi}\meqv{1}{3}$
\begin{gather}
 \varphi(\boldsymbol{\pi}\,u) = 
\varphi(u) \prodp_{\nu\,({\rm mod}\,\boldsymbol{\pi})} 
\varphi(u+\nu/\boldsymbol{\pi}),\quad 
\psi(\boldsymbol{\pi}\,u) = 
\psi(u) \prodp_{\nu\,({\rm mod}\,\boldsymbol{\pi})} 
\psi(u+\nu/\boldsymbol{\pi}) \\
 \varphi(\boldsymbol{\pi}\,u) = 
\varphi(u)\cdot \dfrac{U(\varphi(u))}{R(\varphi(u))},\quad 
\psi(\boldsymbol{\pi}\,u) = 
\psi(u)\cdot \dfrac{U(\psi(u))}{R(\psi(u))} \label{cm1} \\
 U(x)=\prodp_{\nu\,({\rm mod}\,\boldsymbol{\pi})} 
(x-\varphi(\nu/\boldsymbol{\pi})),\quad R(x)=x^{p-1}\,U(x^{-1}) \\
 U(x) \in \Or[x],\ \ \ U(x)\meqv{x^{p-1}}{\boldsymbol{\pi}},\ \ 
U(0)=\boldsymbol{\pi} 
\end{gather}

\medskip

\noindent
\textit{Proof of} (10), (11)\ :\ \ 
\small 
By comparing the divisors we have (10).\ \ By using the first form of (1),
\[ \prod_{\varepsilon\in W'}\varphi(u+\varepsilon\,v)
=\dfrac{\psi(v)^3\left(\varphi(u)^3-\varphi(-v)^3\right)}
{1-\varphi(u)^3\,\varphi(v)^3},\quad 
\prod_{\varepsilon\in W'}\varphi(u-\varepsilon\,v)
=\dfrac{\varphi(u)^3-\varphi(v)^3}
{\psi(v)^3\left(1-\varphi(u)^3\,\varphi(-v)^3\right)}. \]
\[ \therefore\ \ \prod_{\varepsilon\in W} \varphi(u+\varepsilon\,v)
=\dfrac{\left(\varphi(u)^3
-\varphi(-v)^3\right)\left(\varphi(u)^3-\varphi(v)^3\right)}
{\left(1-\varphi(u)^3\,\varphi(v)^3\right)
\left(1-\varphi(u)^3\,\varphi(-v)^3\right)}
=\prod_{\varepsilon\in W} \dfrac{\varphi(u)-\varphi(\varepsilon\,v)}
{1-\varphi(u)\,\varphi(\varepsilon\,v)}. \]

This combined with (10) leads to (11) as follows ; here 
$U$ is an arbitrary $\tfrac{1}{6}$-set mod $\boldsymbol{\pi}$.
\[ \dfrac{\varphi(\boldsymbol{\pi}\,u)}{\varphi(u)}=
\prod_{\nu\in U}\prod_{\ve\in W} \varphi(u+\ve\,\nu/\boldsymbol{\pi})=
\prod_{\nu\in U}\prod_{\ve\in W} 
\dfrac{\varphi(u)-\varphi(\ve\nu/\boldsymbol{\pi})}
{1-\varphi(\ve\nu/\boldsymbol{\pi})\,\varphi(u)}=
\prodp_{\nu\,({\rm mod}\,\boldsymbol{\pi})} 
\dfrac{\varphi(u)-\varphi(\nu/\boldsymbol{\pi})}
{1-\varphi(\nu/\boldsymbol{\pi})\,\varphi(u)}.
 \]

\newpage


\renewcommand{\theequation}{\Roman{section}.\arabic{equation}}
\setcounter{equation}{0}
\setcounter{Lemma}{0} \setcounter{Theorem}{0} 

\normalsize
\setlength{\baselineskip}{14pt}

\section{The Quartic Character Case}
\medskip
\bigskip

Throughout the part {\bf II}, the field $\Q(i)$ and the ring $\Z[i]$ 
are abbreviated to $F$ and $\Or$, respectively. The unit group is 
denoted by\ $W=\{\pm 1,\,\pm i\}$. 
The set $\Or$ or its constant multiple appears also as a period 
lattice for elliptic functions. 
Though we don't treat the octic case, we shall come on a scene 
to need the eighth root of unity\ 
$\zeta_8=e^{2\pi i/8}$\ and so it is not strange 
to meet $\sqrt{2}=(1-i)\,\zeta_8$\ \ or\ \ 
$i\sqrt{2}=(1+i)\,\zeta_8$\ \ in some formulas. 

\subsection{Special elliptic functions with complex multiplication}

\noindent
{\bf 1.1.}\quad We shall define some special functions which play 
the leading role in our argument. 
Let\ $\wp(u)$\ denote the Weierstrass function respect to the period lattice\  
$\varpi\Or$\ so that\ \ $\wp'(u)^2=4\wp(u)^3-4\wp(u)$ ;\ the real period 
$\varpi$ is given by
\[ \varpi \doteq 2\int_0^1\dfrac{d\,x}{\sqrt{1-x^4}}
=2.62205755 \cdots. \]
It is obvious\ $\wp(\varpi u)$\ and\ $\wp'(\varpi u)$\ 
are elliptic functions of periods\ $\Or$. 
Further, we can get another doubly periodic function 
by a slight modification of Weierstrass' $\zeta$ of $\varpi\Or$ :

\noindent
{\bf Definition.}\quad The non-analytic but doubly periodic function\ 
$\Ze(u)$\ is defined by
\begin{equation}
 \Ze(u) \doteq 
\zeta(\varpi u)-\dfrac{\pi}{\varpi}\,\overline{u}, \label{Z}
\end{equation}
Double periodicity of $\Ze$ relative to $\Or$ is easily verified 
by a usual formula of $\zeta$. 
The addition formula of\ $\Ze$,\ that follows immediately from one 
of $\zeta$,\ is useful :
\[ \Ze(u+v)=\Ze(u)+\Ze(v)+\dfrac{1}{2}\,\dfrac{\wp'(\varpi u)-
\wp'(\varpi v)}{\wp(\varpi u)-\wp(\varpi v)}. \]
In particular, such a function\ \ $\sum_{k=1}^r c_k\,\Ze(u+\gamma_k)$ 
is an elliptic function if $\sum_{k=1}^rc_k=0$.

\medskip

The following two functions are specially important, really 
which are nothing but the old lemniscatic sine and cosine functions of Gauss. 

\medskip

\noindent
{\bf Definition.}\quad The elliptic functions\ $\varphi(u)$\ and\ 
$\psi(u)$\ to the period $\Or$ are defined by
\begin{eqnarray} 
\varphi(u) & \doteq & -\dfrac{1-i}{2}\left\{
\mathsf{Z}\left(u-\dfrac{1}{2}\right)-\mathsf{Z}\left(u-\dfrac{i}{2}\right)
 \right\}, \label{phi} \\
\psi(u) & \doteq & -\dfrac{1-i}{2}\left\{
\mathsf{Z}\left(u-\dfrac{1-i}{4}\right)
-\mathsf{Z}\left(u+\dfrac{1-i}{4}\right) \right\}. \label{psi}
\end{eqnarray}

\medskip

From the addition formula of $\Ze(u)$\ \  
we can easily derive the other expressions :
\begin{equation}
\varphi(u)=-2(1-i)\cdot\dfrac{\wp(\varpi u)}{\wp'(\varpi u)},\quad 
\psi(u)=\dfrac{\wp(\varpi u)+i}{\wp(\varpi u)-i}. 
\end{equation}

We need also the following formula.
\begin{equation}
\varphi(u)^{-1}=
\dfrac{1+i}{2} \left\{\mathsf{Z}\left(u \right) 
- \mathsf{Z}\left(u-\dfrac{1+i}{2}\right) \right\}. \label{peinv}
\end{equation}

\medskip

Here are some basic properties of these functions, which are all classical 
or easily deduced from the definitions 
and by usual theory of elliptic functions. 
\[ {\rm div}\,(\varphi)=(0)+((1+i)/2)-(1/2)-(i/2), \]
\[ {\rm div}\,(\psi)=((1+i)/4)+(-(1+i)/4)-((1-i)/4)-(-(1-i)/4), \]
\[ \Ze(i u)=-i\,\Ze(u),\quad 
\varphi(i\,u)=i\,\varphi(u),\quad \psi(i\,u)=\psi(u)^{-1}, \]
\[ \psi(u)=\varphi(u+(1+i)/4),\quad \varphi(u)^{-1}=-i\,\varphi(u+1/2), \]
\[ \varphi'(u)=(1-i)\varpi\,\psi(u)\,(1+\varphi(u)^2),\quad 
\psi'(u)=-(1-i)\varpi\,\varphi(u)\,(1+\psi(u)^2), \]
\[ \varphi(u)^2\psi(u)^2+\varphi(u)^2+\psi(u)^2=1. \]

\medskip

In particular\ $(((1-i)\varpi)^{-1}\varphi'(u))^2=1-\varphi(u)^4$\ holds, 
and so we can ascertain that\ $\varphi(u)={\rm sl}((1-i)\varpi u)$\ and\ 
$\psi(u)={\rm cl}((1-i)\varpi u)$\ by Gauss' 
lemniscatic sine and cosine. 

\medskip

We can refer to the survey monograph [L] by F.\,Lemmermeyer for further 
general facts and some background of these elliptic functions.

\bigskip

\bigskip

\noindent
{\bf 1.2.}\quad 
Let $\bpi$ be a complex prime in $\Or$ so that 
$p=\bpi\,\overline{\bpi}\meqv{1}{4}$, 
and we assume also $\bpi$ is primary : 
$\bpi\meqv{1}{(1+i)^3}$. 
Then we have\ \ $(\Or/(\bpi))^{\times} 
\cong (\Z/p\Z)^{\times}$. 
We often abbreviate like as\ \ 
$\nu\,({\rm mod}\,\bpi)$\ \ in such a case 
when $\nu$ runs over $(\Or/(\bpi))^{\times}$.

\medskip

It is well known that such a division value\ 
$\varphi(1/\bpi)$\ or\ $\psi(1/\bpi)$\ generates 
an abelian extension of the imaginary quadratic field $F=\Q(i)$. In fact, 

\begin{Lemma} \ \ 
Let\ \ $L=F(\varphi(1/\bpi))$\ \ and\ \ 
$L_1=F(\psi(1/\bpi))$. One has
\begin{enumerate}
\item[{\rm (i)}] $L/F$ is a cyclic extension of degree $p-1$, and 
$L_1$ is a subfield of $[L : L_1]=2$.
\item[{\rm (ii)}] ${\rm Gal}(L/F) \cong 
(\Or/(\bpi))^{\times}$\ 
by corresponding\ $\sigma_{\mu}$ to $\mu$, and it holds\ 
$\varphi(\nu/\bpi)^{\sigma_{\mu}}
=\varphi(\mu\nu/\bpi)$,\ \ 
$\psi(\nu/\bpi)^{\sigma_{\mu}}=\psi(\mu\nu/\bpi)$\ \ 
for arbitrary\ \ $\mu, \nu\in (\Or/(\bpi))^{\times}$.
\item[{\rm (iii)}] $\varphi(1/\bpi),\ \psi(1/\bpi)$ are algebraic integers, 
and particularly\ \ $\psi(1/\bpi)$ is a unit.
\item[{\rm (iv)}] The prime ideal $(\bpi)$ splits completely in $L$ :\ 
$(\bpi)=\mathfrak{P}^{p-1}$ where\ \ $\mathfrak{P}=(\varphi(1/\bpi))$.
\end{enumerate}
\end{Lemma}

$L$ and $L_1$ are called the ray class fields of the conductors\ 
$((1+i)^3\bpi)$ and $((1+i)^2\bpi)$, respectively. 
The proof is omitted, but it may be found mostly in [L]. 
We only give some numerical examples in the end of this section. 

\medskip

The following is a minor property of division values and perhaps 
classically known, but it is important for us and hence 
we state it as a lemma with a brief proof : 

\begin{Lemma} \ For arbitrary\ \ 
$\mu, \nu\in (\Or/(\bpi))^{\times}$, 
\begin{equation}
\varphi(\mu/\bpi)\,\varphi(\nu/\bpi)\meqv{1}{(1+i)}. 
\label{Lemm2}
\end{equation}
\end{Lemma}

\noindent
\textit{Proof.}\ \ Recall the $(1+i)$-multiplication formula of 
$\varphi(u)$ : 
$\varphi((1+i)\,u)=(1+i)\,\varphi(u)\,\psi(u)\cdot (1-\varphi(u)^2)^{-1}$.\ \ 
Substituting\ \ $u=\nu/\bpi$\ \ and cancelling by $\mathfrak{P}$, 
we obtain an ideal equality 
$(1-\varphi(\nu/\bpi)^2)=(1+i)$, which implies (\ref{Lemm2}) 
of the case $\mu=\nu$. For the other case, it suffices to show the following.
\[ \varphi(\mu/\bpi) \equiv \varphi(\nu/\bpi)
\pmod{(1+i)}. \]
To prove this, we may assume $\nu=1$. Let $\boldsymbol{\lambda}$ be 
a primary prime such that 
$\boldsymbol{\lambda}\meqv{\mu}{\bpi}$. 
Then the complex multiplication formula (cf. Example 1) says :\ \ 
$\varphi(\boldsymbol{\lambda}\,u)
=\varphi(u)\,U(\varphi(u))\,R(\varphi(u))^{-1}$, 
where $U(x),\ R(x)$ are the polynomials of $x^4$ reciprocal to each other 
over $\Z[i]$, so that $U(1)=R(1)$. Combining this with the fact\ 
$\varphi(1/\bpi)^2\meqv{1}{(1+i)}$\ by (\ref{Lemm2}), we can deduce \ 
$\varphi(\boldsymbol{\lambda}/\bpi)\meqv{\varphi(1/\bpi)}{(1+i)}$. 

\bigskip

As for the division value $\Ze(1/\bpi)$, we have 
\begin{Lemma} \ 
The value $\Ze(1/\bpi)$ generates 
the same extension $L$. Namely, 
\[ L=F(\varphi(1/\bpi))=F(\Ze(1/\bpi))\quad 
\textit{and} \quad \Ze(\nu/\bpi)^{\sigma_{\mu}}
=\Ze(\mu\nu/\bpi)\quad 
(\,\mu, \nu \in (\Or/(\bpi))^{\times}\,). \]
\end{Lemma}

\noindent
\textit{Proof.}\ \ The following is obtained 
by substituting $v=iu$ in the addition formula of $\Ze$.
\begin{equation}
\mathsf{Z}((1+i)u)=(1-i)\,\mathsf{Z}(u) + i\,\varphi(u)^{-1}. \label{Z-phi}
\end{equation}
Using formula (\ref{Z-phi}) repeatedly 
and in view of $(1+i)^{p-1}\meqv{1}{\bpi}$, we have 
\begin{equation*}
(1-(-4)^{(p-1)/4})\,\Ze(\nu/\bpi)=i\sum_{k=1}^{p-1} 
(1-i)^{p-1-k}\,\varphi((1+i)^{k-1}\nu/\bpi)^{-1}.
\end{equation*}
Now it is easy to see the assertion of Lemma II.3.

\begin{Example}{\normalfont \rmfamily 
\ In each case tabulated below, we have the 
$\bpi$-multiplication formula :
\[ \varphi(\bpi\,u)
=\varphi(u)\cdot \dfrac{U(\varphi(u))}{R(\varphi(u))},\quad  
R(x)=x^{p-1}\,U(x^{-1}),\]
\[ U(x)=\prodp_{\nu\,({\rm mod}\,\bpi)} (x-\varphi(\nu/\bpi)),\quad 
V(x)^2=\prodp_{\nu\,({\rm mod}\,\bpi)} (x-\psi(\nu/\bpi)), \]
\[\textit{where}\quad\quad x\,U(x)-R(x)=(x-1)\,V(x)^2. \]
Furthermore it is easy to check 
\[ U(x),\ V(x) \in \Or[x],\quad U(x)\meqv{x^{p-1}}{(\bpi)},\quad 
U(0)=\bpi,\ \ V(0)=1. \]
In fact,\ \ $U(x)$ and $V(x)$ are the minimal polynomials of 
$\varphi(1/\bpi)$ and $\psi(1/\bpi)$, respectively.

\small
\[ \begin{array}{cc|l}
p & \bpi & \quad U(x)\\ \hline
5 & -1+2i & x^4+\bpi \\
13 & 3+2i & x^{12}-(1-4i)\,\bpi\,x^8+(1-2i)\,\bpi\,x^4+\bpi \\
17 & 1+4i & x^{16}-(4+4i)\,\bpi\,x^{12}+(6+4i)\,\bpi\,x^8-(4-4i)\,\bpi\,x^4+\bpi \\[2mm]
p & \bpi & \quad V(x) \\ \hline
5 & -1+2i & x^2+(1-i)\,x+1 \\
13 & 3+2i & x^6-(1+i)\,x^5-(1+2i)\,x^4-4i\,x^3-(1+2i)\,x^2-(1+i)\,x+1 \\
17 & 1+4i & x^8-2i\,x^7+(2-2i)\,x^6+(4+2i)\,x^5+2\,x^4+(4+2i)\,x^3+(2-2i)\,x^2-2i\,x+1 \\
\end{array} \]
\normalsize
}\end{Example}

\vspace{5mm}

\subsection{\boldmath{$L$}-series for Hecke characters of weight one}

\noindent
{\bf 2.1.}\ \ 
Let $\widetilde{\chi}$ denote a Hecke character of weight 1 relative to 
the modulus $(\beta)\subset \Or$, namely it is a multiplicative function 
on the ideal group of $\Or$ 
of the following form :
\[ \widetilde{\chi}((\nu))=\chi_1(\nu)\overline{\nu},\ \ 
\chi_1 : (\Or/(\beta))^{\times} \rightarrow {\bf C}^{\times},\ \ 
\chi_1(\ve)=\ve\ \ (\,\ve\in W\,), \]
where $\chi_1$ is an ordinary residue class character to the modulus 
$(\beta)$, and $(\beta)$ is called the coductor 
of\ \ $\widetilde{\chi}$\ \ if\ \ $\chi_1$\ \ is a primitive character 
to the modulus $(\beta)$. 

It is well known the associated $L$-series has the analytic continuation 
and satisfies a functional equation. 
We follow Weil's argument and his notation.
\begin{eqnarray}
L(s,\ \widetilde{\chi}) &=& \sum_{\mathfrak{a}} 
\widetilde{\chi}(\mathfrak{a})\,N\mathfrak{a}^{-s}
= \dfrac{1}{4} \sum_{\nu\in \Or} \chi_1(\nu)\,\overline{\nu}\,|\nu|^{-2s} 
\notag \\
&=& \dfrac{1}{4} \sum_{\lambda\,({\rm mod}\,\beta)} \chi_1(\lambda)\,
\sum_{\mu\in\Or} (\overline{\lambda}+\overline{\mu}\overline{\beta})\,
|\lambda + \mu\beta|^{-2s} \notag \\
\therefore\ \ L(s,\ \widetilde{\chi}) &=& \beta^{-1}\,
|\beta|^{2-2s}\cdot \dfrac{1}{4}
\sum_{\lambda\,({\rm mod}\,\beta)} \chi_1(\lambda)\,
K_1(\lambda/\beta,\ 0,\ s) \label{L1g}
\end{eqnarray}
Here the function $K_1$ is defined for ${\rm Re}\,s>3/2$ as follows, and 
it is analytically continued to the whole $s$-plane and satisfies 
the own functional equation :  (cf. [W], VIII)
\begin{gather*}
K_1(u, u_0, s) = \sum_{\mu \in \Or}
e^{\pi\,(\overline{u}_0\mu-u_0\overline{\mu})}\,
(\overline{u}+\overline{\mu})\,|u+\mu|^{-2s},  \\
\pi^{-s}\,
\varGamma(s)\, K_1(u, u_0, s)
= e^{\pi\,(\overline{u}_0 u-u_0\overline{u})}\, \pi^{s-2}\,
\varGamma(2-s)\, K_1(u_0, u, 2-s).
\end{gather*}
When $(\beta)$ is the conductor of $\widetilde{\chi}$, 
a usual computation of 
Gauss sum works. From this combined with the above, 
the functional equation of Hecke $L$-series is derived : 
\[ \varLambda(s,\ \widetilde{\chi})
=C(\widetilde{\chi})\,\varLambda(2-s,\ \overline{\widetilde{\chi}}), \]
\[ \textit{where}\quad \varLambda(s,\ \widetilde{\chi})
=\left(\dfrac{2\pi}{\sqrt{4\cdot N(\beta)}}\right)^{-s}\,
\varGamma(s)\,L(s,\ \widetilde{\chi}), \]
\[ \textit{and}\quad 
C(\widetilde{\chi})=-i\,\beta^{-1}\,\sum_{\lambda\,({\rm mod}\,\beta)}
 \chi_1(\lambda)\,e^{2\pi i\,{\rm Re}(\lambda /\beta)}. \]

\bigskip

In particular, we have a simple equality of Hecke $L$-values at $s=1$ : 
\begin{Lemma}\ Let $\widetilde{\chi}$ be a Hecke character 
of weight $1$ with the conductor $(\beta)$. Then
\begin{equation}
L(1,\ \widetilde{\chi})=C(\widetilde{\chi})\,
\overline{L(1,\ \widetilde{\chi})}, \label{centeq} 
\end{equation}
\begin{equation}
\mbox{where}\quad 
C(\widetilde{\chi})=-i\,\beta^{-1}\,
\sum_{\lambda\,({\rm mod}\,\beta)}
 \chi_1(\lambda)\,e^{2\pi i\,{\rm Re}(\lambda /\beta)}. \label{the_root}
\end{equation}
\end{Lemma}

\medskip

\noindent
{\bf Remark.}\quad $L(1,\ \overline{\widetilde{\chi}})=
\overline{L(1,\ \widetilde{\chi})}$.

\medskip

\noindent
We call (\ref{centeq}) \textit{the central value equation}, 
and the constant $C(\widetilde{\chi})$ \textit{the root number}.

\bigskip

\bigskip

\noindent
{\bf 2.2.}\ \ On the other hand, we can notice that 
the value $L(1,\ \widetilde{\chi})$ relates to some elliptic functions. 
As is remarked in [W] (VIII,\,\S 14) or in others, the following is valid. 
\[ E_1^*(u) \doteq K_1(u,\ 0,\ 1)=\varpi\,\zeta(\varpi u)
-\pi\,\overline{u}. \]
By the definition (\ref{Z}) the right-hand side is nothing but our function\ 
$\varpi\,\Ze(u)$. Combining this with the equation (\ref{L1g}) at\ $s=1$,\ 
we obtain the following formula : 
\begin{Lemma}\ Let $\widetilde{\chi}$ be a Hecke character 
of weight $1$ with the conductor $(\beta)$. Then
\begin{equation}
\varpi^{-1}\,L(1,\ \widetilde{\chi})=\dfrac{1}{4\,\beta} 
\sum_{\lambda\,({\rm mod}\,\beta)} \chi_1(\lambda)\,\Ze(\lambda/\beta).
\label{EGSg}
\end{equation}
\end{Lemma}

\medskip

\begin{Example}{\normalfont \rmfamily 
\ The following is probably the simplest case and the derived formula\ \ 
$L(1,\,\widetilde{\chi}_0)=\tfrac{\varpi}{4}$\ \ may be compared with 
the classical formula : 
$1-\tfrac{1}{3}+\tfrac{1}{5}-\tfrac{1}{7}+\cdots=\tfrac{\pi}{4}$.

\medskip

The Hecke character\ $\widetilde{\chi}_0$\ of the conductor 
$(1+i)^3$\ is given as follows :
\[ \widetilde{\chi}_0((\nu)) \doteq \chi_0(\nu)\,\overline{\nu},\ \ 
\textit{where}\ \ \chi_0\ :\ (\Or/(1+i)^3)^{\times} \cong W\ \ 
\textit{is the natural isomorphism.} \]
Then we can evaluate 
the $L$-value at $s=1$\ directly by (\ref{EGSg}) :
\[ \varpi^{-1}\,L(1,\ \widetilde{\chi}_0)=\dfrac{1}{4(1+i)^3}
\sum_{\varepsilon\in W} \varepsilon\,\Ze(-(1+i)\varepsilon/4)
=\dfrac{1+i}{4}\,\Ze((1+i)/4)=\dfrac{1}{4}. \]
Also we can easily check\ \ $\displaystyle{
C(\widetilde{\chi}_0)=-i\,(1+i)^{-3}\,\sum_{\varepsilon\in W}\,
\varepsilon\,e^{2\pi i\,{\rm Re}(-(1+i)\varepsilon/4)}=1}$\ \ as is expected.
}\end{Example}

\vspace{5mm}

\subsection{Elliptic Gauss sums for quartic characters}

\noindent
{\bf 3.1.}\ \ Let $\bpi$ be a primary prime in $\Or$ ; $\bpi\meqv{1}{(1+i)^3}$. Let $\chi_{\bpi}$ be the quartic residue character to the modulus $(\bpi)$ and  the notation will be fixed throughout : 
\[ \chi_{\bpi}(\nu)=\qr{\nu}{\bpi}_{\!\!4}\ :\ 
\chi_{\bpi}(\nu)^4=1\ \ \textit{and}\ \ 
\chi_{\bpi}(\nu)\meqv{\nu^{(p-1)/4}}{\bpi}\ \ 
(\nu\in (\Or/(\bpi))^{\times}). \]

Let $f(u)$ be a doubly periodic function of the period $\Or$, 
which we specify below.

\medskip

\noindent
{\bf Definition.}\quad The following is called an \textit{elliptic Gauss sum}.
\begin{equation}
\G_{\bpi}(\chi_{\bpi},\,f) \doteq 
\dfrac{1}{4}\,\sum_{\nu\,({\rm mod}\,\bpi)}\,\chi_{\bpi}(\nu)\,f(\nu/\bpi).
\end{equation}

In the part {\bf II}, we deal with the three types of elliptic Gauss sums\ 
$\G_{\bpi}(\chi_{\bpi},\,\varphi)$,\ 
$\G_{\bpi}(\chi_{\bpi},\,\Ze)$,\ 
$\G_{\bpi}(\chi_{\bpi},\ \psi)$, and one more 
$\G_{\bpi}(\chi_{\bpi},\,\varphi^{-1})$\ 
for a supplementary use. So we understand that 
$f(u)$ denotes one of these functions $\varphi(u)$, $\Ze(u)$, 
$\psi(u)$, and $\varphi(u)^{-1}$ in the subsequence. 
In these cases, if ``the parity condition" is not satisfied, 
$\G_{\bpi}(\chi_{\bpi},\,f)$ vanishes trivially. 
Since $\varphi(i\,u)=i\,\varphi(u),\ \Ze(i\,u)=-i\,\Ze(u),\ 
\psi(-u)=\psi(u)$\ and\ 
$\chi_{\bpi}(i)=i^{(p-1)/4}$, we can easily check that 
$\G_{\bpi}(\chi_{\bpi},\,f)$ is not trivial only 
in the following cases. 
The parity condition, however, is not sufficient for non-vanishing of 
the elliptic Gauss sum as we shall see later.

\medskip

\textit{The elliptic Gauss sums that we shall consider are the following} : 

\medskip

(a)\ \ \ $\G_{\bpi}(\chi_{\bpi},\,\varphi)$\quad 
\textit{for the case}\ \ 
$p=\bpi\,\overline{\bpi}\meqv{13}{16}$,

\medskip

(b)\ \ \ $\G_{\bpi}(\chi_{\bpi},\,\Ze),\ \ 
\G_{\bpi}(\chi_{\bpi},\,\varphi^{-1})$\quad 
\textit{for the case}\ \ 
$p=\bpi\,\overline{\bpi}\meqv{5}{16}$\ 
\textit{and}\ $p>5$,

\medskip

(c)\ \ \ $\G_{\bpi}(\chi_{\bpi},\,\psi)$\quad 
\textit{for the case}\ \ 
$p=\bpi\,\overline{\bpi}\meqv{1}{8}$.

\medskip

As noted in Lemma II.1 and Lemma II.3, 
$\G_{\bpi}(\chi_{\bpi},\,f) \in L$ and 
$f(\mu\nu/\bpi)=f(\nu/\bpi)^{\sigma_{\mu}}$ are valid, 
and hence we can immediately deduce the property of Lagrange's resolvent : 
\begin{equation}
\G_{\bpi}(\chi_{\bpi},\,f)^{\sigma_{\mu}}=\overline{\chi}_{\bpi}(\mu)\,
\G_{\bpi}(\chi_{\bpi},\,f)\quad (\,\mu\in (\Or/(\bpi))^{\times}\,) \label{lagr}
\end{equation}
In particular, $\G_{\bpi}(\chi_{\bpi},\,f)^4$ 
is an element of $F$, and furthermore we have

\begin{Lemma}\ 
$\G_{\bpi}(\chi_{\bpi},\,f)^4$ is an algebraic integer in $\Or$\ 
for each\ $f=\varphi,\ \varphi^{-1}$\ and\ \ $\psi$. 
\end{Lemma}

\noindent
\textit{Proof.}\ \ We must show the algebraic integrity of 
$\G_{\bpi}(\chi_{\bpi},\,f)$ for each case. 
Let $S$ be an arbitrary quarter subset mod $(\bpi)$, 
namely\ $(\Or/(\bpi))^{\times}=S\cup -S\cup i\,S\cup -i\,S$\ by definition. 

\medskip

\noindent
(a)\ \ The case $p=\bpi\,
\overline{\bpi}\meqv{13}{16}$. 

\noindent 
We have 
$\G_{\bpi}(\chi_{\bpi},\,\varphi)
=\sum_{\nu\in S} \chi_{\bpi}(\nu)\,
\varphi(\nu/\bpi)$. 
Since $\varphi(\nu/\bpi)$'s are algebraic integers, 
the integrity is obvious in this case.

\medskip

\noindent
(b)\ \ The case $p=\bpi\,
\overline{\bpi}\meqv{5}{16},\ \ p>5$. 

\noindent
We have 
$\G_{\bpi}(\chi_{\bpi},\,\varphi^{-1})
=\sum_{\nu\in S} \chi_{\bpi}(\nu)\,
\varphi(\nu/\bpi)^{-1}$. 
Since $\varphi(1/\bpi)\cdot \varphi(\nu/\bpi)^{-1}$ 
is a unit,\ 
$\mathfrak{P}^4\,(\G_{\bpi}(\chi_{\bpi},\,
\varphi^{-1})^4)$ is an integral ideal, where $\mathfrak{P}
=(\varphi(1/\bpi))$. On the other hand, 
$(\bpi)=\mathfrak{P}^{p-1}$ and $p-1>4$,\ \ 
and hence $(\G_{\bpi}(\chi_{\bpi},\,\varphi^{-1})^4)$ 
is already integral. We need the condition $p>5$ in this case ; 
in fact $\G_{\bpi}(\chi_{\bpi},\,\varphi^{-1})^4=-\bpi^{-1} \not\in \Or$ 
in the case $\bpi=-1+2i\ \ (p=5)$.

\medskip

\noindent
(c)\ \ The case $p=\bpi\,
\overline{\bpi}\meqv{1}{8}$. 

\noindent
First of all we have 
$\G_{\bpi}(\chi_{\bpi},\,\psi)=\dfrac{1}{2}
\sum_{\nu\in S} \chi_{\bpi}(\nu)\{\psi(\nu/\bpi)+\chi_{\bpi}(i)\,
\psi(i\nu/\bpi)\}$. 
We must show the summation of the right side is divisible by $2=(1-i)(1+i)$. 
Note $\chi_{\bpi}(i)=\pm 1$. 
By using $(1+i)$-multiplication formula of $\varphi(u)$, it is easy to verify 
\[ \psi(u)-\psi(iu)=-(1-i)\,\varphi((1+i)\,u)\,\varphi(u),\quad 
\psi(u)+\psi(iu)=(1-i)\,\varphi((1+i)\,u)\,\varphi(u)^{-1}. \]
By Lemma II.2,\ \ $\varphi((1+i)\,\nu/\bpi)\,
\varphi(\nu/\bpi)\equiv 
\varphi((1+i)\,\nu/\bpi)\,
\varphi(\nu/\bpi)^{-1}\meqv{1}{(1+i)}$\ \ holds. 
Thus $\sum_{\nu\in S}\chi_{\bpi}(\nu)\,
\varphi((1+i)\,\nu/\bpi)\,
\varphi(\nu/\bpi)^{\pm 1}\equiv \sum_{\nu\in S} 1\equiv 
(p-1)/4\meqv{0}{(1+i)}$, which implies the integrity of 
$\G_{\bpi}(\chi_{\bpi},\,\psi)$. 
Thus the proof of Lemma II.6 is completed. 

\medskip

As a matter of fact, it is also valid 

\medskip

\noindent
{\bf Claim ($\Ze$)}\ :\quad 
$\G_{\bpi}(\chi_{\bpi},\,\Ze)^4\ \ 
\textit{is an algebraic integer in}\ \Or$. \ \ 

\medskip

\noindent
Although the proof is a bit indirect and will be completed 
after Theorem II.4, 
we shall often assume the claim for convenience' sake. 
We here only prepare the following. 
\begin{Lemma}\ One has
\begin{equation}
((1+i)-i\,\overline{\chi}_{\bpi}(1+i))\,
\G_{\bpi}(\chi_{\bpi},\,\Ze)
=\G_{\bpi}(\chi_{\bpi},\,\varphi^{-1}).
\end{equation} \end{Lemma}

\noindent
\textit{Proof.}\ \ This is immediately derived from 
the addition formula (\ref{Z-phi}) of $\Ze(u)$. 

\begin{Lemma}\ We have\ \ 
$\G_{\bpi}(\chi_{\bpi},\,\varphi)^4\meqv{1}{2}$\ \ 
and\ \ 
$\G_{\bpi}(\chi_{\bpi},\,\varphi^{-1})^4\meqv{1}{2}$,\ \ for\ \ 
$p=\bpi\,\overline{\bpi}\meqv{13}{16}$\ \ and\ \ 
$p=\bpi\,\overline{\bpi}\meqv{5}{16}$,\ \ 
respectively. 
\end{Lemma}
\noindent
\textit{Proof.}\quad By Lemma II.2 we can deduce\ \ 
$\varphi(\nu/\bpi)^4\meqv{1}{(1+i)^2}$. Hence we obtain 
\[ \G_{\bpi}(\chi_{\bpi},\,\varphi)^4\equiv 
\sum_{\nu\in S} \varphi(\nu/\bpi)^4\equiv (p-1)/4\meqv{1}{2}. \]
Similarly, using\ \ $\bpi\,\varphi(\nu/\bpi)^{-4}
\meqv{1}{2}$,\ \ we can prove 
$\G_{\bpi}(\chi_{\bpi},\,\varphi^{-1})^4\meqv{1}{2}$.

\bigskip

\bigskip

\noindent
{\bf 3.2.}\quad Obviously from (\ref{lagr}), the value 
$\G_{\bpi}(\chi_{\bpi},\,f)$ belongs to the quartic extension field over $F$. 
So it is convenient to give a suitable quartic root of\ $-\bpi$\ 
or\ $\bpi$\ for the precise investigation of the value of 
$\G_{\bpi}(\chi_{\bpi},\,f)$. 
We shall even define ``the canonical quartic root" of\ $-\bpi$. 
Let $S$ be an arbitrary quarter subset of $(\Or/(\bpi))^{\times}$ ; hence 
$(\Or/(\bpi))^{\times}=\bigcup_{\varepsilon\in W}\varepsilon S$. 
First of all, we notice the following two equations : 
\[ \prod_{\nu\in S}\nu^4\equiv \chi_{\bpi}(-1)\cdot (p-1)!
\meqv{-\chi_{\bpi}(-1)}{\bpi}, \]
\[ \prod_{\nu\in S}\varphi(\nu/\bpi)^4=
\chi_{\bpi}(-1)\prod_{\nu\in (\Or/(\bpi))^{\times}}
\varphi(\nu/\bpi)=\chi_{\bpi}(-1)\bpi. \]

\noindent
By the first equation we can define a quartic or octic root of unity 
according to each $S$ as follows. 
In the case $p\meqv{5}{8}$, i.e. $\chi_{\bpi}(-1)=-1$, we put and 
denote by $\gamma(S)$ the quartic root of unity determined by the property 
$\gamma(S)\equiv \prod_{\nu\in S} \nu \pmod{(\bpi)}$. 
In the case $p\meqv{1}{8}$, i.e. $\chi_{\bpi}(-1)=1$, 
the prime $\bpi$ is decomposable in $\Z[\zeta_8]$. Let $\Pi$ be a prime 
once chosen and fixed such that $\bpi=\Pi\,\Pi'$\ \ in $\Z[\zeta_8]$. 
We denote by $\gamma(S)$ the quartic root of $-1$ such that 
$\gamma(S)\equiv \prod_{\nu\in S} \nu \pmod{(\Pi)}$. Note 
$\gamma(S)\not\in W$, but $\zeta_8\gamma(S)\in W$ in this case. 
Unfortunately $\gamma(S)$ depends on either choice of $(\Pi)$, and the sign 
will be changed when another $(\Pi')$ is chosen, that is, 
only $\gamma(S)^2$ is an invariant of $\bpi$ and $S$. 
We should also remark that\ \ 
$\gamma(S)^2\equiv \prod_{\nu\in S} \nu^2 \pmod{(\bpi)}$\ is valid 
for both cases. 

\medskip

\noindent
{\bf Definition.}\quad The following is called the \textit{canonical quartic 
root} of\ $-\bpi$.
\begin{equation}
\widetilde{\bpi} \doteq \gamma(S)^{-1}\,\prod_{\nu\in S}
\varphi(\nu/\bpi). \label{canonic}
\end{equation}

We have\ \ $\widetilde{\bpi}^4=-\bpi$,\ \ and 
$\widetilde{\bpi}$ is independent of the choice of $S$ 
because of the property $\varphi(iu)=i\varphi(u)$. 
As is remarked in the above, there is an ambiguity of the sign of 
$\widetilde{\bpi}$ in the case $p\meqv{1}{8}$. Also we should 
remark that $\widetilde{\bpi} \not\in L=F(\varphi(1/\bpi))$ but 
$\zeta_8\,\widetilde{\bpi} \in L$, and\ $(\zeta_8\,\widetilde{\bpi})^4=\bpi$\ 
in the case $p\meqv{1}{8}$, while 
$\widetilde{\bpi} \in L$ holds in the case $p\meqv{5}{8}$.

\medskip

\noindent
{\bf Remark.}\quad In the case $p\meqv{5}{8}$, we can take 
$S=\ker \chi_{\bpi}$ as a quarter subset of $(\Or/(\bpi))^{\times}$ ; 
then $S$ is the subgroup consisting of all 
quartic residues mod $(\bpi)$. 
This choice has some advantages. Particularly, it is valid 
\[ \G_{\bpi}(\chi_{\bpi},\,f)=\sum_{\nu\in S} f(\nu/\bpi),\quad 
\gamma(S)=1, \quad \widetilde{\bpi}=\prod_{\nu\in S} \varphi(\nu/\bpi). \]

\medskip

The following is a fundamental property of the quartic residue symbol, 
which is also directly verified in view of\ 
$\gamma(\mu\,S)=\chi_{\bpi}(\mu)\,\gamma(S)$. 
\begin{equation}
\widetilde{\bpi}^{\sigma_{\mu}} = \chi_{\bpi}(\mu)\,\widetilde{\bpi}, \quad 
(\,\mu\in (\Or/(\bpi))^{\times}\,) \label{quartres}
\end{equation}
where we mean $\widetilde{\bpi}^{\sigma_{\mu}}=
\zeta_8^{-1}(\zeta_8\widetilde{\bpi})^{\sigma_{\mu}}$ in the strict meaning 
when $p\meqv{1}{8}$.

\medskip

\noindent
{\bf Definition.}\quad 
The following is called the \textit{coefficient} of the elliptic Gauss sum\ 
$\G_{\bpi}(\chi_{\bpi},\,f)$.
\begin{equation}
\alpha_{\bpi} \doteq \widetilde{\bpi}^{-3}\,
\G_{\bpi}(\chi_{\bpi},\,f). \label{coeff}
\end{equation}

\begin{Theorem}\ 
The elliptic Gauss sum is expressible as follows :
\[ \G_{\bpi}(\chi_{\bpi},\,f)=\alpha_{\bpi}\,\widetilde{\bpi}^3, \]
where the coefficient $\alpha_{\bpi}$ is an algebraic integer in 
$\Or$ or in $\zeta_8\Or$, for $p\meqv{5}{8}$ or $p\meqv{1}{8}$, respectively. 
Futher, one has
\begin{equation}
\alpha_{\bpi}\meqv{1}{(1+i)}\quad 
\textit{if}\ \ p=\bpi\,\overline{\bpi}\meqv{5}{8}.
\end{equation}
\end{Theorem}

\noindent
{\textit{Proof.}}\quad By the definition and by virtue of (\ref{lagr}) and 
(\ref{quartres}), we have \ $\alpha_{\bpi}^{\sigma_{\mu}}=\alpha_{\bpi}$\ \ 
for an arbitrary $\mu\in (\Or/(\bpi))^{\times}$, 
and hence $\alpha_{\bpi}\in F$. 
For the integrity, we can check it in similar manner to the proof (b) 
of Lemma II.6. We here assume that {\rm Claim\,($\Ze$)} is valid. 
Suppose first $p\meqv{5}{8},\ p>5$, then 
$\alpha_{\bpi}\,\widetilde{\bpi}^3$ is 
an algebraic integer, though 
$\bpi=-\widetilde{\bpi}^4$ is a prime in $\Or$; 
it means $\alpha_{\bpi}$ itself is already an integer. 
In the case $p\meqv{1}{8}$ we need some modification, but the essence is 
the very same. The last assertion of Theorem II.1 is immediately deduced 
from the lemmas II.7 and II.8.

\bigskip

\begin{Example}{\normalfont \rmfamily \ 
This example is based on an idea of Y.\,$\widehat{\rm O}$nishi. 
Consider the case $\bpi=3+2i\ (p=13)$, 
and we shall show $\alpha_{\bpi}=1$.\ \ 
Take $S=\ker \chi_{\bpi}$. 
Then $S=\{ 1,\,3,\,9\,\}=\{ 1,\,-2i,\,1-i\,\}$,\ \ and so we have 
\[ \widetilde{\bpi}=\prod_{\nu\in S} \varphi(\nu/\bpi)
=\varphi(1/\bpi)\,\varphi(-2i/\bpi)\,
\varphi((1-i)/\bpi)
=\dfrac{(-2-2i)\,\varphi(1/\bpi)^3}
{1+\varphi(1/\bpi)^4}, \]
by using suitable multiplication formulas. Since\ \ 
$-2-2i=1-\bpi=1+\widetilde{\bpi}^4$,\ \ 
$\varphi(1/\bpi)$ is a solution of the following equation.
\[ \widetilde{\bpi}\,x^4-(1+\widetilde{\bpi}^4)\,x^3 
+\widetilde{\bpi} = 0. \]
The equation is decomposed as follows.
\[ (\widetilde{\bpi}\,x-1)
(x^3-\widetilde{\bpi}^3\,x^2-\widetilde{\bpi}^2\,x 
- \widetilde{\bpi}) = 0. \]
The second factor must be the minimal polynomial of 
$\varphi(1/\bpi)$ over 
$F(\widetilde{\bpi})$,\ \ and hence\ the sum of the three roots 
$\varphi(\nu/\bpi)\ \ (\nu\in S)$\ \ is 
$\widetilde{\bpi}^3$, namely, 
\[ \G_{\bpi}(\chi_{\bpi},\,\varphi)=\sum_{\nu\in S}\varphi(\nu/\bpi)=\widetilde{\bpi}^3. \] 

In general, it seems pretty hard to compute the value of the coefficient 
$\alpha_{\bpi}$ by hand. 
More examples by computor will be given in the table in Appendix II.
}\end{Example}

\vspace{5mm}

\subsection{The quartic Hecke characters and 
{\boldmath$L$}-values at {\boldmath$s=1$}}

\noindent
{\bf 4.1.}\quad 
We introduce a Hecke character $\widetilde{\chi}_{\bpi}$ induced by 
the quartic residue character $\chi_{\bpi}$.\ \ 
As mentioned before, it is of the form\ \ $\widetilde{\chi}_{\bpi}((\nu))
=\chi_1(\nu)\,\overline{\nu}$\ \ with a residue class character $\chi_1$.\ \ 
For the purpose we first modify the character $\chi_{\bpi}$ 
into $\chi_1$ satisfying $\chi_1(i)=i$.\ \ 
After the preparation of supplementary simple characters\ 
$\chi_0$ and $\chi_0'$,\ \ we shall treat the four cases separately 
in view of\ $\chi_{\bpi}(i)=i^{(p-1)/4}$.

\medskip

Let $\chi_0$ be the character with conductor $(1+i)^3$ that gives 
the natural isomorphism\ \ 
$(\Or/(1+i)^3)^{\times}\ \cong\ W$,\ \ namely,
\[ \chi_0(\nu) \doteq \varepsilon\ \ \mbox{for}\ \ 
\nu\meqv{\varepsilon}{(1+i)^3},\ \ \varepsilon\in W=\{\pm 1,\,\pm i\}. \]

Let $\chi_0'$ be the character with conductor $(1+i)^2$ that gives 
the natural isomorphism\ \ 
$(\Or/(1+i)^2)^{\times}=(\Or/(2))^{\times}\ \cong\ \{\pm 1\}$,\ \ 
namely,
\[ \chi_0'(\nu) \doteq \delta^2\ \ \mbox{for}\ \ 
\nu\meqv{\delta}{(1+i)^2},\ \ \delta\in \{1,\,i \}. \]

Let $\bpi$ be a primary prime in $\Or$ ; $\bpi\meqv{1}{(1+i)^3}$. 
and let $\chi_{\bpi}$ be the quartic residue character to the modulus $(\bpi)$.

\medskip

\noindent
{\bf Definition.}\quad For each $\bpi$, the Hecke character 
$\widetilde{\chi}_{\bpi}$ is fixed throughout as follows.
\begin{equation}
\widetilde{\chi}_{\bpi}((\nu)) \doteq 
\chi_1(\nu)\,\overline{\nu},\quad \chi_1 \doteq 
\begin{cases}
\chi_{\bpi}\cdot \chi_0' \quad 
& \mbox{for}\quad p=\bpi\,
\overline{\bpi}\meqv{13}{16}, \\
\chi_{\bpi} \quad 
& \mbox{for}\quad p=\bpi\,
\overline{\bpi}\meqv{5}{16}, \\
\chi_{\bpi}\cdot \chi_0 \quad 
& \mbox{for}\quad p=\bpi\,
\overline{\bpi}\meqv{1}{16}, \\
\chi_{\bpi}\cdot \overline{\chi}_0 \quad 
& \mbox{for}\quad p=\bpi\,
\overline{\bpi}\meqv{9}{16}. 
\end{cases}
\end{equation}

\medskip

For later use, we summarize these circumstances as a brief list : 

\medskip

(a)\ \ \ The case\ \ $p=\bpi\,\overline{\bpi}\meqv{13}{16}$.\ \ 
The conductor of $\widetilde{\chi}_{\bpi}$ is 
$(\beta)
=(2\bpi)$.
\begin{equation} \begin{split} 
& (\Or/(\beta))^{\times} \cong (\Or/(\bpi))^{\times}\,\times\,\{\pm 1\}\ \ 
\mbox{by}\ \ \lambda \ \ \mbox{to}\ \ (\kappa,\ \delta^2) :\ 
\lambda\equiv 2\,\kappa + \bpi\,\delta\!\!\pmod{\beta}\ ; \\
& \chi_1(\lambda)=\chi_{\bpi}(2)\,\chi_{\bpi}(\kappa)\,\delta^2. 
\end{split} \label{13} \end{equation}

\medskip

(b)\ \ \ The case\ \ $p=\bpi\,
\overline{\bpi}\meqv{5}{16}$.\ \ 
The conductor of $\widetilde{\chi}_{\bpi}$ is $(\beta)
=(\bpi)$.
\begin{equation}
(\Or/(\beta))^{\times} \cong (\Or/(\bpi))^{\times}\ \ 
\mbox{by}\ \ 
\lambda\equiv \kappa\!\!\pmod{\beta}\ ;\ \ 
\chi_1(\lambda)=\chi_{\bpi}(\kappa). \label{5}
\end{equation}

\medskip

(c)\ \ \ The case\ \ $p=\bpi\,\overline{\bpi}\meqv{1}{8}$.\ \ 
The conductor of $\widetilde{\chi}_{\bpi}$ is $(\beta)=((1+i)^3\bpi)$. 
\begin{equation} \begin{split}
& (\Or/(\beta))^{\times} \cong (\Or/(\bpi))^{\times}\,\times\,W\ \ 
 \mbox{by}\ \lambda \ \mbox{to}\ \ (\kappa,\ \varepsilon) :\ 
\lambda\equiv (1+i)^3\,\kappa + \bpi\,\varepsilon
\!\!\pmod{\beta}\ ; \\
& \chi_1(\lambda)= \begin{cases} 
\overline{\chi}_{\bpi}(1+i)\,
\chi_{\bpi}(\kappa)\,\varepsilon\ \ (p\meqv{1}{16}), \\
\overline{\chi}_{\bpi}(1+i)\,
\chi_{\bpi}(\kappa)\,\overline{\varepsilon} \ \ (p\meqv{9}{16}). \end{cases} 
\end{split} \label{1} \end{equation}

\bigskip

\bigskip

\noindent
{\bf 4.2.}\quad Now we can evaluate the value of the associated\ $L$-series, 
especially at\ $s=1$,\ and we shall show that $L(1,\,\widetilde{\chi}_{\bpi})$ 
is expressed by the corresponding elliptic Gauss sum. 

\begin{Theorem}\ Let $\widetilde{\chi}_{\bpi}$ be 
the Hecke character for a primary prime\ $\bpi$. Then 
\begin{equation}
\varpi^{-1}\,L(1,\,\widetilde{\chi}_{\bpi}) = 
\begin{cases}
-\dfrac{1+i}{2} \chi_{\bpi}(2)\,\bpi^{-1}\,
\G_{\bpi}(\chi_{\bpi},\,\varphi)\ \ 
& \mbox{if}\quad p=\bpi\,
\overline{\bpi}\meqv{13}{16}, \\
\bpi^{-1}\,\G_{\bpi}(\chi_{\bpi},\,\Ze)\ \ 
& \mbox{if}\quad p=\bpi\,
\overline{\bpi}\meqv{5}{16}, \\
\overline{\chi}_{\bpi}(1+i)\,\bpi^{-1}\,
\G_{\bpi}(\chi_{\bpi},\,\psi)\ \ 
& \mbox{if}\quad p=\bpi\,
\overline{\bpi}\meqv{1}{16}, \\
-\overline{\chi}_{\bpi}(1+i)\,\bpi^{-1}\,
\G_{\bpi}(\chi_{\bpi},\,\psi)\ \ 
& \mbox{if}\quad p=\bpi\,
\overline{\bpi}\meqv{9}{16}. 
\end{cases} \end{equation}
\end{Theorem}

\noindent
{\textit{Proof.}}\quad We follow the formula (\ref{EGSg}) of Lemma II.5. 

\medskip

(a)\ \ \ The case $p=\bpi\,\overline{\bpi}
\meqv{13}{16}$.\ \ In view of (\ref{13}), we have
\begin{eqnarray*}
\varpi^{-1}\,L(1,\ \widetilde{\chi}_{\bpi}) & = & 
\dfrac{1}{4\,\beta}\cdot \sum_{\lambda\,({\rm mod}\,\beta)}\,
\chi_1(\lambda)\,\Ze(\lambda/\beta) \\
&=& \dfrac{1}{8\,\bpi}\,\chi_{\bpi}(2)\cdot
\sum_{\kappa\,({\rm mod}\,\bpi)} \chi_{\bpi}(\kappa)
\sum_{\delta=1,\,i} \delta^2\,\Ze(\kappa/\bpi+\delta/2) \\
{} &=& -\dfrac{1+i}{2}\cdot 
\dfrac{\chi_{\bpi}(2)}{\bpi} \cdot 
\dfrac{1}{4} \sum_{\kappa\,({\rm mod}\,\bpi)} 
\chi_{\bpi}(\kappa)\varphi(\kappa/\bpi),
\end{eqnarray*}
since\ \ 
$\sum_{\delta=1,\,i} \delta^2\,\Ze(u + \delta/2)=-(1+i)\,\varphi(u)$\ \ 
by the definition (\ref{phi}). 

\medskip

(b)\ \ \ The case $p=\bpi\,\overline{\bpi}
\meqv{5}{16}$.\ \ In view of (\ref{5}), we obtain directly 
\begin{equation*}
\varpi^{-1}\,L(1,\ \widetilde{\chi}_{\bpi}) =
\dfrac{1}{\bpi}\cdot \dfrac{1}{4}
\sum_{\kappa\,({\rm mod}\,\bpi)}\,
\chi_{\bpi}(\kappa)\,\Ze(\kappa/\bpi).
\end{equation*}

\medskip

(c)\ \ \ The case $p=\bpi\,\overline{\bpi}
\meqv{1}{16}$.\ \ 
In view of (\ref{1}), we have
\begin{eqnarray*}
\varpi^{-1}\,L(1,\ \widetilde{\chi}_{\bpi}) & = & 
-\dfrac{1+i}{16}\cdot 
\dfrac{\overline{\chi}_{\bpi}(1+i)}{\bpi}
\sum_{\kappa\,({\rm mod}\,\bpi)}\,
\chi_{\bpi}(\kappa)
\sum_{\varepsilon\in W}\varepsilon\,
\Ze(\kappa/\bpi - \varepsilon(1+i)/4) \\
&=& \dfrac{\overline{\chi}_{\bpi}(1+i)}{\bpi}
\cdot \dfrac{1}{4}\sum_{\kappa\,({\rm mod}\,\bpi)}\,
\chi_{\bpi}(\kappa)\,\psi(\kappa/\bpi),
\end{eqnarray*}
since\ \ $\sum_{\varepsilon\in W}\varepsilon\,\Ze(u-\varepsilon(1+i)/4)
=-(1-i)\{\psi(u)+\psi(iu)\}$\ \ by the definition (\ref{psi}).

\medskip

(c')\ \ \ The case $p=\bpi\,\overline{\bpi}
\meqv{9}{16}$.\ \ 
In view of (\ref{1}), we have
\begin{eqnarray*}
\varpi^{-1}\,L(1,\ \widetilde{\chi}_{\bpi}) & = & 
-\dfrac{1+i}{16}\cdot 
\dfrac{\overline{\chi}_{\bpi}(1+i)}{\bpi}
\sum_{\kappa\,({\rm mod}\,\bpi)}\,
\chi_{\bpi}(\kappa)
\sum_{\varepsilon\in W}\,\overline{\varepsilon}\,
\Ze(\kappa/\bpi - \varepsilon(1+i)/4) \\
&=& -\dfrac{\overline{\chi}_{\bpi}(1+i)}{\bpi}
\cdot \dfrac{1}{4}\sum_{\kappa\,({\rm mod}\,\bpi)}\,
\chi_{\bpi}(\kappa)\,\psi(\kappa/\bpi),
\end{eqnarray*}
since\ \ $\sum_{\varepsilon\in W}\,\overline{\varepsilon}\,
\Ze(u-\varepsilon(1+i)/4)
=(1-i)\{\psi(u)-\psi(iu)\}$\ \ also by the definition (\ref{psi}).

\medskip

Thus we have completed the proof of Theorem II.2.

\vspace{5mm}

\subsection{The explicit formula of the root number\ \ 
{\boldmath$C(\widetilde{\chi}_{\bpi})$}}

\noindent
{\bf 5.1.}\quad 
We require an important formula about 
the classical quartic Gauss sum. 
Let $\bpi$ be a primary prime in $\Or$ 
and set $p=\bpi\,\overline{\bpi}$.\ \ 
The quartic residue character $\chi_{\bpi}$ may be considered as 
a character on $(\Z/p\Z)^{\times}$.\ \ 
Then the quartic Gauss sum is defined by 
\begin{equation}
G_4(\bpi) 
\doteq \sum_{r=1}^{p-1} \chi_{\bpi}(r)\,e^{2\pi ir/p}.
\end{equation}

Also we here should recall the definition (\ref{canonic}) 
of the canonical quartic root $\widetilde{\bpi}$ 
of $-\bpi$.
\begin{Lemma}[$G_4$-formula]\ \ 
\begin{equation}
G_4(\bpi) = \chi_{\bpi}(-2)\,
\widetilde{\bpi}^3\,\overline{\widetilde{\bpi}} 
\label{G4}
\end{equation}
\end{Lemma}

\noindent
{\bf Remark.}\ \ As an immediate consequence we have a famous formula 
(cf. [I-R], Prop.9.10.1):\ 
$G_4(\bpi)^2=-\chi_{\bpi}(-1)\,\bpi\,
\sqrt{p}$,\ \ and also we obtain\ \ $G_4(\bpi)^4
=\bpi^3\,\overline{\bpi}$.\ \ 
The ambiguity of the definition of 
$\widetilde{\bpi}$ does not matter, for the right side 
of our $G_4$-formula depends only on $\widetilde{\bpi}^2$. 

\medskip

\noindent
{\textit{Proof.}}\ \ This is only a slight modification of the celebrated 
formula of Matthews. He used the lattice\ $\theta\,\Or$\ 
instead of our\ $\varpi\,\Or$,\ 
where\ \ $\theta=\sqrt{2}\,\varpi=3.70814935\cdots$.\ \ Let\ $\wp_1(u)$ 
denote Weierstrass' $\wp$ with the period lattice\ $\theta\,\Or$. 
Hence the relation $\wp(\varpi\,u)=2\,\wp_1(\theta\,u)$ holds, 
so that $\wp_1'(u)^2=4\,\wp_1(u)^3-\wp_1(u)$\ ;\ further we have
\[ \varphi(u)={\rm sl}((1-i)\,\varpi\,u)
=-2\,(1-i)\,\dfrac{\wp(\varpi\,u)}{\wp'(\varpi\,u)}
=-\sqrt{2}\,(1-i)\,\dfrac{\wp_1(\theta\,u)}{\wp_1'(\theta\,u)}
=\zeta_8^{-1}\,T(\theta\,u), \]
where $T(u)=-2\,\wp_1(u)\wp_1'(u)^{-1}$ after his notation.\ \ 
Matthews' formula states 

\medskip

\noindent
{\bf Formula} (\,[M2], esp. p.51\,)\ \ 
\begin{equation}
G_4(\bpi)=-\beta(\bpi)\,
\chi_{\bpi}(2i) 
\prod_{r\in N} T(\theta\,r/\bpi)\cdot p^{1/4}, \label{Matthews}
\end{equation}
where\ $N=\{1,\,2,\,\cdots,\,(p-1)/2\}$, 
and the constant $\beta(\bpi)$ is uniquely determined 
by the conditions\ $\beta(\bpi)\equiv 
\prod_{r\in N}\,r\pmod{\bpi}$\ and $\beta(\bpi)^2=-1$. 

\medskip

Now we can derive the formula (\ref{G4}) from (\ref{Matthews}); 
indeed they are equivalent. We first 
note that since $N$ is a half subset mod $(\bpi)$, 
there is a quarter subset $S_0$ such that $N=S_0\cup iS_0$. 
Then $\beta(\bpi)\equiv \prod_{r\in N}\,r\equiv 
\chi_{\bpi}(i)\prod_{r\in S_0}\,r^2\pmod{(\bpi)}$ ; 
namely, $\beta(\bpi)=\chi_{\bpi}(i)\,\gamma(S_0)^2$.\ \ 
Also we can observe $\gamma(S_0)^4=-\chi_{\bpi}(-1)$. Hence we have 
\begin{eqnarray*}
G_4(\bpi) &=& -\chi_{\bpi}(-2)\,\gamma(S_0)^2
\prod_{r\in N}\{\zeta_8\,\varphi(r/\bpi)\}\cdot p^{1/4} \\
&=& -\chi_{\bpi}(-2)\,\gamma(S_0)^4\,
\chi_{\bpi}(-1)\cdot \gamma(S_0)^{-2}\prod_{\nu\in S_0}\,
\varphi(\nu/\bpi)^2 \cdot p^{1/4} 
= \chi_{\bpi}(-2)\,\widetilde{\bpi}^3\,
\overline{\widetilde{\bpi}}.
\end{eqnarray*}
This finishes the proof of $G_4$-formula.

\bigskip

\bigskip

\noindent
{\bf 5.2.}\quad We are now ready to give the explicit value of the root number 
$C(\widetilde{\chi}_{\bpi})$. 
\begin{Theorem}\ Let $\widetilde{\chi}_{\bpi}$ be 
the Hecke character for a primary prime $\bpi$. Then 
\begin{equation}
C(\widetilde{\chi}_{\bpi})=
\begin{cases}
-i\,\widetilde{\bpi}^{-1}\,
\overline{\widetilde{\bpi}}\quad\ \ \ & \mbox{if}\quad 
p=\bpi\,\overline{\bpi}\meqv{13}{16}, \\
i\,\chi_{\bpi}(2)\,\widetilde{\bpi}^{-1}\,
\overline{\widetilde{\bpi}}\quad\ \ \ & \mbox{if}\quad 
p=\bpi\,\overline{\bpi}\meqv{5}{16}, \\
-\overline{\chi}_{\bpi}(1+i)\,
\widetilde{\bpi}^{-1}\,
\overline{\widetilde{\bpi}}\quad\ \ \ & \mbox{if}\quad 
p=\bpi\,\overline{\bpi}\meqv{1}{16}, \\
-i\,\chi_{\bpi}(1+i)\,\widetilde{\bpi}^{-1}\,
\overline{\widetilde{\bpi}}\quad\ & \mbox{if}\quad 
p=\bpi\,\overline{\bpi}\meqv{9}{16}.
\end{cases} \end{equation}
\end{Theorem}

\noindent
\textit{Proof.}\ \ We first evaluate some simple Gauss sums. 
The first three are easily verified by direct calculation :
\begin{equation} \begin{split}
g(\chi_0)\doteq 
\sum_{\varepsilon\in W} \varepsilon\,e^{2\pi i\,{\rm Re}(\varepsilon/(1+i)^3)}
&=-2-2i,\ \ g(\overline{\chi}_0)\doteq \sum_{\varepsilon\in W} 
\overline{\varepsilon}\,e^{2\pi i\,{\rm Re}(\varepsilon/(1+i)^3)}
=2-2i \\
\mbox{\textit{and}}\quad & 
g(\chi_0')\doteq 
\sum_{\delta=1,\,i} \delta^2\,e^{2\pi i\,{\rm Re}(\delta/2)}=-2.
\end{split} \end{equation}
The next sum is essentially nothing but the quartic Gauss sum :
\begin{equation}
g(\chi_{\bpi}) \doteq \sum_{\kappa\,({\rm mod}\,\bpi)} 
\chi_{\bpi}(\kappa)\,e^{2\pi i\,{\rm Re}(\kappa/\bpi)}
=\chi_{\bpi}(-2)\,\widetilde{\bpi}^3\,
\overline{\widetilde{\bpi}}\cdot 
\begin{cases}
1 & \textit{if}\quad p\meqv{13,\,1}{16} \\
-1 & \textit{if}\quad p\meqv{5,\,9}{16} \end{cases}
 \label{gchi}
\end{equation}
In fact, we first replace the sum over\ 
$\kappa\,({\rm mod}\,\bpi)$\ 
by one over\ $r\,({\rm mod}\,p)$,\ and then, by using\ 
${\rm Re}(r\,\overline{\bpi}/p)=ar/p$\ 
where $\bpi=a+b\,i\ (a,\,b\in \Z)$, 
we can calculate as follows :
\[ g(\chi_{\bpi})=\sum_{r\,({\rm mod}\,p)} 
\chi_{\bpi}(r)\,
e^{2\pi i\,{\rm Re}(r\,\overline{\bpi}/p)}
=\sum_{r\,({\rm mod}\,p)} \chi_{\bpi}(r)\,e^{2\pi i a r/p} 
= \overline{\chi}_{\bpi}(a)\,G_4(\bpi). \]
Furthermore, we know $\overline{\chi}_{\bpi}(a)=1$\ or\ $-1$\ 
for\ $p\equiv {13,\,1}$\ or ${5,\,9}\pmod{16}$, respectively 
(cf. [IR] Chap.\,9, Exerc.\,34), and finally, 
by applying $G_4$-formula we have (\ref{gchi}).

We return to the proof of Theorem II.3. Using 
the formula (\ref{the_root}) of Lemma II.4 :
\[ C(\widetilde{\chi}_{\bpi}) = -i\,\beta^{-1}\,
\sum_{\lambda\,({\rm mod}\,\beta)}\,
\chi_1(\lambda)\,e^{2\pi i\,{\rm Re}(\lambda /\beta)}, \]
we treat each of the four cases according to the definition of $\chi_1$, 
especially in view of (\ref{13}), (\ref{5}) and (\ref{1}).\ \ 
Thus we can easily obtain 

\medskip

(a)\ \ The case\ \ $p=\bpi\,\overline{\bpi}
\meqv{13}{16}$. 
\[ C(\widetilde{\chi}_{\bpi}) = 
-i\,2^{-1}\,\bpi^{-1}\,\chi_{\bpi}(2)\,
g(\chi_0')\,g(\chi_{\bpi})=-i\,\widetilde{\bpi}^{-1}\,
\overline{\widetilde{\bpi}}. \]

(b)\ \ The case\ \ $p=\bpi\,\overline{\bpi}
\meqv{5}{16}$. 
\[ C(\widetilde{\chi}_{\bpi}) = 
-i\,\bpi^{-1}\,g(\chi_{\bpi})=
i\,\chi_{\bpi}(2)\cdot \widetilde{\bpi}^{-1}\,
\overline{\widetilde{\bpi}}. \]

(c)\ \ The case\ \ $p=\bpi\,\overline{\bpi}
\meqv{1}{16}$. 
\[ C(\widetilde{\chi}_{\bpi}) = 
-i\,(1+i)^{-3}\,\bpi^{-1}\,
\overline{\chi}_{\bpi}(1+i)\,
g(\chi_0)\,g(\chi_{\bpi})=
-\chi_{\bpi}(1+i)\,\widetilde{\bpi}^{-1}\,
\overline{\widetilde{\bpi}}. \]

(c')\ \ The case\ \ $p=\bpi\,\overline{\bpi}
\meqv{9}{16}$. 
\[ C(\widetilde{\chi}_{\bpi}) = 
-i\,(1+i)^{-3}\,\bpi^{-1}\,
\overline{\chi}_{\bpi}(1+i)\,
g(\overline{\chi}_0)\,g(\chi_{\bpi})=
-i\,\chi_{\bpi}(1+i)\,\widetilde{\bpi}^{-1}\,
\overline{\widetilde{\bpi}}. \]

These complete the proof of Theorem II.3.

\vspace{5mm}

\subsection{Rationality of the elliptic Gauss sum coefficient}

\noindent
{\bf 6.1.}\quad Now we can mention about the rationality 
of the coefficients of elliptic Gauss sums. More precisely, 
the coefficient itself is not always rational 
but it will be seen that the essential factor of this 
is certainly a rational integer. In other words, we shall extract 
a rational integer factor from the elliptic Gauss sum, i.e. 
from the Hecke $L$-value at $s=1$, 
which seems also the most important part in respect of 
arithmetical nature. The following theorem, together with 
corollaries, is the main result of the part {\bf II}. 
\begin{Theorem}\ 
Let $\alpha_{\bpi}$ be the coefficient 
of the elliptic Gauss sum {\rm (\ref{coeff})}. Then 
\begin{equation}
\alpha_{\bpi}= \begin{cases}
\overline{\alpha}_{\bpi}
 & \mbox{if}\quad p
=\bpi\,\overline{\bpi}\meqv{13}{16}, \\
i\,\chi_{\bpi}(2)\,
\overline{\alpha}_{\bpi} & \mbox{if}\quad p
=\bpi\,\overline{\bpi}\meqv{5}{16}, \\
-\,\overline{\chi}_{\bpi}(1+i)\,
\overline{\alpha}_{\bpi} & \mbox{if}\quad p
=\bpi\,\overline{\bpi}\meqv{1}{16}, \\
-i\,\overline{\chi}_{\bpi}(1+i)\,
\overline{\alpha}_{\bpi} & \mbox{if}\quad p
=\bpi\,\overline{\bpi}\meqv{9}{16}.
\end{cases} \label{ratio} \end{equation}
\end{Theorem}

\noindent
\textit{Proof.}\ \ By the theorems II.1,\,II.2 and II.3 we have already known 
both the explicit values of $L(1,\,\widetilde{\chi}_{\bpi})$ 
and $C(\widetilde{\chi}_{\bpi})$. 
To prove Theorem II.4, we have only to substitute them for the both sides of 
the central value equation\ :\ 
$L(1,\,\widetilde{\chi}_{\bpi})
=C(\widetilde{\chi}_{\bpi})\,
\overline{L(1,\,\widetilde{\chi}_{\bpi})}$\ \ (cf. Lemma II.4).

\noindent
For example, suppose that $p=\bpi\,\overline{\bpi}\meqv{13}{16}$.\ \ 
In this case we have
\[ \varpi^{-1}\,L(1,\,\widetilde{\chi}_{\bpi})
=\dfrac{1+i}{2}\,\chi_{\bpi}(2)\,\widetilde{\bpi}^{-1}\,\alpha_{\bpi}\quad \mbox{and}\quad 
C(\widetilde{\chi}_{\bpi})=-i\,\widetilde{\bpi}^{-1}\,
\overline{\widetilde{\bpi}}, \]
and hence the central value equation 
implies\ \ $\chi_{\bpi}(2)\,\alpha_{\bpi}=
-\overline{\chi}_{\bpi}(2)\,
\overline{\alpha}_{\bpi}$.\ \ This immediately proves 
$\alpha_{\bpi}=\overline{\alpha}_{\bpi}$\ \ 
since $\chi_{\bpi}(2)^2=\chi_{\bpi}(i)^2=-1$.\ \ 
In other cases, the argument is quite similar, so we omit the details. 
Thus the proof is finished.

\bigskip

Before stating the corollaries of Theorem II.4, we give a proof of the 
integrity of $\G_{\bpi}(\chi_{\bpi},\,\Ze)$, 
which has been postponed until now. 

\medskip

\noindent
\textit{Proof of} {\bf Claim ($\Ze$)}.\ \ Assume that 
$p=\bpi\,\overline{\bpi}\meqv{5}{16}$ and 
$\G_{\bpi}(\chi_{\bpi},\,\Ze)=
\alpha_{\bpi}\,\widetilde{\bpi}^3$.\ \ 
We know $\alpha_{\bpi} \in F$. Further, Theorem II.4 shows 
$\overline{\chi}_{\bpi}(1+i)\,\alpha_{\bpi}
=\chi_{\bpi}(1+i)\,\overline{\alpha}_{\bpi}$, and hence\ 
$\alpha_{\bpi}=\chi_{\bpi}(1+i)\, a_{\bpi}$\ for some $a_{\bpi} \in \Q$. 
On the other hand, by Lemma II.7 and the subsequent discussions, we know 
$((1+i)-i\,\overline{\chi}_{\bpi}(1+i))\,
\alpha_{\bpi} \in \Or$. Namely, 
$a_{\bpi},\,-(1+2i)\,a_{\bpi},\,-a_{\bpi}$ or $(1-2i)\,a_{\bpi}$ 
is algebraic integer in $\Or$, when $\chi_{\bpi}(1+i)=1,\,-1,\,i$ or 
$-i$, accordingly. This means that $a_{\bpi} \in \Z$ holds already. 
Thus the proof of Claim $(\Ze)$ is completed.

\begin{Corollary}\ Suppose that $p=\bpi\,\overline{\bpi}\meqv{5}{8}$. 
There exists a rational integer $a_{\bpi}$ such that 
$a_{\bpi}\meqv{1}{2}$, and the coefficient 
$\alpha_{\bpi}$ of the elliptic Gauss sum is expressed by $a_{\bpi}$ 
as follows. In particular,\ \ $|\alpha_{\bpi}|^2=a_{\bpi}^2$.
\begin{equation}
\alpha_{\bpi}= \begin{cases}
a_{\bpi} & \mbox{if} \quad p=\bpi\,\overline{\bpi}\meqv{13}{16}, \\
a_{\bpi}\,\chi_{\bpi}(1+i) & 
\mbox{if} \quad p=\bpi\,\overline{\bpi}\meqv{5}{16}.
\end{cases} \end{equation}
\end{Corollary}

\noindent
\textit{Proof.}\ \ 
We can derive the integrity of the coefficient\ 
$\alpha_{\bpi}$\ from Theorem II.1, and the rationality from Theorem II.4. 
The property\ $a_{\bpi}\meqv{1}{2}$\ follows from\ 
$\alpha_{\bpi}\meqv{1}{(1+i)}$\ in Theorem II.1.

\begin{Corollary}\ Suppose that $p=\bpi\,\overline{\bpi}\meqv{1}{8}$. 
There exists a rational integer $a_{\bpi}$, and the coefficient 
$\alpha_{\bpi}$ of the elliptic Gauss sum is expressed by $a_{\bpi}$ 
as follows. In particular,\ \ 
$|\alpha_{\bpi}|^2=2\,a_{\pi}^2$\ \ or\ \ $|\alpha_{\bpi}|^2=a_{\bpi}^2$\ \ 
according as\ \ $\chi_{\bpi}(2)=1$\ or\ $\chi_{\bpi}(2)=-1$.
\begin{eqnarray}
\alpha_{\bpi}= \begin{cases}
a_{\bpi}\cdot i\,\sqrt{2} & \mbox{if} \quad \chi_{\bpi}(1+i)=1, \\
a_{\bpi}\cdot \sqrt{2} & \mbox{if} \quad \chi_{\bpi}(1+i)=-1, \\
a_{\bpi}\cdot \zeta_8 & \mbox{if} \quad \chi_{\bpi}(1+i)=i, \\
a_{\bpi}\cdot i\,\zeta_8 & \mbox{if} \quad \chi_{\bpi}(1+i)=-i, 
\end{cases} \quad \mbox{and}\ \ p\meqv{1}{16}, \\
\alpha_{\bpi}= \begin{cases}
a_{\bpi}\cdot i\,\zeta_8 & \mbox{if} \quad \chi_{\bpi}(1+i)=1, \\
a_{\bpi}\cdot \zeta_8 & \mbox{if} \quad \chi_{\bpi}(1+i)=-1, \\
a_{\bpi}\cdot i\,\sqrt{2} & \mbox{if} \quad \chi_{\bpi}(1+i)=i, \\
a_{\bpi}\cdot \sqrt{2} & \mbox{if} \quad \chi_{\bpi}(1+i)=-i, 
\end{cases} \quad \mbox{and}\ \ p\meqv{9}{16}. 
\end{eqnarray}
\end{Corollary}

\noindent
\textit{Proof.}\ \ In this case we have $\alpha_{\bpi} \in \zeta_8\,\Or$, 
which combined with the rationality relation (\ref{ratio}) will immediately 
give an explicit form of the coefficient. For example, consider the case of 
$p\meqv{1}{16}$ and $\chi_{\bpi}(1+i)=1$. Put $\alpha_{\bpi}=(c+di)\,\zeta_8$ 
with $c,\,d\in \Z$. By (\ref{ratio}) we see $\alpha_{\bpi}
=-\overline{\alpha}_{\bpi}$, and hence $c+di=i(c-di)$, which means $c=d$. 
Namely, we have\ \ 
$\alpha_{\bpi}=c\,(1+i)\,\zeta_8=a_{\bpi}\cdot i\,\sqrt{2}$\ \ by 
putting\ \ $a_{\bpi}=c$.\ \ 
We omit the details for other seven cases.

\bigskip

\bigskip

\noindent
{\bf 6.2.}\ \ The substance of Theorem II.4 and the corollaries 
can be stated by the language of Hecke $L$-values in various ways. 
The following is one of them. It shows that there is a close relation 
between the values of Hecke's $L(1,\,\widetilde{\chi}_{\bpi})$ and 
the quartic Gauss sum $G_4(\bpi)$, especially between their arguments. 
Roughly speaking, the argument of $L(1,\,\widetilde{\chi}_{\bpi})$ is 
parallel to one of $\widetilde{\bpi}^{-1}$, and the argument of 
$G_4(\bpi)$ is parallel to one of $\widetilde{\bpi}^2$. Hence by 
eliminating the factors of $\widetilde{\bpi}$ from their formulas, 
the result can be obtained. 
\begin{Theorem}\ 
Let $a_{\bpi}\ (p=\bpi\,\overline{\bpi})$ be a rational integer 
given in Corollary II.1 or II.2. 
\begin{equation*}
\varpi^{-2}\,L(1,\,\widetilde{\chi}_{\bpi})^2=
\begin{cases}
2^{-1}\,i\,\chi_{\bpi}(2)\,p^{1/4}\,a_{\bpi}^2\cdot G_4(\bpi)^{-1} 
 & \textit{if}\quad p\meqv{13}{16}, \\
i\,p^{1/4}\,a_{\bpi}^2\cdot G_4(\bpi)^{-1} 
 & \textit{if}\quad p\meqv{5}{16}, \\
-2\,\overline{\chi}_{\bpi}(1+i)\,p^{1/4}\,a_{\bpi}^2\cdot G_4(\bpi)^{-1} 
 & \textit{if}\quad p\meqv{1}{16},\ \chi_{\bpi}(2)=1, \\
-\overline{\chi}_{\bpi}(1+i)\,p^{1/4}\,a_{\bpi}^2\cdot G_4(\bpi)^{-1} 
 & \textit{if}\quad p\meqv{1}{16},\ \chi_{\bpi}(2)=-1, \\
2\,i\,\overline{\chi}_{\bpi}(1+i)\,p^{1/4}\,a_{\bpi}^2\cdot G_4(\bpi)^{-1} 
 & \textit{if}\quad p\meqv{9}{16},\ \chi_{\bpi}(2)=1, \\
i\,\overline{\chi}_{\bpi}(1+i)\,p^{1/4}\,a_{\bpi}^2\cdot G_4(\bpi)^{-1} 
 & \textit{if}\quad p\meqv{9}{16},\ \chi_{\bpi}(2)=-1.
\end{cases} \end{equation*}
\end{Theorem}

\noindent
{\textit{Proof.}}\ \ Consider the case $p=\bpi\,\overline{\bpi}\meqv{13}{16}$. 
By Theorems II.1 and II.2, 
we have $\varpi^{-1}\,L(1,\,\widetilde{\chi}_{\bpi})
=2^{-1}\,(1+i)\,\chi_{\bpi}(2)\,\alpha_{\bpi}\,\widetilde{\bpi}^{-1}$. 
Theorem II.4 shows $\alpha_{\bpi}^2=a_{\bpi}^2$. These combined with 
$G_4$-formula : $G_4(\bpi)=-\chi_{\bpi}(2)\,\widetilde{\bpi}^2\,p^{1/4}$ 
implies the result. We omit the similar discussions for other cases. 
Obviously the formula for 
$|L(1,\,\widetilde{\chi}_{\bpi})|^2$ is very simple.

\begin{Corollary}\ $L(1,\,\widetilde{\chi}_{\bpi})\ne 0$\quad 
\textit{if}\quad $p=\bpi\,
\overline{\bpi}\meqv{5}{8}$.
\end{Corollary}

\noindent
{\textit{Proof.}}\ \ Because $|\alpha_{\bpi}|^2=a_{\bpi}^2\meqv{1}{2}$. 

\bigskip

\noindent
{\bf Remark.}\quad The case $p=\bpi\,\overline{\bpi}\meqv{1}{8}$. 
According to some observation 
it is plausible that $a_{\bpi}\meqv{1}{2}$ and so 
$L(1,\,\widetilde{\chi}_{\bpi})$ never vanishes 
if\ \ $p\meqv{1}{16}$\ and\ $\chi_{\bpi}(1+i)\ne 1$, or if\ \ 
$p\meqv{9}{16}$\ and\ $\chi_{\bpi}(1+i)\ne -i$. 
On the other hand, we can observe that 
$L(1,\,\widetilde{\chi}_{\bpi})$ happens to vanish very often 
in the contrary cases. 
For examples it seems \ $L(1,\,\widetilde{\chi}_{\bpi})=0$\  holds 
for each prime as follows : 
\begin{eqnarray*}
&& p=113,\,257,\,593,\,1201,\,1217,\,2129,\,2593,\,\ldots,\ \ 
( p\meqv{1}{16},\ \ \chi_{\bpi}(1+i)=1 ), \\
&& p=89,\,601,\,1097,\,1193,\,1433,\,1481,\,1721,\,\ldots,\ \ 
( p\meqv{9}{16},\ \ \chi_{\bpi}(1+i)=-i ). 
\end{eqnarray*}

\bigskip

\bigskip

\noindent
{\bf Appendix II.}\ \ 

\medskip

For convenience' and interest's sake, we append a small table 
of the coefficients of elliptic Gauss sums in the following pages.
The computation and the table was made by using {\sc ubasic}. 

In the table (1),\,(2), the coefficient is given as\ \ 
$\alpha_{\bpi}=a_{\bpi}$\ \ or\ \ $\alpha_{\bpi}=
a_{\bpi}\cdot (\chi_{\bpi}(1+i))$,\ \ for the case\ \ 
$p\meqv{13}{16}$\ or\ \ $p\meqv{5}{16}$,\ \ respectively. (cf. Corollary II.1)

For the case of $p\meqv{1}{8}$,\ we need to add a few remarks. 
In this case, as was mentioned before, 
there is an ambiguity caused by the choice of $\Pi$ or $\Pi'$ 
in defining the quartic root $\widetilde{\bpi}$ of $-\bpi$. 
In our computation, we take the quarter subset $S_0$ 
such as $S_0 \cup i\,S_0=\{1,\,2,\,\ldots,\,(p-1)/2\}$. 
And we set\ \ $\gamma(S_0)=\zeta_8$\ or\ $\overline{\zeta}_8$\ 
when\ $\gamma(S_0)^2\equiv i$\ or\ $-i\pmod{\bpi}$, respectively. 
This means that we have chosen an appropriate $\Pi$ for temporary convenience. 
Anyway in the case $p\meqv{1}{8}$,\ 
only the quantity $\alpha_{\bpi}^2$ has an invariant meaning. 
In the table (3),\,(4), we can see and check that the value of the coefficient
 is exactly in conformity with the statement of Corollary II.2. 

It is very notable that the magnitude of $a_{\bpi}$ seems to be 
remarkable small. 
In fact, thanks to Mr.\,Naruo Kanou's computation by {\sc pari/gp}, we know 
\begin{eqnarray*}
&&-49\le a_{\bpi} \le 49\quad \mbox{for}\quad 13\le p \le 3999949,\ \ 
p\meqv{13}{16}, \\
&&-43\le a_{\bpi} \le 47\quad \mbox{for}\quad 37\le p \le 3999893,\ \ 
p\meqv{5}{16}.
\end{eqnarray*}

\newpage

\small
\noindent
{} \hfill
{\bf Appendix II.\ \ Table of the coefficients of elliptic Gauss sums (1)} 
\hfill {}

\vfill

\scriptsize
\[ p=\boldsymbol{\pi} \overline{\boldsymbol{\pi}}\meqv{13}{16},\quad 
\G_{\bpi}(\chi_{\bpi},\,\varphi)=\alpha_{\bpi}\,\widetilde{\bpi}^3  \]
\[ \begin{array}{rrl|crl|crl}
p  & \boldsymbol{\pi}\quad  & \ \ 
 \alpha_{\boldsymbol{\pi}} \quad & 
p & \boldsymbol{\pi}\quad \ \ & \ \ 
 \alpha_{\boldsymbol{\pi}} \quad & 
p & \boldsymbol{\pi}\quad \ \ & \ \ 
 \alpha_{\boldsymbol{\pi}} \quad \\
\hline
 13 & 3+2i &\ \ \, 1 & 2477 & 19+46i &\ \ \, 5 & 5741 &-29+70i &-1  \\
 29 &-5+2i &\ \ \, 1 & 2557 &-21+46i &\ \ \, 5 & 5821 & 75+14i &-5  \\
 61 &-5+6i &-1 & 2621 & 11+50i &\ \ \, 3 & 5869 &-45+62i &-1  \\
 109 & 3+10i &\ \ \, 1 & 2749 & 43+30i &-1 & 5981 & 59+50i &\ \ \, 1  \\
 157 & 11+6i &\ \ \, 1 & 2797 & 51+14i &\ \ \, 1 & 6029 &-77+10i &\ \ \, 1  \\
 173 &-13+2i &-1 & 2861 & 19+50i &-1 & 6173 &-53+58i &-1  \\
 269 &-13+10i &\ \ \, 3 & 2909 &-53+10i &-5 & 6221 &-61+50i &\ \ \, 1  \\
 317 & 11+14i &-1 & 2957 &-29+46i &\ \ \, 1 & 6269 &-37+70i &\ \ \, 1  \\
 349 &-5+18i &\ \ \, 1 & 3037 & 11+54i &-1 & 6301 & 75+26i &-5  \\
 397 & 19+6i &-1 & 3181 &-45+34i &-5 & 6317 &-29+74i &-1  \\
 461 & 19+10i &-1 & 3229 & 27+50i &\ \ \, 7 & 6397 & 59+54i &\ \ \, 5  \\
 509 &-5+22i &-1 & 3373 & 3+58i &\ \ \, 1 & 6637 &-61+54i &\ \ \, 1  \\
 541 &-21+10i &\ \ \, 1 & 3389 &-5+58i &-5 & 6653 &-53+62i &-1  \\
 557 & 19+14i &-1 & 3469 &-45+38i &\ \ \, 5 & 6701 & 35+74i &-5  \\
 653 &-13+22i &-1 & 3517 & 59+6i &-3 & 6733 & 3+82i &-5  \\
 701 &-5+26i &\ \ \, 1 & 3533 &-13+58i &-1 & 6781 & 75+34i &-5  \\
 733 & 27+2i &\ \ \, 1 & 3581 & 59+10i &\ \ \, 1 & 6829 &-77+30i &-1  \\
 797 & 11+26i &-1 & 3613 & 43+42i &\ \ \, 5 & 7069 & 75+38i &\ \ \, 3  \\
 829 & 27+10i &-5 & 3677 & 59+14i &\ \ \, 3 & 7213 & 83+18i &\ \ \, 5  \\
 877 &-29+6i &-1 & 3709 &-53+30i &-3 & 7229 &-85+2i &-1  \\
 941 &-29+10i &\ \ \, 3 & 3821 &-61+10i &\ \ \, 1 & 7309 & 35+78i &\ \ \, 1  \\
 1021 & 11+30i &-3 & 3853 & 3+62i &\ \ \, 1 & 7517 & 11+86i &-1  \\
 1069 &-13+30i &-1 & 3917 &-61+14i &-1 & 7549 &-85+18i &-1  \\
 1117 &-21+26i &-5 & 4013 &-13+62i &-3 & 7741 & 75+46i &-5  \\
 1181 &-5+34i &\ \ \, 1 & 4093 & 27+58i &-3 & 7757 & 19+86i &-5  \\
 1213 & 27+22i &-3 & 4157 & 59+26i &-1 & 7789 & 83+30i &-1  \\
 1229 & 35+2i &\ \ \, 3 & 4253 &-53+38i &\ \ \, 1 & 7853 & 67+58i &-1  \\
 1277 & 11+34i &-1 & 4349 & 43+50i &-1 & 7901 &-85+26i &-1  \\
 1373 &-37+2i &-3 & 4397 &-61+26i &\ \ \, 3 & 7933 & 43+78i &\ \ \, 3  \\
 1453 & 3+38i &-3 & 4493 & 67+2i &-3 & 7949 & 35+82i &-1  \\
 1549 & 35+18i &\ \ \, 1 & 4621 &-61+30i &-1 & 8093 &-37+82i &\ \ \, 3  \\
 1597 &-21+34i &-5 & 4637 & 59+34i &-1 & 8221 & 11+90i &-1  \\
 1613 &-13+38i &-1 & 4733 &-37+58i &\ \ \, 3 & 8237 &-29+86i &\ \ \, 1  \\
 1693 &-37+18i &-1 & 4813 & 67+18i &\ \ \, 1 & 8269 &-13+90i &-7  \\
 1709 & 35+22i &\ \ \, 1 & 4861 &-69+10i &\ \ \, 3 & 8317 & 91+6i &-1  \\
 1741 &-29+30i &\ \ \, 3 & 4877 &-61+34i &\ \ \, 5 & 8429 &-77+50i &\ \ \, 1  \\
 1789 &-5+42i &\ \ \, 1 & 4909 & 3+70i &\ \ \, 1 & 8461 & 19+90i &\ \ \, 1  \\
 1901 & 35+26i &-1 & 4957 &-69+14i &-1 & 8573 & 43+82i &\ \ \, 3  \\
 1933 &-13+42i &\ \ \, 1 & 4973 & 67+22i &-1 & 8669 &-85+38i &\ \ \, 1  \\
 1949 & 43+10i &\ \ \, 1 & 5021 & 11+70i &-1 & 8861 &-5+94i &\ \ \, 5  \\
 1997 &-29+34i &\ \ \, 1 & 5101 & 51+50i &-1 & 8893 &-53+78i &-3  \\
 2029 &-45+2i &\ \ \, 1 & 5197 &-29+66i &-1 & 8941 &-29+90i &-3  \\
 2141 &-5+46i &-3 & 5261 & 19+70i &\ \ \, 3 & 9133 &-93+22i &\ \ \, 5  \\
 2221 &-45+14i &-5 & 5309 &-53+50i &-7 & 9181 & 91+30i &-1  \\
 2237 & 11+46i &-1 & 5437 &-69+26i &\ \ \, 1 & 9277 &-21+94i &-5  \\
 2269 &-37+30i &\ \ \, 1 & 5501 &-5+74i &\ \ \, 1 & 9293 &-77+58i &-5  \\
 2333 & 43+22i &-1 & 5581 & 35+66i &\ \ \, 1 & 9341 &-85+46i &\ \ \, 5  \\
 2381 & 35+34i &\ \ \, 1 & 5693 & 43+62i &\ \ \, 5 & 9421 &-45+86i &-5  \\
\end{array} \]

\vfill

\vfill

\newpage

\small
\noindent
{} \hfill
{\bf Appendix II.\ \ Table of the coefficients of elliptic Gauss sums (2)} 
\hfill {}

\vfill

\scriptsize
\[ p=\boldsymbol{\pi} \overline{\boldsymbol{\pi}}\meqv{5}{16},\quad 
\G_{\bpi}(\chi_{\bpi},\,\Ze)=\alpha_{\bpi}\,\widetilde{\bpi}^3  \]
\[ \begin{array}{rrl|crl|crl}
p  & \boldsymbol{\pi}\quad  & \quad
 \alpha_{\boldsymbol{\pi}} \quad & 
p & \boldsymbol{\pi}\quad \ \ & \quad
 \alpha_{\boldsymbol{\pi}} \quad & 
p & \boldsymbol{\pi}\quad \ \ & \quad
 \alpha_{\boldsymbol{\pi}} \quad \\
\hline
 37 &-1+6i &\ \ \, 1 \cdot(i)& 2693 & 47+22i &-1 \cdot(i)& 6037 &-41+66i &-1 \cdot( 1 ) \\
 53 & 7+2i &\ \ \, 1 \cdot( 1 )& 2741 &-25+46i &-3 \cdot(i)& 6053 & 47+62i &\ \ \, 1 \cdot(-i) \\
 101 &-1+10i &\ \ \, 1 \cdot( 1 )& 2789 &-17+50i &-1 \cdot(-1 )& 6101 &-25+74i &-1 \cdot(-1 ) \\
 149 & 7+10i &\ \ \, 1 \cdot(-1 )& 2837 &-41+34i &\ \ \, 5 \cdot( 1 )& 6133 & 7+78i &-1 \cdot(i) \\
 181 &-9+10i &-1 \cdot(-1 )& 2917 &-1+54i &-3 \cdot(i)& 6197 & 71+34i &-1 \cdot( 1 ) \\
 197 &-1+14i &-1 \cdot(-i)& 3061 & 55+6i &-1 \cdot(-i)& 6229 &-73+30i &\ \ \, 5 \cdot(i) \\
 229 & 15+2i &-1 \cdot(-1 )& 3109 & 47+30i &\ \ \, 1 \cdot(-i)& 6277 & 79+6i &-5 \cdot(i) \\
 277 &-9+14i &-1 \cdot(i)& 3221 & 55+14i &-1 \cdot(i)& 6373 &-17+78i &-1 \cdot(-i) \\
 293 &-17+2i &\ \ \, 1 \cdot(-1 )& 3253 &-57+2i &-3 \cdot( 1 )& 6389 & 55+58i &\ \ \, 1 \cdot(-1 ) \\
 373 & 7+18i &-1 \cdot( 1 )& 3301 &-49+30i &-1 \cdot(-i)& 6421 & 39+70i &-1 \cdot(-i) \\
 389 &-17+10i &\ \ \, 1 \cdot( 1 )& 3413 & 7+58i &-1 \cdot(-1 )& 6469 & 63+50i &-1 \cdot(-1 ) \\
 421 & 15+14i &\ \ \, 1 \cdot(-i)& 3461 & 31+50i &-1 \cdot(-1 )& 6581 &-41+70i &-1 \cdot(-i) \\
 613 &-17+18i &-1 \cdot(-1 )& 3541 &-25+54i &\ \ \, 1 \cdot(-i)& 6661 &-81+10i &-3 \cdot( 1 ) \\
 661 &-25+6i &-1 \cdot(-i)& 3557 &-49+34i &\ \ \, 1 \cdot(-1 )& 6709 &-25+78i &-1 \cdot(i) \\
 677 &-1+26i &\ \ \, 3 \cdot( 1 )& 3637 & 39+46i &\ \ \, 5 \cdot(i)& 6869 & 55+62i &-3 \cdot(i) \\
 709 & 15+22i &\ \ \, 1 \cdot(i)& 3701 & 55+26i &-1 \cdot(-1 )& 6917 & 79+26i &\ \ \, 3 \cdot( 1 ) \\
 757 &-9+26i &\ \ \, 1 \cdot(-1 )& 3733 &-57+22i &\ \ \, 1 \cdot(-i)& 6949 & 15+82i &\ \ \, 1 \cdot(-1 ) \\
 773 &-17+22i &\ \ \, 1 \cdot(i)& 3797 &-41+46i &\ \ \, 5 \cdot(i)& 6997 & 39+74i &-1 \cdot(-1 ) \\
 821 &-25+14i &\ \ \, 1 \cdot(i)& 3877 & 31+54i &\ \ \, 1 \cdot(i)& 7013 &-17+82i &-1 \cdot(-1 ) \\
 853 & 23+18i &\ \ \, 1 \cdot( 1 )& 3989 &-25+58i &\ \ \, 1 \cdot(-1 )& 7109 & 47+70i &\ \ \, 1 \cdot(i) \\
 997 & 31+6i &\ \ \, 1 \cdot(i)& 4021 & 39+50i &\ \ \, 3 \cdot( 1 )& 7237 &-81+26i &-7 \cdot( 1 ) \\
 1013 & 23+22i &\ \ \, 1 \cdot(-i)& 4133 &-17+62i &\ \ \, 1 \cdot(-i)& 7253 & 23+82i &-5 \cdot( 1 ) \\
 1061 & 31+10i &\ \ \, 1 \cdot( 1 )& 4229 &-65+2i &-1 \cdot(-1 )& 7333 & 63+58i &\ \ \, 3 \cdot( 1 ) \\
 1093 &-33+2i &-1 \cdot(-1 )& 4261 &-65+6i &\ \ \, 3 \cdot(i)& 7349 &-25+82i &-1 \cdot( 1 ) \\
 1109 &-25+22i &\ \ \, 1 \cdot(-i)& 4357 &-1+66i &-1 \cdot(-1 )& 7477 &-9+86i &-1 \cdot(-i) \\
 1237 &-9+34i &\ \ \, 1 \cdot( 1 )& 4373 & 23+62i &-1 \cdot(i)& 7541 & 71+50i &\ \ \, 1 \cdot( 1 ) \\
 1301 &-25+26i &\ \ \, 1 \cdot(-1 )& 4421 &-65+14i &-1 \cdot(-i)& 7573 & 87+2i &\ \ \, 5 \cdot( 1 ) \\
 1381 & 15+34i &-1 \cdot(-1 )& 4517 &-49+46i &\ \ \, 1 \cdot(-i)& 7589 &-65+58i &\ \ \, 1 \cdot( 1 ) \\
 1429 & 23+30i &-1 \cdot(i)& 4549 &-65+18i &-1 \cdot(-1 )& 7621 & 15+86i &\ \ \, 3 \cdot(i) \\
 1493 & 7+38i &-1 \cdot(-i)& 4597 &-41+54i &\ \ \, 1 \cdot(-i)& 7669 & 87+10i &\ \ \, 3 \cdot(-1 ) \\
 1621 & 39+10i &-1 \cdot(-1 )& 4789 & 55+42i &\ \ \, 1 \cdot(-1 )& 7717 &-81+34i &\ \ \, 1 \cdot(-1 ) \\
 1637 & 31+26i &-1 \cdot( 1 )& 4933 &-33+62i &-1 \cdot(-i)& 7829 &-73+50i &\ \ \, 5 \cdot( 1 ) \\
 1669 & 15+38i &-3 \cdot(i)& 5077 & 71+6i &-1 \cdot(-i)& 7877 &-49+74i &-1 \cdot( 1 ) \\
 1733 &-17+38i &-3 \cdot(i)& 5189 &-17+70i &-1 \cdot(i)& 8053 & 87+22i &\ \ \, 1 \cdot(-i) \\
 1861 & 31+30i &-1 \cdot(-i)& 5237 & 71+14i &-1 \cdot(i)& 8069 &-65+62i &-1 \cdot(-i) \\
 1877 &-41+14i &\ \ \, 3 \cdot(i)& 5333 &-73+2i &\ \ \, 1 \cdot( 1 )& 8101 &-1+90i &\ \ \, 1 \cdot( 1 ) \\
 1973 & 23+38i &-1 \cdot(-i)& 5381 &-65+34i &-1 \cdot(-1 )& 8117 &-89+14i &\ \ \, 1 \cdot(i) \\
 2053 &-17+42i &-1 \cdot( 1 )& 5413 & 63+38i &\ \ \, 1 \cdot(i)& 8293 & 47+78i &-1 \cdot(-i) \\
 2069 &-25+38i &\ \ \, 1 \cdot(-i)& 5477 &-1+74i &-5 \cdot( 1 )& 8389 &-17+90i &\ \ \, 1 \cdot( 1 ) \\
 2213 & 47+2i &\ \ \, 1 \cdot(-1 )& 5557 &-9+74i &\ \ \, 1 \cdot(-1 )& 8501 & 55+74i &-1 \cdot(-1 ) \\
 2293 & 23+42i &-1 \cdot(-1 )& 5573 & 47+58i &\ \ \, 1 \cdot( 1 )& 8581 &-65+66i &-1 \cdot(-1 ) \\
 2309 & 47+10i &\ \ \, 5 \cdot( 1 )& 5653 &-73+18i &\ \ \, 3 \cdot( 1 )& 8597 &-89+26i &\ \ \, 1 \cdot(-1 ) \\
 2341 & 15+46i &-1 \cdot(-i)& 5669 &-65+38i &\ \ \, 1 \cdot(i)& 8629 & 23+90i &-1 \cdot(-1 ) \\
 2357 &-41+26i &-1 \cdot(-1 )& 5701 & 15+74i &\ \ \, 1 \cdot( 1 )& 8677 &-81+46i &\ \ \, 3 \cdot(-i) \\
 2389 &-25+42i &\ \ \, 1 \cdot(-1 )& 5717 & 71+26i &-1 \cdot(-1 )& 8693 &-73+58i &-1 \cdot(-1 ) \\
 2437 &-49+6i &-1 \cdot(i)& 5749 &-57+50i &-1 \cdot( 1 )& 8741 & 79+50i &\ \ \, 1 \cdot(-1 ) \\
 2549 & 7+50i &-5 \cdot( 1 )& 5813 &-73+22i &\ \ \, 1 \cdot(-i)& 8821 &-89+30i &\ \ \, 1 \cdot(i) \\
 2677 & 39+34i &-3 \cdot( 1 )& 5861 & 31+70i &\ \ \, 5 \cdot(i)& 8837 &-1+94i &\ \ \, 1 \cdot(-i) \\
\end{array} \]

\vfill

\vfill

\newpage

\small
\noindent
{} \hfill
{\bf Appendix II.\ \ Table of the coefficients of elliptic Gauss sums (3)}
\hfill {}

\vfill

\scriptsize
\[ p=\boldsymbol{\pi} \overline{\boldsymbol{\pi}}\meqv{1}{16},\quad 
\G_{\bpi}(\chi_{\bpi},\,\psi)=\alpha_{\bpi}\,\widetilde{\bpi}^3 \]
\[ \begin{array}{rrlc|crlc}
p\ \   & \boldsymbol{\pi}\quad & \quad
 \alpha_{\boldsymbol{\pi}} & \chi_{\pi}(1+i) & 
p & \boldsymbol{\pi}\quad & \quad
 \alpha_{\boldsymbol{\pi}} & \chi_{\pi}(1+i) \\
\hline
 17 & 1+4i &\ \ \, 1 \cdot i\zeta_8 & -i & 2897 &-31+44i &-1 \cdot \zeta_8 & \ \ \,i  \\
 97 & 9+4i &-1 \cdot \zeta_8 & \ \ \,i & 3041 &-55+4i &-1 \cdot \zeta_8 & \ \ \,i  \\
 113 &-7+8i &\ \ \, 0 \cdot i\sqrt{2} & \ \ \,1  & 3089 &-55+8i &\ \ \, 2 \cdot i\sqrt{2} & \ \ \,1  \\
 193 &-7+12i &-1 \cdot i\zeta_8 & -i & 3121 &-39+40i &\ \ \, 0 \cdot i\sqrt{2} & \ \ \,1  \\
 241 &-15+4i &\ \ \, 1 \cdot i\zeta_8 & -i & 3137 & 1+56i &\ \ \, 1 \cdot \sqrt{2} &-1  \\
 257 & 1+16i &\ \ \, 0 \cdot i\sqrt{2} & \ \ \,1 & 3169 &-55+12i &-1 \cdot i\zeta_8 & -i  \\
 337 & 9+16i &\ \ \, 1 \cdot \sqrt{2} &-1 & 3217 & 9+56i &\ \ \, 0 \cdot i\sqrt{2} & \ \ \,1  \\
 353 & 17+8i &\ \ \, 1 \cdot \sqrt{2} &-1 & 3313 & 57+8i &-2 \cdot i\sqrt{2} & \ \ \,1  \\
 401 & 1+20i &\ \ \, 1 \cdot i\zeta_8 & -i & 3329 & 25+52i &-1 \cdot \zeta_8 & \ \ \,i  \\
 433 & 17+12i &-1 \cdot \zeta_8 & \ \ \,i & 3361 &-15+56i &\ \ \, 1 \cdot \sqrt{2} &-1  \\
 449 &-7+20i &\ \ \, 1 \cdot \zeta_8 & \ \ \,i & 3457 &-39+44i &-5 \cdot i\zeta_8 & -i  \\
 577 & 1+24i &\ \ \, 1 \cdot \sqrt{2} &-1 & 3617 & 41+44i &-5 \cdot i\zeta_8 & -i  \\
 593 &-23+8i &\ \ \, 0 \cdot i\sqrt{2} & \ \ \,1 & 3697 & 49+36i &\ \ \, 1 \cdot i\zeta_8 & -i  \\
 641 & 25+4i &\ \ \, 1 \cdot \zeta_8 & \ \ \,i & 3761 & 25+56i &\ \ \, 0 \cdot i\sqrt{2} & \ \ \,1  \\
 673 &-23+12i &\ \ \, 1 \cdot i\zeta_8 & -i & 3793 & 33+52i &-1 \cdot i\zeta_8 & -i  \\
 769 & 25+12i &\ \ \, 1 \cdot i\zeta_8 & -i & 3889 & 17+60i &\ \ \, 1 \cdot \zeta_8 & \ \ \,i  \\
 881 & 25+16i &\ \ \, 1 \cdot \sqrt{2} &-1 & 4001 & 49+40i &-1 \cdot \sqrt{2} &-1  \\
 929 &-23+20i &\ \ \, 1 \cdot \zeta_8 & \ \ \,i & 4049 &-55+32i &\ \ \, 1 \cdot \sqrt{2} &-1  \\
 977 &-31+4i &\ \ \, 1 \cdot i\zeta_8 & -i & 4129 &-23+60i &-1 \cdot i\zeta_8 & -i  \\
 1009 &-15+28i &-1 \cdot \zeta_8 & \ \ \,i & 4177 & 9+64i &-1 \cdot \sqrt{2} &-1  \\
 1153 & 33+8i &\ \ \, 1 \cdot \sqrt{2} &-1 & 4241 & 65+4i &\ \ \, 1 \cdot i\zeta_8 & -i  \\
 1201 & 25+24i &\ \ \, 0 \cdot i\sqrt{2} & \ \ \,1 & 4273 & 57+32i &\ \ \, 1 \cdot \sqrt{2} &-1  \\
 1217 &-31+16i &\ \ \, 0 \cdot i\sqrt{2} & \ \ \,1 & 4289 & 65+8i &\ \ \, 1 \cdot \sqrt{2} &-1  \\
 1249 &-15+32i &\ \ \, 2 \cdot i\sqrt{2} & \ \ \,1 & 4337 & 49+44i &\ \ \, 1 \cdot \zeta_8 & \ \ \,i  \\
 1297 & 1+36i &\ \ \, 3 \cdot i\zeta_8 & -i & 4481 & 65+16i &\ \ \, 0 \cdot i\sqrt{2} & \ \ \,1  \\
1361 &-31+20i &\ \ \, 1 \cdot i\zeta_8 & -i & 4513 &-47+48i &\ \ \, 0 \cdot i\sqrt{2} & \ \ \,1  \\
 1409 & 25+28i &-1 \cdot i\zeta_8 & -i & 4561 &-31+60i &-1 \cdot \zeta_8 & \ \ \,i  \\
 1489 & 33+20i &\ \ \, 5 \cdot i\zeta_8 & -i & 4657 &-39+56i &-2 \cdot i\sqrt{2} & \ \ \,1  \\
 1553 &-23+32i &\ \ \, 1 \cdot \sqrt{2} &-1 & 4673 &-7+68i &-1 \cdot \zeta_8 & \ \ \,i  \\
 1601 & 1+40i &\ \ \, 1 \cdot \sqrt{2} &-1 & 4721 & 25+64i &-1 \cdot \sqrt{2} &-1  \\
 1697 & 41+4i &-3 \cdot \zeta_8 & \ \ \,i & 4801 & 65+24i &\ \ \, 1 \cdot \sqrt{2} &-1  \\
 1777 &-39+16i &-1 \cdot \sqrt{2} &-1 & 4817 & 41+56i &\ \ \, 0 \cdot i\sqrt{2} & \ \ \,1  \\
 1873 & 33+28i &-1 \cdot \zeta_8 & \ \ \,i & 4993 &-63+32i &\ \ \, 0 \cdot i\sqrt{2} & \ \ \,1  \\
 1889 & 17+40i &\ \ \, 1 \cdot \sqrt{2} &-1 & 5009 & 65+28i &\ \ \, 1 \cdot \zeta_8 & \ \ \,i  \\
 2017 & 9+44i &-1 \cdot i\zeta_8 & -i & 5153 &-23+68i &\ \ \, 1 \cdot \zeta_8 & \ \ \,i  \\
 2081 & 41+20i &-1 \cdot \zeta_8 & \ \ \,i & 5233 &-7+72i &\ \ \, 0 \cdot i\sqrt{2} & \ \ \,1  \\
 2113 & 33+32i &-2 \cdot i\sqrt{2} & \ \ \,1 & 5281 & 41+60i &\ \ \, 1 \cdot i\zeta_8 & -i  \\
 2129 &-23+40i &\ \ \, 0 \cdot i\sqrt{2} & \ \ \,1 & 5297 &-71+16i &\ \ \, 1 \cdot \sqrt{2} &-1  \\
 2161 &-15+44i &-1 \cdot \zeta_8 & \ \ \,i & 5393 & 73+8i &-2 \cdot i\sqrt{2} & \ \ \,1  \\
 2273 &-47+8i &\ \ \, 1 \cdot \sqrt{2} &-1 & 5441 &-71+20i &\ \ \, 1 \cdot \zeta_8 & \ \ \,i  \\
 2417 & 49+4i &\ \ \, 1 \cdot i\zeta_8 & -i & 5521 & 65+36i &\ \ \, 1 \cdot i\zeta_8 & -i  \\
 2593 & 17+48i &\ \ \, 0 \cdot i\sqrt{2} & \ \ \,1 & 5569 &-63+40i &-1 \cdot \sqrt{2} &-1  \\
 2609 &-47+20i &\ \ \, 5 \cdot i\zeta_8 & -i & 5857 & 9+76i &\ \ \, 1 \cdot i\zeta_8 & -i  \\
 2657 & 49+16i &\ \ \, 0 \cdot i\sqrt{2} & \ \ \,1 & 5953 & 57+52i &\ \ \, 1 \cdot \zeta_8 & \ \ \,i  \\
 2689 & 33+40i &-3 \cdot \sqrt{2} &-1 & 6113 & 73+28i &-1 \cdot i\zeta_8 & -i  \\
 2753 &-7+52i &-1 \cdot \zeta_8 & \ \ \,i & 6257 &-79+4i &-1 \cdot i\zeta_8 & -i  \\
 2801 & 49+20i &-1 \cdot i\zeta_8 & -i & 6337 &-71+36i &\ \ \, 1 \cdot \zeta_8 & \ \ \,i  \\
 2833 &-23+48i &\ \ \, 1 \cdot \sqrt{2} &-1 & 6353 & 73+32i &\ \ \, 1 \cdot \sqrt{2} &-1  \\
\end{array} \]

\vfill

\vfill

\newpage

\small
\noindent
{} \hfill
{\bf Appendix II.\ \ Table of the coefficients of elliptic Gauss sums (4)}
\hfill {}

\vfill

\scriptsize
\[ p=\boldsymbol{\pi} \overline{\boldsymbol{\pi}}\meqv{9}{16},\quad 
\G_{\bpi}(\chi_{\bpi},\,\psi)=\alpha_{\bpi}\,\widetilde{\bpi}^3 \]
\[ \begin{array}{rrlc|crlc}
p\ \   & \boldsymbol{\pi}\quad & \quad
 \alpha_{\boldsymbol{\pi}} & \chi_{\pi}(1+i) & 
p & \boldsymbol{\pi}\quad & \quad
 \alpha_{\boldsymbol{\pi}} & \chi_{\pi}(1+i) \\
\hline
 41 & 5+4i &-1 \cdot i\zeta_8 & \ \ \,1& 2777 & 29+44i &-1 \cdot i\zeta_8 & \ \ \,1 \\
 73 &-3+8i &\ \ \, 1 \cdot i\sqrt{2} & \ \ \,i & 2857 &-51+16i &\ \ \, 0 \cdot \sqrt{2} & -i  \\
 89 & 5+8i &\ \ \, 0 \cdot \sqrt{2} & -i & 2953 & 53+12i &-1 \cdot \zeta_8 &-1  \\
 137 &-11+4i &-1 \cdot i\zeta_8 & \ \ \,1& 2969 & 37+40i &\ \ \, 0 \cdot \sqrt{2} & -i  \\
 233 & 13+8i &\ \ \, 1 \cdot i\sqrt{2} & \ \ \,i & 3001 &-51+20i &\ \ \, 1 \cdot \zeta_8 &-1  \\
 281 & 5+16i &\ \ \, 1 \cdot i\sqrt{2} & \ \ \,i & 3049 & 45+32i &\ \ \, 0 \cdot \sqrt{2} & -i  \\
 313 & 13+12i &\ \ \, 1 \cdot i\zeta_8 & \ \ \,1& 3209 & 53+20i &-3 \cdot i\zeta_8 & \ \ \,1 \\
 409 &-3+20i &-1 \cdot \zeta_8 &-1 & 3257 &-11+56i &\ \ \, 0 \cdot \sqrt{2} & -i  \\
 457 & 21+4i &\ \ \, 1 \cdot i\zeta_8 & \ \ \,1& 3433 &-27+52i &-1 \cdot i\zeta_8 & \ \ \,1 \\
 521 &-11+20i &-1 \cdot i\zeta_8 & \ \ \,1& 3449 &-43+40i &\ \ \, 0 \cdot \sqrt{2} & -i  \\
 569 & 13+20i &-3 \cdot \zeta_8 &-1 & 3529 &-35+48i &\ \ \, 0 \cdot \sqrt{2} & -i  \\
 601 & 5+24i &\ \ \, 0 \cdot \sqrt{2} & -i & 3593 & 53+28i &\ \ \, 1 \cdot \zeta_8 &-1  \\
 617 &-19+16i &\ \ \, 2 \cdot \sqrt{2} & -i & 3673 & 37+48i &-1 \cdot i\sqrt{2} & \ \ \,i  \\
 761 &-19+20i &-1 \cdot \zeta_8 &-1 & 3769 & 13+60i &\ \ \, 1 \cdot i\zeta_8 & \ \ \,1 \\
 809 & 5+28i &-1 \cdot \zeta_8 &-1 & 3833 & 53+32i &-1 \cdot i\sqrt{2} & \ \ \,i  \\
 857 & 29+4i &\ \ \, 1 \cdot \zeta_8 &-1 & 3881 &-59+20i &-1 \cdot i\zeta_8 & \ \ \,1 \\
 937 &-19+24i &\ \ \, 1 \cdot i\sqrt{2} & \ \ \,i & 3929 &-35+52i &\ \ \, 1 \cdot \zeta_8 &-1  \\
 953 & 13+28i &-1 \cdot i\zeta_8 & \ \ \,1& 4057 &-59+24i &-2 \cdot \sqrt{2} & -i  \\
 1033 &-3+32i &\ \ \, 2 \cdot \sqrt{2} & -i & 4073 & 37+52i &-3 \cdot i\zeta_8 & \ \ \,1 \\
 1049 & 5+32i &\ \ \, 1 \cdot i\sqrt{2} & \ \ \,i & 4153 &-43+48i &-1 \cdot i\sqrt{2} & \ \ \,i  \\
 1097 & 29+16i &\ \ \, 0 \cdot \sqrt{2} & -i & 4201 &-51+40i &\ \ \, 3 \cdot i\sqrt{2} & \ \ \,i  \\
 1129 &-27+20i &\ \ \, 3 \cdot i\zeta_8 & \ \ \,1& 4217 &-11+64i &-1 \cdot i\sqrt{2} & \ \ \,i  \\
 1193 & 13+32i &\ \ \, 0 \cdot \sqrt{2} & -i & 4297 & 61+24i &\ \ \, 1 \cdot i\sqrt{2} & \ \ \,i  \\
 1289 &-35+8i &\ \ \, 1 \cdot i\sqrt{2} & \ \ \,i & 4409 & 53+40i &-2 \cdot \sqrt{2} & -i  \\
 1321 & 5+36i &-1 \cdot i\zeta_8 & \ \ \,1& 4441 & 29+60i &\ \ \, 1 \cdot i\zeta_8 & \ \ \,1 \\
 1433 & 37+8i &\ \ \, 0 \cdot \sqrt{2} & -i & 4457 &-19+64i &\ \ \, 0 \cdot \sqrt{2} & -i  \\
 1481 &-35+16i &\ \ \, 0 \cdot \sqrt{2} & -i & 4649 & 5+68i &\ \ \, 1 \cdot i\zeta_8 & \ \ \,1 \\
 1609 &-3+40i &\ \ \, 1 \cdot i\sqrt{2} & \ \ \,i & 4729 & 45+52i &-1 \cdot \zeta_8 &-1  \\
 1657 &-19+36i &-3 \cdot \zeta_8 &-1 & 4793 & 13+68i &-1 \cdot \zeta_8 &-1  \\
 1721 &-11+40i &\ \ \, 0 \cdot \sqrt{2} & -i & 4889 &-67+20i &-5 \cdot \zeta_8 &-1  \\
 1753 &-27+32i &-1 \cdot i\sqrt{2} & \ \ \,i & 4937 & 29+64i &\ \ \, 0 \cdot \sqrt{2} & -i  \\
 1801 &-35+24i &\ \ \, 1 \cdot i\sqrt{2} & \ \ \,i & 4969 & 37+60i &-1 \cdot \zeta_8 &-1  \\
 1913 &-43+8i &\ \ \, 2 \cdot \sqrt{2} & -i & 5081 &-59+40i &\ \ \, 0 \cdot \sqrt{2} & -i  \\
 1993 &-43+12i &-1 \cdot \zeta_8 &-1 & 5113 & 53+48i &-1 \cdot i\sqrt{2} & \ \ \,i  \\
 2089 & 45+8i &-1 \cdot i\sqrt{2} & \ \ \,i & 5209 & 5+72i &-2 \cdot \sqrt{2} & -i  \\
 2137 & 29+36i &\ \ \, 1 \cdot \zeta_8 &-1 & 5273 &-67+28i &\ \ \, 1 \cdot i\zeta_8 & \ \ \,1 \\
 2153 & 37+28i &-1 \cdot \zeta_8 &-1 & 5417 &-59+44i &\ \ \, 1 \cdot \zeta_8 &-1  \\
 2281 & 45+16i &\ \ \, 2 \cdot \sqrt{2} & -i & 5449 &-43+60i &-1 \cdot \zeta_8 &-1  \\
 2297 &-19+44i &-1 \cdot i\zeta_8 & \ \ \,1& 5641 &-75+4i &-1 \cdot i\zeta_8 & \ \ \,1 \\
 2377 & 21+44i &-5 \cdot \zeta_8 &-1 & 5657 & 61+44i &-3 \cdot i\zeta_8 & \ \ \,1 \\
 2393 & 37+32i &-1 \cdot i\sqrt{2} & \ \ \,i & 5689 &-75+8i &\ \ \, 0 \cdot \sqrt{2} & -i  \\
 2441 & 29+40i &\ \ \, 1 \cdot i\sqrt{2} & \ \ \,i & 5737 &-51+56i &-1 \cdot i\sqrt{2} & \ \ \,i  \\
 2473 & 13+48i &\ \ \, 0 \cdot \sqrt{2} & -i & 5801 & 5+76i &\ \ \, 1 \cdot \zeta_8 &-1  \\
 2521 &-35+36i &\ \ \, 1 \cdot \zeta_8 &-1 & 5849 &-35+68i &-1 \cdot \zeta_8 &-1  \\
 2617 &-51+4i &-1 \cdot \zeta_8 &-1 & 5881 &-75+16i &-1 \cdot i\sqrt{2} & \ \ \,i  \\
 2633 &-43+28i &\ \ \, 1 \cdot \zeta_8 &-1 & 5897 &-11+76i &\ \ \, 1 \cdot \zeta_8 &-1  \\
 2713 &-3+52i &-1 \cdot \zeta_8 &-1 & 6073 & 77+12i &\ \ \, 1 \cdot i\zeta_8 & \ \ \,1 \\
 2729 & 5+52i &-1 \cdot i\zeta_8 & \ \ \,1& 6089 &-67+40i &-1 \cdot i\sqrt{2} & \ \ \,i  \\
\end{array} \]

\vfill

\vfill

\end{document}